\documentclass[leqno,11pt]{amsart}
\usepackage{amssymb, amsmath}

\usepackage{xcolor,ulem,textcomp,MnSymbol}
\usepackage{url}
\usepackage[mathscr]{euscript}
\usepackage{enumerate}
\usepackage{graphicx}
\usepackage{cancel}

\theoremstyle{definition}
\theoremstyle{plain}
\newtheorem{thm}{Theorem}[section]

 \newtheorem{Thm}{Theorem}[section]
 \newtheorem{Rmk}[thm]{Remark}

  \definecolor{darkgreen}{rgb}{0.0, 0.5, 0.0}

\def\bbE {\mathbb {E}}

\def\bbP {\mathbb{P}}
\def\bbR {\mathbb{R}}
\def\bbS {\mathbb{S}}
\def\bbT {\mathbb{T}}

\def\cA {\mathcal{A}}
\def\cB {\mathcal{B}}
\def\cC {\mathcal{C}}

\def\cF {\mathcal{F}}
\def\cH {\mathcal{H}}
\def\cG {\mathcal{G}}
\def\cG {\mathcal{G}}
\def\cI  {\mathcal{I}}
\def\cL {\mathcal{L}}

\def\cN {\mathcal{N}}

\def\cP {\mathcal{P}}
\def\cR {\mathcal{R}}

\def\cZ {\mathcal{Z}}

\def\gp {{\varphi}}

\def\eps {{\varepsilon}}
\def\e {{\varepsilon}}

\def\indc {\mathbf{1}}

\def\la {\langle}
\def\ra {\rangle}

\def\d {{\partial}}

\newcommand{\Cov}{\operatorname{Cov}}

\newcommand{\Var}{\operatorname{Var}}

\newcommand{\ba}{\begin{aligned}}
\newcommand{\ea}{\end{aligned}}

\newcommand{\be}{\begin{equation}}
\newcommand{\ee}{\end{equation}}

\numberwithin{equation}{section}

\begin{document}


\title{Dynamics of dilute gases~: a~statistical approach}

%
%

\author[T. Bodineau]{Thierry Bodineau}
\address[T. Bodineau]%
{CMAP, CNRS, Ecole Polytechnique, I.P. Paris
\\
Route de Saclay, 91128 Palaiseau Cedex, FRANCE}
\email{thierry.bodineau@polytechnique.edu}

\author[I. Gallagher]{Isabelle Gallagher}
\address[I. Gallagher]%
{DMA, \'Ecole normale sup\'erieure, CNRS, PSL Research University\\
45 rue d'Ulm, 75005 Paris, FRANCE \\
and Universit\'e de Paris}
\email{gallagher@math.ens.fr}

\author[L. Saint-Raymond]{Laure Saint-Raymond}
\address[L. Saint-Raymond]
{IHES, Bures sur Yvette, FRANCE
}
\email{Laure.Saint-Raymond@ens-lyon.fr}

\author[S. Simonella]{Sergio Simonella}
\address[S. Simonella]
{UMPA UMR 5669 du CNRS, ENS de Lyon, Universit\'e de Lyon\\
46  all\'ee d'Italie, 69007 Lyon,
FRANCE
}
\email{sergio.simonella@ens-lyon.fr}



\begin{abstract}
The evolution of a  gas can be  described  by different models depending on the observation scale.
A natural question, raised by Hilbert in his sixth problem,  is whether these models provide consistent predictions.
 In particular,  for rarefied gases, it is expected  that continuum laws of kinetic theory can be obtained 
 directly from molecular dynamics governed by the fundamental principles of mechanics.

In the case  of hard sphere gases, Lanford \cite{Lanford} showed that the Boltzmann equation 
emerges as the law of large numbers  in the low density limit, at least for very short 
times. The goal of this survey is to present  recent progress in the understanding of this limiting process, 
providing a  complete  statistical description.

\end{abstract}

\maketitle

\section{AIM\,:\,PROVIDING A STATISTICAL PICTURE OF DILUTE GAS DYNAMICS}

\subsection{A very simple physical model}

Even though at the time Boltzmann published his famous paper \cite{boltzmann}, the atomistic theory was still dismissed by some scientists, it is now well established that matter is composed of atoms, 
which are the elementary constituents of all solid, liquid and gaseous substances. The particularity of dilute gases is that their atoms are very weakly bound and almost independent. In other words, there are very few constraints on their geometric arrangement because their volume is negligible compared to the total volume occupied by the gas.

\begin{figure}[h] 
\centering
\includegraphics[width=5in]{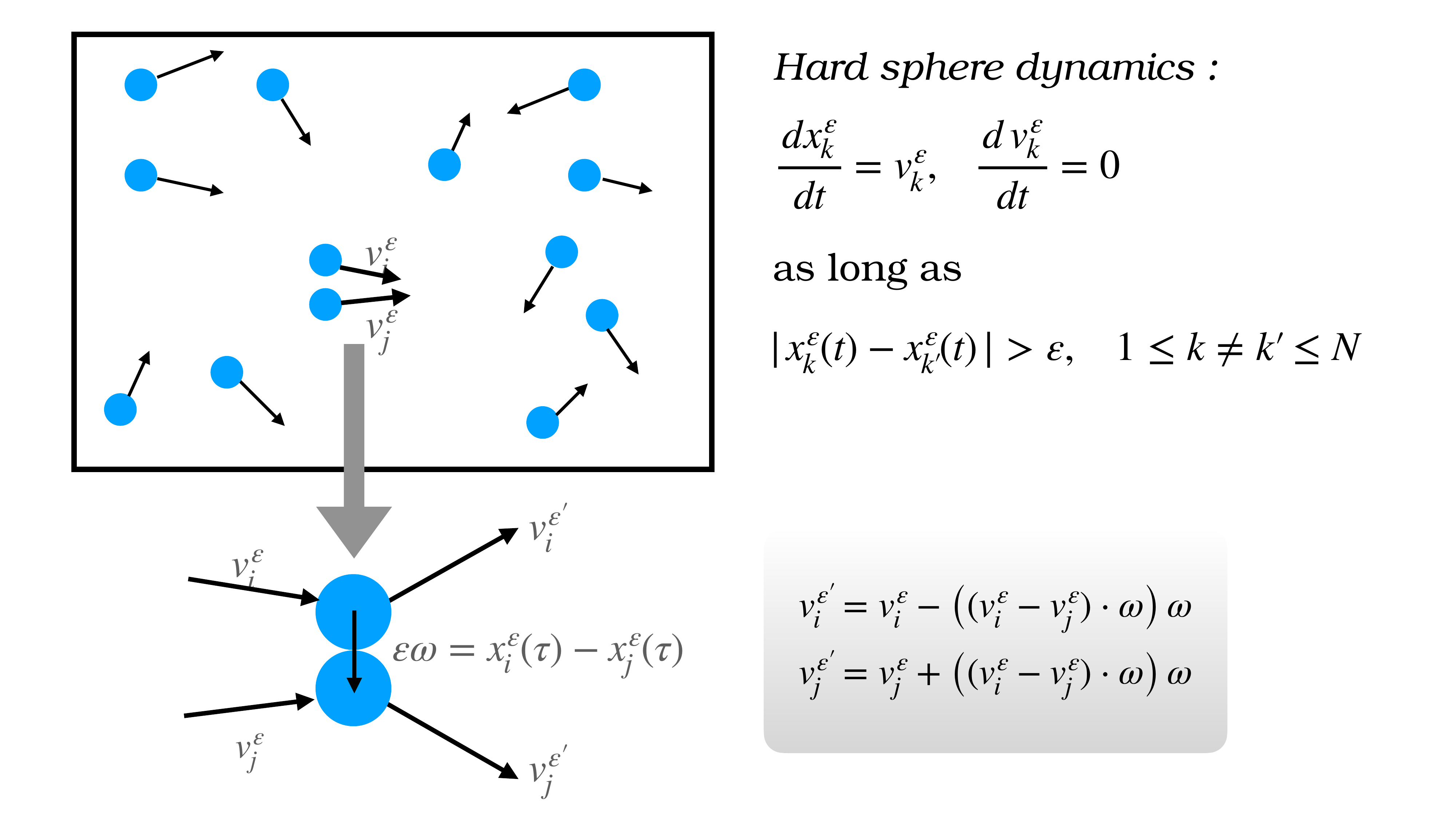} 
\caption{
\color{black}
At time $t$, the hard-sphere system is described  by the positions 
$\big( x^\varepsilon_k (t) \big)_{k \leq N}$ and the velocities 
$\big( v^\varepsilon_k (t) \big)_{k \leq N}$ of the $N$ particles. Particles move in straight lines  and  when two particles touch each other at distance $\varepsilon > 0$ (the diameter of the spheres), they are scattered according to   elastic reflection laws. 
The scattering rules, mapping the precollisional velocities $(v^\eps_i, v^\eps_j)$ to the postcollisional velocities 
$({v^\eps_i}' , {v^\eps_j}' )$, are determined in terms of the relative position
$\omega = \big( x^\varepsilon_i (\tau) - x^\varepsilon_i (\tau) \big) / \varepsilon$ of the particles at the collision time~$\tau$.
The collisions preserve the total momentum $v^\eps_i+v^\eps_j = {v^\eps_i}'+ {v^\eps_j}'$ 
and the kinetic energy $\frac{1}{2} (|v^\eps_i|^2+|v^\eps_j|^2) = \frac{1}{2} (| {v^\eps_i}' |^2+| {v^\eps_j}'|^2)$.
}
\label{HS-fig}
\end{figure}

If we neglect the internal structure of atoms (consisting of a nucleus and electrons) and their possible organisation into molecules, we can represent a gas as a large system of  correlated interacting particles. 
We will also neglect the effect of long range interactions and assume strong interatomic forces at very short distance. Each particle moves freely most of the time and occasionally collides with some other particle leading to an almost instantaneous scattering.
The simplest example of such a model consists in assuming that the particles are identical tiny balls 
of unit mass interacting only by contact (see Figure \ref{HS-fig}). 
We then speak of a \emph{gas of hard spheres}. All the results we will present should nevertheless extend to isotropic, compactly supported stable interaction potentials
\cite{Ruelle_livre,Pulvirenti_Saffirio_Simonella}.

\medskip

\begin{figure}[h] 
\centering
\includegraphics[width=4in]{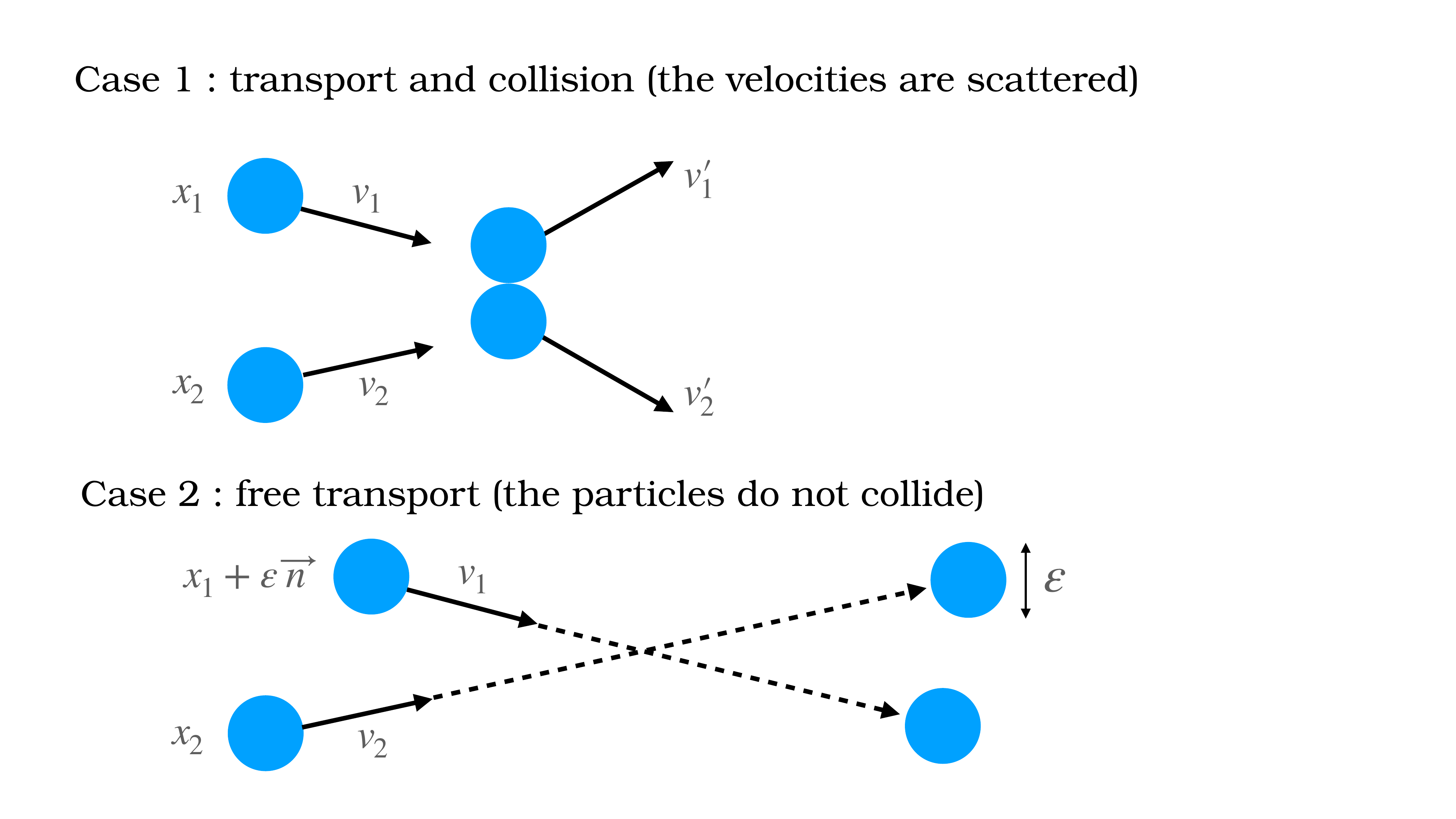} 
\caption{Particles are very small (of diameter $\eps \ll 1$) and therefore the dynamics is very sensitive to small spatial shifts. 
In the first case depicted above, two particles with initial positions $x_1,x_2$ and velocities $v_1,v_2$ collide and are scattered. In the second case, by shifting the first particle by a distance $\eps$ in the direction $\vec{n}$, the two particles 
no longer collide and they move in straight lines. Thus a perturbation of order $\eps$ of the initial conditions can lead to very different trajectories.
}
\label{fig: instability}
\end{figure}

This microscopic description of a gas is daunting because the number of particles involved is extremely large, the individual size of these particles is tiny (of diameter $\eps \ll 1$)  and therefore positions are very sensitive to small spatial shifts  
 (see Figure \ref{fig: instability}). 
In practice, this model is not efficient for making theoretical predictions, and numerical methods are often in favour of Monte Carlo simulations.
The question we would like to address here is a more fundamental one, namely  the consistency of this (simplified) atomic description with the kinetic or fluid models used in applications. This question was formalised by Hilbert at the ICM in 1900, in his sixth problem: "Boltzmann's work on the principles of mechanics suggests the problem of developing mathematically the limiting processes, there merely indicated, which lead from the atomistic view to the laws of motion of continua".

\medskip

The Boltzmann equation, mentioned by Hilbert and which we will present in more detail later, expresses that the distribution of particles evolves under the combined effect of free transport and collisions. For these two effects to be of the same order of magnitude, a simple calculation shows that, in dimension $d \geq 2$, the number of particles $N$ and their diameter size $\eps$ must satisfy the scaling relation $N\eps^{d-1} = O(1)$, the so-called \emph{Boltzmann-Grad scaling} \cite{Grad}.
 \begin{figure}[h] 
\centering
\includegraphics[width=2in]{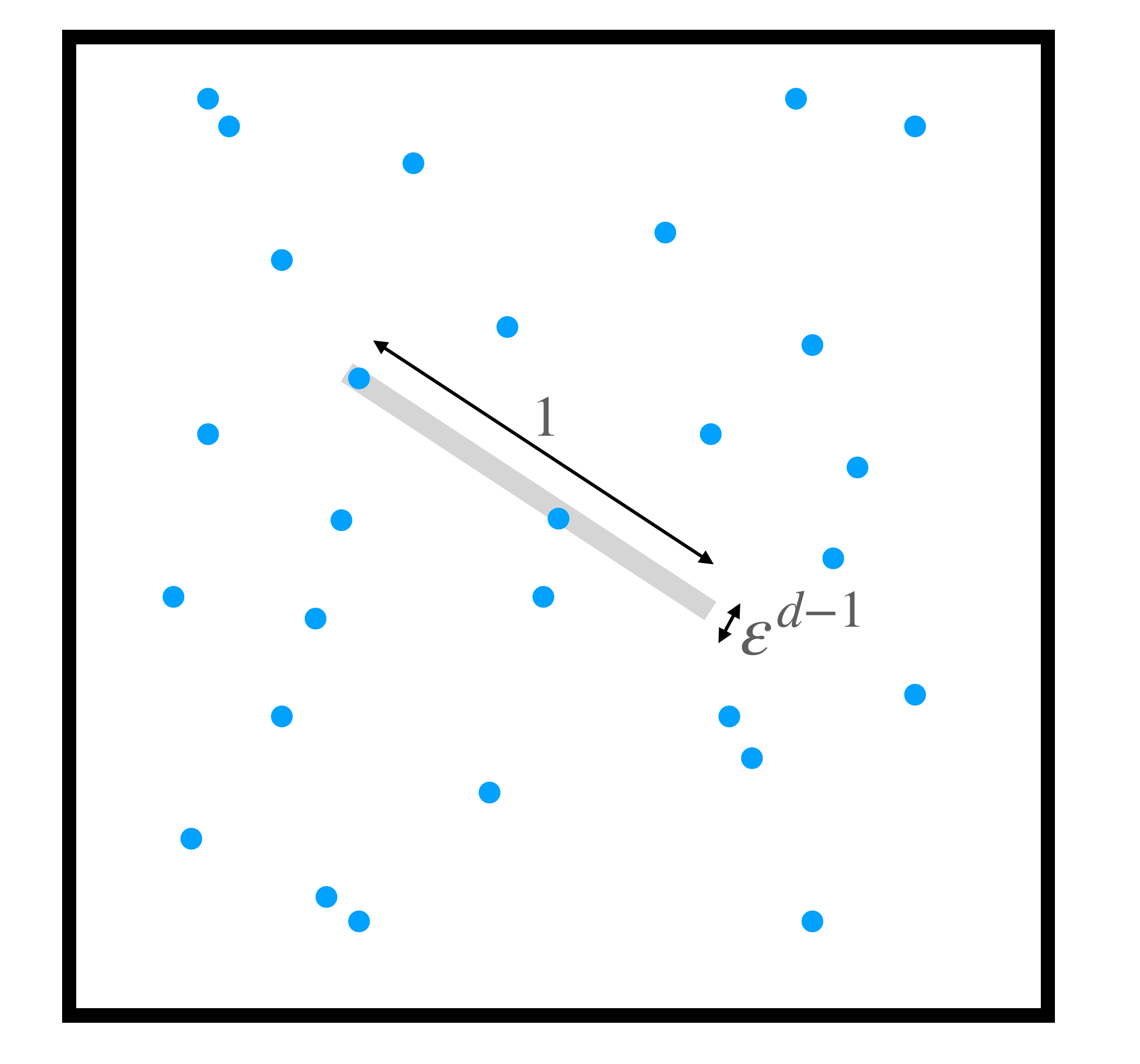} 
\caption{ 
Consider $N$ spheres of diameter $\eps$ uniformly distributed in a box.
If the mean free path is equal to $1$, then the grey tube of length 1 and section area of order $\eps^{d-1}$ represents the volume spanned by a typical particle between two collisions. 
The Boltzmann-Grad scaling~$N\eps^{d-1}= 1$ is tuned such that on average this tube intersects one particle.
}
\label{fig: scaling}
\end{figure} 
Indeed the regime described by the Boltzmann equation is such that the mean free path, namely the average distance covered by a particle travelling in straight line between two collisions, is of order 1. Thus 
a typical particle trajectory should span a tube of volume $1 \times \eps^{d-1}$  between two collisions. This means that on average this tube should intersect the position of one of the other $(N-1)$ particles (see Figure \ref{fig: scaling}). 
Note that in this regime the total volume occupied by the particles at a given time is proportional to $N\eps^d$ and therefore is negligible compared to the total volume occupied by the gas. We speak then of a \emph{dilute gas}.

\subsection{Three levels of averaging}
\label{subsec: Three levels of averaging}

As  already shown in the previous scaling argument, the equations that we want to derive  describe the behavior of ``typical particles''. We   therefore have to introduce several averaging processes, and then  to describe the average dynamics.

For a statistical description of a monoatomic gas, all particles are considered identical (same geometry, same mass, same interaction law,...). This is referred to as the exchangeability assumption. The \emph{empirical distribution} of particles is defined as
\begin{equation}
\label{eq: empirical measure}
\pi^N _t (x,v) = \frac1N \sum_{i = 1}^N \delta _{x - x^\eps_i(t)} \delta _{v- v^\eps_i(t)}\,,
\end{equation}
where $(x^\eps_i(t),v^\eps_i(t))_{i \leq N}$ stands for the positions 
and velocities of the $N$ particles at time~$t$ and~$\delta _x$ stands for the Dirac mass at $x=0$.
 This measure is completely symmetric (i.e. invariant under any permutation of the particle labels) 
 due to the exchangeability assumption.
However this first averaging is not enough to obtain a simple description of the dynamics when $N$ is large, because of the instabilities mentioned in the previous section (see Figure \ref{fig: instability}) which lead to a strong dependency in $\eps$ of the particle trajectories.
We will therefore introduce a second averaging with respect to initial configurations.

From the physical point of view, this averaging is natural as  only fragmentary information on the initial configuration is available. 
A natural starting point is the particle distribution $f^0 = f^0(x,v)$ 
 which prescribes the probability for a particle to be at position~$x$ with velocity $v$. As $N$ is large, we  assume that the initial data $(X_N, V_N) =(x_i, v_i)_{1\leq i \leq N}$ are independent random variables identically distributed according to $f^0$. This assumption has however to be slightly corrected in order to take into account the exclusion between particles~$|x_i - x_j|>\eps$ for $i\neq j$. This statistical framework is referred to as the \emph{canonical} ensemble~\cite{Ruelle_livre}. 
 This is a simple framework to derive rigorous foundations for the kinetic theory, i.e. 
 to characterise, in the large $N$ asymptotics,  the average dynamics and more precisely
the evolution equation governing  the distribution $f(t,x,v)$  at time $t$ of a typical particle.

\medskip

In this paper, our goal is actually to go beyond this average dynamics, and to understand in a  fine way the correlations arising dynamically inside the gas. Fixing a priori the number $N$ of particles induces additional correlations and thus  technical difficulties. 
To bypass them, we introduce a third level of averaging, by assuming that the number $N$ of particles  is also a random variable, and that only its average tuned by $\mu_\eps = \eps^{- (d-1)}$ is determined according to the Boltzmann-Grad scaling. 
Roughly speaking, $N$ is chosen according to  a Poisson law of mean close to $\mu_\eps$, and then   for any fixed $N$, the variables $(X_N, V_N) $ are identically distributed, and independent up to the spatial exclusion. 
More precisely, the variables~$(N, X_N, V_N)$ are chosen jointly under the so-called \emph{grand canonical}  measure which will be introduced later in \eqref{eq: initial measure}. This is referred to as the \emph{grand canonical ensemble} 
and from now on, we will use this setting.

We therefore seek to understand the statistical behavior of the empirical measure 
\begin{equation}
\label{pi-def}
 \pi^\eps _t (x,v) = \frac1{\mu_\eps}  \sum_{i = 1} ^N \delta _{x - x^\eps_i(t)} \delta_ {v- v^\eps_i(t)}\,,
 \end{equation}
where the initial configuration $(N, (X^{\eps0}_N = X_N , V^{\eps0}_N = V_N))$ is a random variable, but the microscopic dynamics is completely deterministic (governed by the hard sphere equations represented in  Figure~\ref{HS-fig}).

\subsection{A probabilistic approach}
\label{questions-sec}

The first question is to determine the law of large numbers, that is the limiting distribution of a typical particle when $\mu_\eps \to \infty$. In the case of $N$ independent identically distributed variables $(\eta_i)_{1\leq i \leq N}$, the law of large numbers implies in particular that the average converges in probability to its expectation
$$ 
\frac1N \sum_{i = 1} ^N \eta_i  \xrightarrow[N \to \infty] { } \bbE (\eta) .
$$
For the interacting particle system, two difficulties arise.
The first one is that, even at time~0, the variables $(x_i, v_i)_{1\leq i \leq N}$ are weakly correlated due to the exclusion. 
In the low density regime, this problem is well understood by classical methods of equilibrium statistical mechanics
(see e.g. \cite{Ruelle_livre}).
In particular, denoting the average of any continuous test function $h$ under the initial empirical measure by
\begin{equation*}
\la \pi^\eps_0, h\ra = \frac1{\mu_\eps}  \sum_{i = 1} ^N h \big( x^{\eps 0}_i , v^{\eps 0}_i  \big) \,,
\end{equation*}
 the following  convergence in probability holds
$$ 
\la \pi^\eps_0, h\ra - \int f^0 h (x,v) \; dx dv \xrightarrow[\mu_\eps \to \infty] { } 0 
\qquad \text{under the grand-canonical measure}.
$$
We stress the fact that, throughout this paper, the limit $\mu_\eps \to \infty$ implies that 
the sphere diameter $\eps$ tends also to 0 as both parameters are linked by the Boltzmann-Grad scaling~$\mu_\eps \eps^{d-1} = 1$.
The second difficulty which is main challenge is to understand whether the initial quasi-independence is propagated in time so that there exists a function $f(t,x,v)$ such that the following  convergence in probability holds 
\begin{equation}
\label{convergence-empirical}
\begin{aligned} 
\la \pi^\eps_t, h\ra - \int f(t) h \; dx dv \xrightarrow[\mu_\eps \to \infty] { } 0
\qquad   &\text{under the grand-canonical measure}\\
& \text{on initial configurations},
\end{aligned}
\end{equation}
and whether $f(t)$ evolves according to a deterministic  equation, namely the Boltzmann equation.
As we will see, this question is particularly delicate since the Boltzmann equation obtained in the limit is  singular {\color{black} (see \eqref{boltz})}.
The major result  proving this convergence goes back to Lanford \cite{Lanford} and will be explained in Section \ref{lanford-sec}.

\medskip
The approximation \eqref{convergence-empirical} of the empirical measure neglects two types of errors. The first one is the fact that there are corrector terms which converge to 0 as $\mu_\eps \to +\infty$. The second one  is related to the vanishing probability of the initial configurations for which the convergence does not hold. 
A classical question in statistical physics is to quantify more precisely these errors, by studying fluctuations, i.e. deviations between the empirical measure and its expectation. In the case of $N$ independent and identically distributed random variables $(\eta_i)_{1\leq i \leq N}$, 
the central limit theorem implies that the fluctuations are of order $O(1/\sqrt{N})$ 
and   the following convergence in law holds
$$ 
\sqrt{N} \Big(  \frac1N \sum_{i = 1} ^N \eta_i  - \bbE(\eta) \Big) 
 \xrightarrow[N \to \infty] {(law) } \cN (0, \Var(\eta) ) ,
$$
where $\cN (0, \Var(\eta) )$ is the normal law of variance 
$\Var(\eta) = \bbE ((\eta- \bbE(\eta))^2)$. 
In particular, at this scale, some randomness is retrieved. 
Investigating the same fluctuation regime for the dynamics of hard sphere gases
 consists in considering the scaled fluctuation field  $\zeta^\eps_t$ defined by duality 
\begin{equation}
\label{zeta-def}
\la \zeta^\eps_t ,h \ra = \sqrt{\mu_\eps}  \Big(  \la \pi^\eps_t, h\ra   - \bbE_\eps ( \la \pi^\eps_t, h\ra) \Big)\,,
\end{equation}
where $h$ is a continuous test function, and $\bbE_\eps$ denotes the expectation on initial configurations under the grand-canonical measure.
A series of recent works \cite{BGSS1, BGSS2, BGSS3, BGSS4} has allowed to characterize   these dynamical fluctuations, and to derive a stochastic evolution equation governing the limiting process. These results will be presented in Sections \ref{BGSS-sec} and \ref{sec: Close to equilibrium}.

\medskip

The last question generally studied in a classical statistical approach is the one of quantifying rare events, i.e. of estimating the probability of observing an atypical behavior (which deviates macroscopically from the average). For independent and identically distributed random variables, this probability is exponentially small, and it is therefore natural to study the asymptotics
\begin{equation}
\label{eq: LD iid}
I(m) : = \lim_{\delta \to 0} \; \lim_{N\to \infty } 
\; - \frac{1}{N} \log 
\bbP \left( \Big| \frac1N \sum_{i = 1}^N \eta _i - m \Big| < \delta \right)
\quad \text{with} \quad 
m \not = \bbE( \eta) \,.
\end{equation}
The limit $I(m)$ is called the large deviation function and it can be expressed as the Legendre transform of the log-Laplace transform of a single variable $u : \bbR \mapsto \log \bbE \big( \exp ( u \eta) \big)$ \cite{Dembo_Zeitouni}.
To generalise this statement to correlated variables, it is necessary to compute a more global Laplace transform and this requires a control on the correlations with  exponential accuracy. 
The methods of dynamical cumulants introduced in \cite{BGSS1, BGSS2} are a key tool to compute exponential moments of the hard sphere distribution and in this way, to 
control the measure of events up to scales which are vanishing exponentially fast.
 We will   give a flavour of those techniques in Section  \ref{BGSS-sec}.

Note that precise conjectures regarding those three questions are formulated by Rezakhanlou in \cite{RezakhanlouLNM}.


\section{TYPICAL DYNAMICAL BEHAVIOR}

\subsection{Boltzmann's great intuition}       
\label{boltz-sec}

The equation which rules the typical evolution of a hard sphere gas was proposed heuristically by Boltzmann \cite{boltzmann} about one century before its   rigorous derivation by Lanford~\cite{Lanford} as the ``limit'' of the particle system when $\mu_\eps \to +\infty$. The revolutionary idea of Boltzmann was to write an evolution equation for the probability density $f = f(t,x,v)$ giving the proportion of particles at position $x$ with velocity $v$ at time $t$. In the absence of collisions and in a domain without boundary, this density $f$ would be exactly transported along the physical trajectories $x(t) = x(0)+vt$, meaning that $f(t,x,v) = f^0(x-vt, v)$. The difficulty consists then in taking into account the statistical effect of collisions. Insofar as the size of the particles is negligible, one can consider that these collisions are pointwise both in $t$ and~$x$. Boltzmann proposed therefore a rather intuitive counting~:
\begin{itemize}
\item the number of particles with velocity $v$ is increased  when a  particle of velocity~$v'$ collides with a particle of velocity $v'_1$, and jumps to velocity $v$ (see \eqref{eq: scattering rules}).
  Notice that here, $( v ', v_1' )$ play the role of precollisional velocities, while instead in Figure~\ref{HS-fig} 
this notation was used for the postcollisional velocities in the particle system;
\item the number of particles with velocity $v$ is  decreased when a particle of velocity $v$ collides with a particle of velocity   $v_1$, and is deflected into another velocity.
\end{itemize}
The probability of these jumps is described by a transition rate, referred to as the \textit{collision cross section} $b$. The function $b(v,v_1,\omega)$ is non negative, depends only on the relative velocity~$|v-v_1|$ 
  and on the angle between $(v-v_1)$ and $\omega$, a scattering vector which is distributed uniformly in the unit sphere $\bbS^{d-1} \subset \bbR^d$.
For the hard sphere interaction, we shall see that~$\omega$ keeps track of the way two hard spheres collide (see Figure \ref{HS-fig}) and that~$b(v-v_1, \omega) = \big( (v-v_1) \cdot \omega \big)_+$.
 In particular, it is invariant under $(v, v_1) \mapsto (v_1,v)$ (exchangeability) and under $(v,v_1, \omega)\mapsto (v', v'_1, \omega)$ (microscopic reversibility).

\medskip

The fundamental assumption in Boltzmann's theory is that, in a rarefied gas, the correlations between 
 two particles about to collide should be very weak. Therefore the joint probability to have both precollisional particles  of velocities $v$ and $v_1$ at position $x$ at time~$t$ should be  well approximated by $f(t,x,v) f(t,x,v_1)$. This independence property is called the molecular \textit{chaos assumption}.
The equation then states
\begin{equation}
\label{boltz}
\left\{ 
\begin{aligned}
&\d_t f+ \underbrace{v\cdot \nabla_x f  }_{\mbox{\footnotesize{transport} }} = \underbrace{C(f,f)  }_{\mbox{\footnotesize{collision} }} \\
&C(f,f) (t,x,v) = \iint \Big[  \underbrace{f(t,x,v') f(t,x,v'_1)  }_{\mbox{\footnotesize{gain term} }} -  \underbrace{f(t,x,v) f(t,x,v_1) }_{\mbox{\footnotesize{loss term} }}  \Big]  \underbrace{b(v-v_1, \omega) }_{\mbox{\footnotesize{cross section} }}  dv_1 d\omega
\end{aligned}\right.
\end{equation}
where the scattering rules 
\begin{equation}
\label{eq: scattering rules}
v' = v -  \big( (v-v_1) \cdot \omega \big) \, \omega,
\qquad 
v_1' = v_1 +  \big( (v-v_1) \cdot \omega \big) \, \omega
\end{equation}
are analogous to the microscopic collision rules introduced in Figure \ref{HS-fig}, with the important difference that $\omega$ is now a random vector chosen uniformly in the unit sphere $\bbS^{d-1} \subset \bbR^d$.  Indeed the relative position of the colliding particles has been forgotten in the limit $\varepsilon \to 0$.  As a consequence, the Boltzmann equation is singular as it  involves a product  of  densities at the same point $x$.


\medskip

Boltzmann's idea of reducing to a kinetic equation 
the Hamiltonian dynamics describing the  atomistic behaviour,  was revolutionary and   opened the way to 
 the description of non-equilibrium phenomena by mesoscopic equations.
However, the Boltzmann equation 
 was first  heavily criticised as it seems to violate some basic physical principles.  
Indeed, what made Boltzmann's theory such a breakthrough, but also made it unacceptable for many of his contemporaries, is that it predicts a time  irreversible evolution, providing actually a quantitative formulation of the second principle of thermodynamics. The Boltzmann equation (\ref{boltz}) indeed has  a Lyapunov functional defined by  $S(t) = -\iint f\log f (t,x,v) dxdv $ and referred to as the entropy, which can only increase along the evolution $\frac{d}{dt} S(t) \geq 0$, with equality if and only if the gas is at thermal equilibrium.  At first sight, this irreversibility does not seem to be compatible with the fact that the hard sphere dynamics is governed by a Hamiltonian system, i.e. a system of ordinary differential equations which is completely time reversible. 
Soon after Boltzmann postulated his equation, these two different behaviours were considered, by Loschmidt, as a paradox  and an obstruction to Boltzmann's theory. 
A fully satisfactory mathematical explanation of this issue remained open during almost one century, until the role of probability was precisely identified~: the underlying dynamics is reversible, but the description which is given of this dynamics is only partial (obtained by averaging or looking at the most probable path) and therefore is not reversible.
%

\subsection{Lanford's theorem}
\label{lanford-sec}

Lanford's result \cite{Lanford} shows in which sense the Boltzmann equation~(\ref{boltz}) is a good approximation of the hard sphere dynamics. 
Let us first define the initial distribution.

\smallskip

\noindent
\textbf{Initial data.}
\label{Def: initial data}
\textit{
Consider $\bbT^d = [0,1]^d$ the unit domain with periodic boundary conditions and~$f^0= f^0(x,v)$ a  Lipschitz   probability density in $\bbT^d \times \bbR^d$, 
with   Gaussian tails  
at large velocities. 
To define a system of hard spheres which are initially  independent (up to the exclusion) and   identically distributed according to $f^0$, we introduce the grand canonical measure : the probability density of finding $N$ particles with coordinates~$Z_N = (x_i,v_i)_{i \leq N}$ is given by
\begin{equation}
\label{eq: initial measure}
\frac{1}{N!}W^{\eps}_{N}(Z_N) 
= \frac{1}{\cZ^ \eps} \,\frac{\mu_\eps^N}{N!} \, 
 \, \prod_{i =1}^N f^0 (x_i,v_i) \; \prod_{i\neq j}\indc_{ |x_i - x_j| > \eps}  \, , \qquad \hbox{ for } N= 0, 1,2,\dots 
\end{equation} 
where the constant $\cZ^ \eps$ is the normalisation factor of the probability measure.
Once the random initial configuration is chosen, the hard sphere dynamics evolve deterministically 
and the corresponding probability and expectation on the particle trajectories will be denoted by $\bbP_\eps$ and $\bbE_\eps$.}

%

\smallskip

\noindent
Lanford's result can be stated as follows (this is not exactly the original formulation, see in particular Section~\ref{irreversibility} below for comments).
\begin{Thm}[Lanford]
\label{lanford-thm} 
In the Boltzmann-Grad limit ($\mu_\eps \to \infty$ with $\mu_\eps \eps^{d-1} = 1$), the empirical measure $\pi^\eps_t $ of the hard sphere system defined by (\ref{pi-def})   concentrates  on the solution of the Boltzmann equation (\ref{boltz}), i.e. for any bounded and continuous function $h$
$$
\forall \delta >0, \qquad 
\lim_{\mu_\eps \to \infty} 
\bbP_\eps \left( \Big|\la \pi^\eps_t, h\ra  - \int f(t) h dxdv  \Big| \geq \delta \right) = 0\,,
$$
on a time interval $[0,T_{L}]$ depending only on the initial distribution $f^0$.
\end{Thm}
Let us comment on the time of validity $T_L$ of the approximation. This time depends on the initial data $f^0$ and turns out to be  of the order of a fraction of the mean time between two successive collisions for a typical  particle. This time is large enough for the microscopic system to undergo a large number of collisions (of the order $O(\mu_\eps)$), and in particular irreversibility already shows up at this scale. But this time is (far) too small to see phenomena such as  relaxation towards (local) thermodynamic equilibrium, and a fortiori hydrodynamic regimes. Physically we do not expect   this time to be  critical, in the sense that the dynamics would change nature afterwards.
Actually, in practice the Boltzmann equation is used in many applications (such as  calculations for  the reentrance of spatial vehicles in the atmosphere) without time restriction. However, it is important to note that a time restriction may  not  be only technical~: from the mathematical point of view, one cannot exclude that the Boltzmann equation exhibits singularities (typically spatial concentrations which would prevent making sense of the collision term, and which would also contradict locally the low density assumption). 
In order to construct global in time solutions for the Boltzmann equation, one actually has  either to consider small fluctuations around some equilibrium, or to introduce a renormalisation procedure 
\cite{DiPerna_Lions}.
These two approaches rely strongly on entropy production estimates, which do not have any counterpart at the microscopic level (i.e. for fixed $\mu_\eps$, $\eps$). In the current state of our knowledge, the problem of extending Lanford's convergence result to longer times faces serious obstructions, even to  the time of existence and uniqueness of the solution to the Boltzmann equation. This will be discussed later on in Section \ref{sec:maindifficulties}
 (see also Section \ref{sec: open problems}). 
In Section \ref{time-sec}, we  will also present some recent results in this direction, 
providing a global in time convergence for the fluctuation field at equilibrium.

\subsection{Heuristics of the proof}
\label{lanford-proof-heuristics}

Let us now explain informally how the Boltzmann equation \eqref{boltz} can be guessed from the particle dynamics. 
The goal is to transport the initial grand canonical measure, defined in \eqref{eq: initial measure}, along the dynamics and then to project this measure at time $t$ on the 1-point particle phase space. We therefore define by duality
 $F_1^\eps (t,z)$  the density of a typical particle with respect to the test function $h$ as
\begin{equation}
\label{eq: F1}
\int F_1^\eps (t, z) h(z) d z 
= \bbE_\eps \left(  \la \pi_t^\eps , h \ra  \right) ,
\end{equation}
where the empirical measure $\pi^\eps _t$ was introduced in \eqref{pi-def}.
More generally, 
we are going  to introduce $\pi_{k,t}^\eps$, the natural extension of the empirical measure $\pi_t^\eps$ to $k$ distinct particles.
For simplicity, the particle coordinates $\big( x^\eps_i (t), v^\eps_i (t) \big)$ at time $t$  will be denoted by $z^\eps_i (t)$.  
For any test function $h_k$ of~$k$ variables, we define
\begin{equation}
\label{eq: empirique k}
\la \pi_{k,t}^\eps , h_k \ra =  \frac1{\mu_\eps^k} \sum_{(i_1,\dots, i_k)} h_k\big( z^\eps_{i_1}(t), \dots, z^\eps_{i_k}(t)\big)  
\end{equation}
and the sum is over  the $k$-tuples of indices  among all the particles at time $t$.
 We stress the fact that $\pi_{k,t}^\eps$ differs from $(\pi_t^\eps)^{\otimes k}$ as the variables are never repeated.
We will study the $k$-particle correlation functions $F_k^\eps $ which are symmetric finite dimensional projections of the probability measure 
\begin{equation}
\label{eq: Fk}
\int F_k^\eps (t, Z_k) h_k(Z_k) dZ_k 
= \bbE_\eps \left(  \la \pi_{k,t}^\eps , h_k \ra  \right) ,
\end{equation}
denoting $Z_k= (x_i, v_i) _{1\leq i \leq k}$.
The correlation functions are key to describe the kinetic limit. In particular,  Theorem \ref{lanford-thm}
shows that~$F_1^\eps (t,z)$   converges  to the solution of the Boltzmann-equation 
$f(t)$ in  the Boltzmann-Grad limit ($\mu_\eps \to \infty$ with $\mu_\eps \eps^{d-1} = 1$). Let us explain briefly why this holds.

Let $h$ be a bounded smooth test function on~${\mathbb T}^d \times {\mathbb R}^d$.
Consider the evolution of the empirical measure  during a short time interval $[t,t+\delta]$ and split the different contributions according to the number of collisions for each  particle 
\begin{align}
\label{eq: short time delta}
\bbE_\eps & \left[ \la \pi^\eps_{t+\delta} , h\ra \right] 
 -  \bbE_\eps \left[ \la \pi^\eps_t , h\ra \right] =
 \bbE_\eps \left[ \frac{1}{\mu_\eps } \sum_{j \atop \text{no collision}} 
\Big( h \big( z^\eps _j(t + \delta) \big) - h \big( z^\eps _j(t) \big) \Big)\right]\\
& + \bbE_\eps \left[ \frac{1}{2\mu_\eps } \sum_{(i, j) \atop \text{with 1 collision}}  \Big(
h \big( z^\eps_i (t + \delta) \big) + h \big( z^\eps _j(t + \delta) \big) 
- h \big( z^\eps _i (t) \big) - h \big( z^\eps _j(t) \big)\Big)
\right] +   O ( \delta^2 )\,  , \nonumber
\end{align}
and we are going to argue that the error term $\delta^2$ takes into account all the groups of particles undergoing at least 2 collisions in the short time interval $\delta$.

\begin{figure}[h] 
\centering
\includegraphics[width=2.2in]{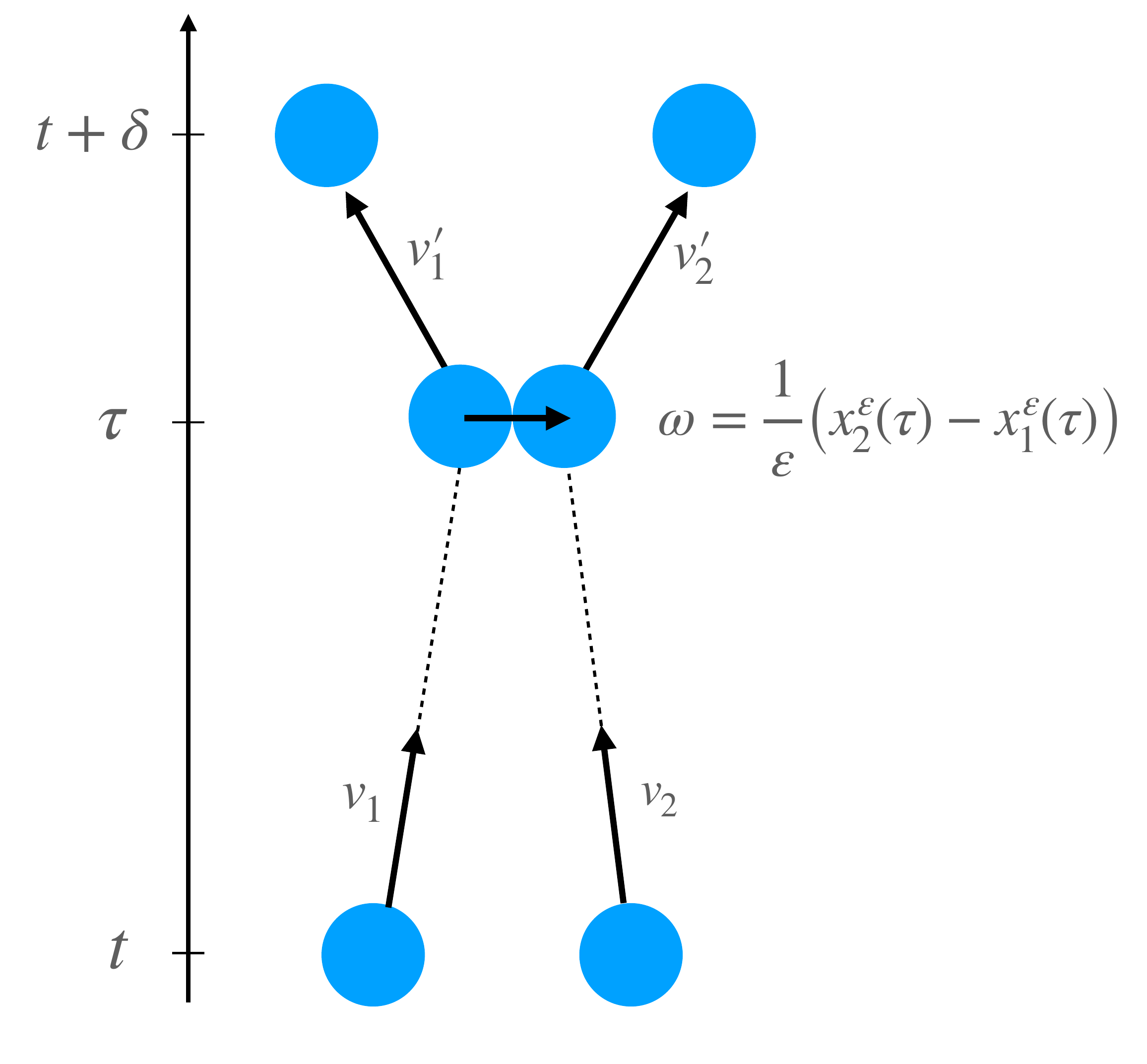} 
\caption{
\color{black}
On the figure, two particles collide  in the time interval $[t,t+\delta]$ according to the scattering rules of Figure \ref{HS-fig}.
The collision occurs at time  $\tau$ if $x_1  - x_2 + (\tau- t) (v_1  - v_2) = -\eps \omega$.
Therefore $x_2$ has to be in a tube with axis $v_1  - v_2$ and the coordinates $z_1,z_2$ at time $t$
 can be parametrised by $( x_1 , v_1,v_2, \tau, \omega)$. 
This change of variables has a Jacobian 
$d z_1 d z _2 = \eps^{d-1}  \; \big( ( v_1 - v_2 ) \cdot \omega \big)_+
d\omega d\tau dx_1 d v_1 d v_2$. 
}
\label{fig: micro_collision}
\end{figure}

The asymptotic behavior when $\delta$ tends to 0 will be analysed now for 
each term in~\eqref{eq: short time delta}.
The transport contribution arises from the particles moving in straight line
without collisions, indeed if the distribution $F_1^\eps$ is smooth enough,
one gets
\begin{align*}
\label{eq: short time delta}
 \bbE_\eps \left[ \frac{1}{\mu_\eps } \sum_{j \atop \text{no collision}} \Big(
h \big( z^\eps _j(t + \delta) \big) - h \big( z^\eps _j(t) \big)\Big) \right]
= \delta \int d z_1    F_1^\eps (t, z_1) v_1 \cdot \nabla_x h (z_1) + o(\delta)\,.
\end{align*}
We turn next to the term involving one collision.
Note first that two particles starting at~$(x_1,v_1)$ and $(x_2,v_2)$ at time $t$  collide at a later time $\tau \leq t + \delta$ if the following geometric condition holds (see  Figure \ref{fig: micro_collision})
\begin{equation}
\label{eq: collision condition tube}
x_1  - x_2 + (\tau- t) (v_1  - v_2)  = -\eps \omega.
\end{equation}
This implies that their relative position must belong to a tube  {\color{black} oriented in the direction
$v_1  - v_2$ with length $\delta |v_1  - v_2|$ and width $\eps$.} This set has a size proportional to
$\delta \eps^{d-1} |v_2 - v_1|$ with respect to the Lebesgue measure. 
More generally, a series of $k-1$ collisions between $k$ particles imposes~$k-1$ constraints of the previous form. Using the Boltzmann-Grad scaling $\mu_\eps \eps^{d-1} = 1$ and neglecting the velocity contribution, one can show that this event has a vanishing probability bounded from above by  
\begin{equation}
\label{eq: size k-1 collisions}
\left( \frac{\delta}{\mu_\eps} \right)^{k-1}. 
\end{equation}
Since there are, on average, $\mu_\eps^k$ ways of choosing these $k$ colliding particles, we deduce that  the occurence of $k-1$ collisions in \eqref{eq: short time delta} has a probability of order $\delta^{k-1} \mu_\eps$. 
 This explains why in \eqref{eq: short time delta} the probability of the terms involving more than 1 collision, i.e. involving $k \geq 3$ colliding particles, has been estimated by $O(\delta^2)$.

This crude estimate is not sufficient to recover the collision operator $C(f,f)$ of the Boltzmann equation \eqref{boltz}.
{\color{black}
We are going now to analyse more carefully the term with one collision in \eqref{eq: short time delta} in order to identify $C(f,f)$. }
As the collision term involves 2 particles, it is no longer a function of the empirical measure.
The correlation function $F_2^\eps$ defined in~\eqref{eq: Fk} will be needed to rewrite it :
\begin{align}
& Coll  = \bbE_\eps \left[  \frac{1}{2\mu_\eps } \sum_{(i, j) \atop \text{with 1 collision}}  \Big(
h \big( z^\eps_i (t + \delta) \big) + h \big( z^\eps _j(t + \delta) \big) 
- h \big( z^\eps _i(t) \big) - h \big( z^\eps _j(t) \big)\Big)
\right]  \nonumber \\
& =  \frac{\mu_\eps}{2}  \int d z_1 d z_2 \; F_2^\eps (t, z_1,z_2) \indc_{\text{$1$ and $2$ collide}}
\; 
\Big[h \big( z_1 (\delta) \big) + h \big( z_2(\delta) \big) - h \big( z_1 \big) - h \big( z_2 \big)\Big]  
+ o( \delta) ,
\label{eq: omega collision 0}
\end{align}
where $z_1 (\delta), z_2(\delta)$ stands for the particle coordinates after a time $\delta$.
After the collision the velocities are scattered to $v_1',v_2'$ according to the deflection parameter $\omega$ (see Figure \ref{fig: micro_collision}), but the positions are almost unchanged as $\delta \ll \eps$. 
Since the function $h$ is smooth, the last term in \eqref{eq: omega collision 0} can be approximated by the velocity jump
\begin{align}
\Delta h (z_1,z_2, \omega) = h \big( x_1, v_1'  \big) + h \big( x_2, v_2' \big) 
- h \big( z_1 \big) - h \big( z_2 \big).
\label{eq: omega collision}
\end{align}
By the condition \eqref{eq: collision condition tube}, it is equivalent to parametrise two colliding particles either by their coordinates $z_1,z_2$ at time $t$ 
or by 
their coordinates at the collision time $\tau$ which are determined by
$x _1,v_1, \tau, \omega, v_2$  (see Figure \ref{fig: micro_collision}). 
This change of variables has a Jacobian $\eps^{d-1} ( (v_1 - v_2) \cdot \omega)_+$.
Since $\eps^{d-1} = 1/ \mu_\eps$ and $\delta \ll \eps$, we deduce from \eqref{eq: omega collision} that 
\begin{align}
\label{eq: jacobien}
Coll =  \frac{1}{2}\int_t^{t + \delta} d \tau \int d z_1 d v_2  d \omega \; F_2^\eps (\tau, z_1, z_2) 
\big( (v_1 - v_2) \cdot \omega\big)_+ \, 
\Delta h (z_1, z_2 , \omega)  + o( \delta) \, ,
\end{align}
with $z_2 = (x_1 + \eps \omega, v_2 )$ as both particles are next to each other at the collision time. The cross section 
{\color{black} $b(v_1 -v_2, \omega) = \big( (v_1 -v_2) \cdot \omega \big)_+$} 
in the Boltzmann equation can be identified from the equation above.
From the previous heuristics, the relation \eqref{eq: short time delta}   provides "almost" a weak formulation  of the collision operator in \eqref{boltz} in the limit $\delta \to 0$
\begin{align}
\label{eq: weak boltzmann}
\partial_t  \int d z_1    \, F_1^\eps (t, z_1) & h(z_1)   
=  \int d z_1    F_1^\eps (t, z_1) v_1 \cdot \nabla h (z_1) \\
& + \frac{1}{2}\int d z_1  \, d \omega \, d v_2 \; \delta_{x_2-x_1 - \eps \omega} \; F_2^\eps (t, z_1, z_2) \big( (v_1 - v_2) \cdot \omega\big)_+
\Delta h (z_1, z_2 , \omega) \nonumber
\end{align}
{\color{black} 
where we used the Dirac notation to stress  that $z_2 = (x_1 + \eps \omega, v_2 )$.}
The key step to close the equation is the  \textit{molecular chaos assumption} postulated by Boltzmann  
  which asserts that  the precollisional particles  remain independently distributed at any time so 
  that 
\begin{align}
\label{eq: molecular chaos}
F_2^\eps (t, z_1, z_2) \simeq F_1^\eps (t,  z_1) F_1^\eps (t, z_2) \, .
\end{align}
When the diameter of the spheres $\eps$ tends to 0, the coordinates $x_1$ and $x_2$ coincide and the scattering parameter $\omega$ becomes a random parameter. 
Assuming that  $F_1^\eps$ converges, then its limit  has to satisfy the Boltzmann equation \eqref{boltz}.
 Establishing rigorously the factorisation~\eqref{eq: molecular chaos} requires   implementing a different and more involved strategy which will be presented in Section \ref{lanford-proof}.

\subsection{Some elements of proof}
\label{lanford-proof}

Lanford's proof \cite{Lanford} has been completed and improved over the years; we refer to the monographs
\cite{spohn2012large, Cercignani_Illner_Pulvirenti, Cercignani_Gerasimenko_Petrina} for accounts of the related results.
In the more recent years, several quantitative convergence results were established, and the proofs extended to the case of compactly supported potentials  \cite{GSRT, Pulvirenti_Saffirio_Simonella, Pulvirenti_Simonella}. 
In the following, we sketch the main steps of the proof for the hard sphere dynamics.

The proof of Lanford's theorem relies on the study of the correlation functions $F^\eps_k$ defined in~(\ref{eq: Fk}), characterising joint probabilities of $k$ particles. In particular, we do not consider directly the empirical measure, but only  its average $F^\eps_1$ under the grand canonical probability $\bbP_\eps$. The starting point is the system of ordinary differential equations for the hard sphere positions and velocities (see Figure \ref{HS-fig}), which provides, by applying Green's formula to the Liouville equation, the following equation on the first correlation function 
\begin{equation}
\label{F1-eq}
\left\{ 
\begin{aligned}
&\d_t F^\eps_1+ \underbrace{v\cdot \nabla_x F^\eps_1  }_{\mbox{\footnotesize{transport} }} = \underbrace{C^\eps (F^\eps_2 )  \,  , }_{\mbox{\footnotesize{collision at distance $\eps$} }}  \\
&C^\eps (F^\eps_2 ) (t,x,v) \\
& \qquad = \iint \Big[  \underbrace{F^\eps_2 (t,x,v', x+\eps \omega, v'_1  )}_{\mbox{\footnotesize{gain term} }} -  \underbrace{F^\eps_2 (t,x,v,x-\eps \omega, v_1) }_{\mbox{\footnotesize{loss term} }}  \Big]  
\underbrace{\big( (v-v_1) \cdot \omega \big)_+ }_{\mbox{\footnotesize{cross section} }}  dv_1 d\omega \, .
\end{aligned}\right.  
\end{equation}
A weak form of this equation has been stated in \eqref{eq: weak boltzmann}.
In the limit $\mu_\eps \to \infty$, 
we expect that it can be closed by the factorisation $F^\eps_2 \sim F^\eps_1 \otimes F^\eps_1$, called the propagation of chaos~\eqref{eq: molecular chaos}.
We are unable to prove it directly, nor will it be shown directly from \eqref{F1-eq} that the limit~$F_1$ of~$F_1^\eps$ satisfies an infinitesimal evolution equation of the previous form. We will rather obtain a series expansion of $F_1$, which will be identified with the solution of the Boltzmann equation by a uniqueness argument.  
The proof is therefore very different from the heuristics presented in Section~\ref{lanford-proof-heuristics}.

\medskip

The proof can  be divided into three steps. The first one is to rewrite $F_1^\eps(t,x,v)$  as an ``average'' (weighted with the initial correlation functions $F_k^{\eps,0}$) of all possible dynamics such that at time $t$, a particle stands at position $x$ with velocity $v$. The analytical way of doing so is to derive  evolution equations similar to (\ref{F1-eq}) for all correlation functions $F_k^\eps$, and then to write the iterated Duhamel formula for this hierarchy of equations, called the BBGKY hierarchy  after Bogoliubov-Born-Green-Kirkwood-Yvon 
  (see \cite{Cercignani_Illner_Pulvirenti} for an account and references).  We will not give the details of these technical computations  here, but will retrieve the final series expansion (formally) using a more probabilistic  perspective based on geometric representations   in terms of pseudo-trajectories.

\begin{figure}[h] 
\centering
\includegraphics[width=4in]{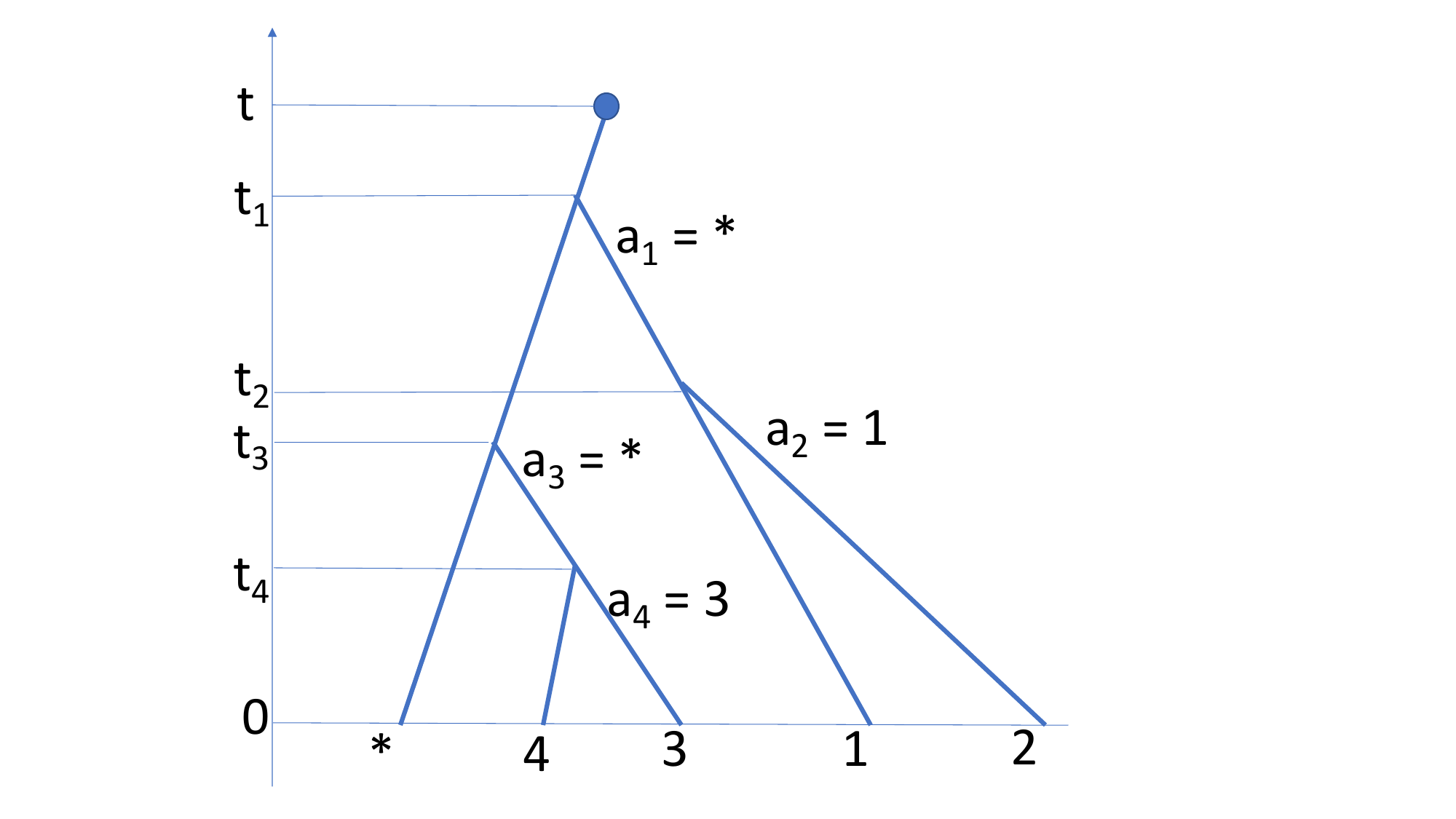} 
\caption{The sequence of collisions in the backward history can be encoded in a tree with the root indexed by the particle $*$ and $n$ branchings (here $n= 4$).
At each creation time, the label of the particle colliding with the fresh particle is indicated. For example at time $t_3$, the particle $*$ collides with particle $3$ so that $a_3 = *$.}
\label{fig: collisiontree}
\end{figure}

The idea is to track back the history of the particle sitting at position $x$ with velocity~$v$ at time $t$, referred to as particle $*$, in order to characterize all initial configurations which contribute to $F_1^\eps(t,x,v)$. We start by following (backward in time) this particle, which has a uniform rectilinear motion $x(t') = x- v (t-t')$ until it collides with another particle, called particle $1$, say at time $t_1$. Note that this collision can actually be either a physical collision (with scattering) or a mathematical artefact coming from the loss term of the equation \eqref{F1-eq} (particles touch each other but are not deflected). 
Thus in order to understand the history of particle $*$, we need to track back the history of both particles $*$ and $1$ before time $t_1$. 
From time $t_1$,  both particles are then transported by the  $2$-particle backward flow until the next collision, say with particle $2$ at time $t_2$,... and we iterate this procedure until  time~$0$. 
Notice that in between the creations of new particles, the particles may collide between themselves as they are transported by the backward hard sphere flow : this will be called    \textit{recollision}.
The history of the particle $*$ can be reconstructed (see Figure \ref{fig: collisiontree}) by prescribing
\begin{itemize}
\item the total number of collisions $n$;
\item the combinatorics of collisions, encoded in a tree $a \in \cA_{1,n}$ with root indexed by the label $*$ and $n$ branchings ($a_i \in \{*, 1,\dots,  i-1\}$ for $1\leq i \leq n$);
\item the collision parameters $(T_n, V_n, \Omega_n) = (t_i, v_i, \omega_i)_{1\leq i \leq n} $ with $0<t_n <\dots <t_1 <t$.
\end{itemize}
We then define the pseudo-trajectory $\Psi^\eps_{1, n} $  starting from $z = (x,v)$ at time $t$  as follows
\begin{itemize}
\item on $]t_i, t_{i-1}[$ the group of $i$ particles is transported by the backward flow;
\item at time $t_i$, particle $i$ is added at position $x_{a_i} (t_i) +\eps \omega_i$, with velocity $v_i$;
\item 
 if the velocities $(v_i, v_{a_i} (t_i^+))$ are postcollisional,
meaning that $(v_{a_i} (t_i^+) - v_i)\cdot \omega_i >0$, then they are instantaneously scattered as in Figure \ref{HS-fig} (with deflection angle~$\omega_i$).
\end{itemize}
We stress the fact that pseudo-trajectories are not particle trajectories of the physical system, but a geometric interpretation of an iterated Duhamel expansion. In particular, pseudo-trajectories do not involve a fixed number of particles, they are coded in terms of random trees (with creation of particles at random times as in Figure \ref{fig: collisiontree}) 
and of signs associated with the gain and loss terms of the collision operator.

Note that not all collision parameters $(T_n, V_n, \Omega_n)  $ are admissible since particles should never overlap. We denote by $\cG^\eps$ the set of admissible parameters.
With these notations, we obtain the following representation of $F^\eps_1 $
\begin{equation}
\label{F1-series}
F^\eps_1 (t,x,v) = \sum_{n= 0}^{+\infty} \sum _{a \in \cA_{1,n} } \int _{\cG^\eps} dT_n dV_n d\Omega_n  \cC(\Psi_{1,n} ^\eps ) F_{1+n} ^{\eps, 0} ( \Psi_{1,n} ^\eps(0))\,,
\end{equation}
where $\Psi_{1,n} ^\eps(0)$ stands for the particle configuration at time~0 of the pseudo-trajectory and 
the term $ \cC(\Psi_{1,n} ^\eps )$ comes from the collision cross-sections 
$$
\cC(\Psi_{1,n} ^\eps ) = \prod_{i = 1} ^n \big( (v_i -  v_{a_i} (t_i^+))\cdot \omega_i \big) \,.
$$
The elementary factor indexed by $i$  is positive if the addition of particle $i$  corresponds to a physical collision (with scattering), and negative if not.

 \begin{Rmk}
 \label{Fk-rmk}
           A similar formula holds for the $k$ point correlation function $F^\eps_k$, except that collision trees $a \in \cA_{k,n}$ have $k$ roots and $n$ branchings.           
\end{Rmk}

\medskip
The formula (\ref{F1-series}) for the first correlation function has been obtained in a rather formal way. In order to study the convergence as $\mu_\eps$ tends to infinity, we need to establish the uniform convergence of   the series (\ref{F1-series}). We  actually use very rough estimates 
{\color{black} (forgetting in particular  the signs of the gain and loss terms in \eqref{F1-eq}, although the cancellations between these different contributions should improve the estimates)} and prove that the series is  absolutely convergent for short times uniformly with respect to $\eps$. Note that this is the only argument in the proof which requires a restriction on short kinetic times.

Let us now estimate the size of the term in \eqref{F1-series}  corresponding to $n$ branchings.
The different contributions are ~:  
\begin{itemize}
\item a combinatorial factor taking into account all the branching choices $|\cA_{1,n}| = n!$;
\item  the volume $t^n/n!$ of the simplex in time $\{ t_n < \dots < t_1 < t \}$;
\item  the $L^\infty$ norm of $F_{1+n} ^{\eps, 0}$ which grows like $ \| f^0\|_\infty ^n$. 
\end{itemize}
This leads to an upper bound of the form $(C\| f^0\|_\infty t)^n$ which implies that the series is absolutely convergent uniformly in $\eps$ on a small time interval depending only on a (weighted)~$L^\infty$ norm of $f^0$.

\begin{Rmk}
 For the sake of simplicity, we do not discuss here the problem of large velocities which create a divergence in the collision cross-section $  \cC(\Psi_{1,n} ^\eps ) $. It can be dealt with similar but more technical arguments, introducing  weighted functional spaces encoding the exponential decay of correlation functions $F_{1+n} ^{\eps, 0} $ at large energies. 
 \end{Rmk}

\begin{figure}[h] 
\centering
\includegraphics[width=4in]{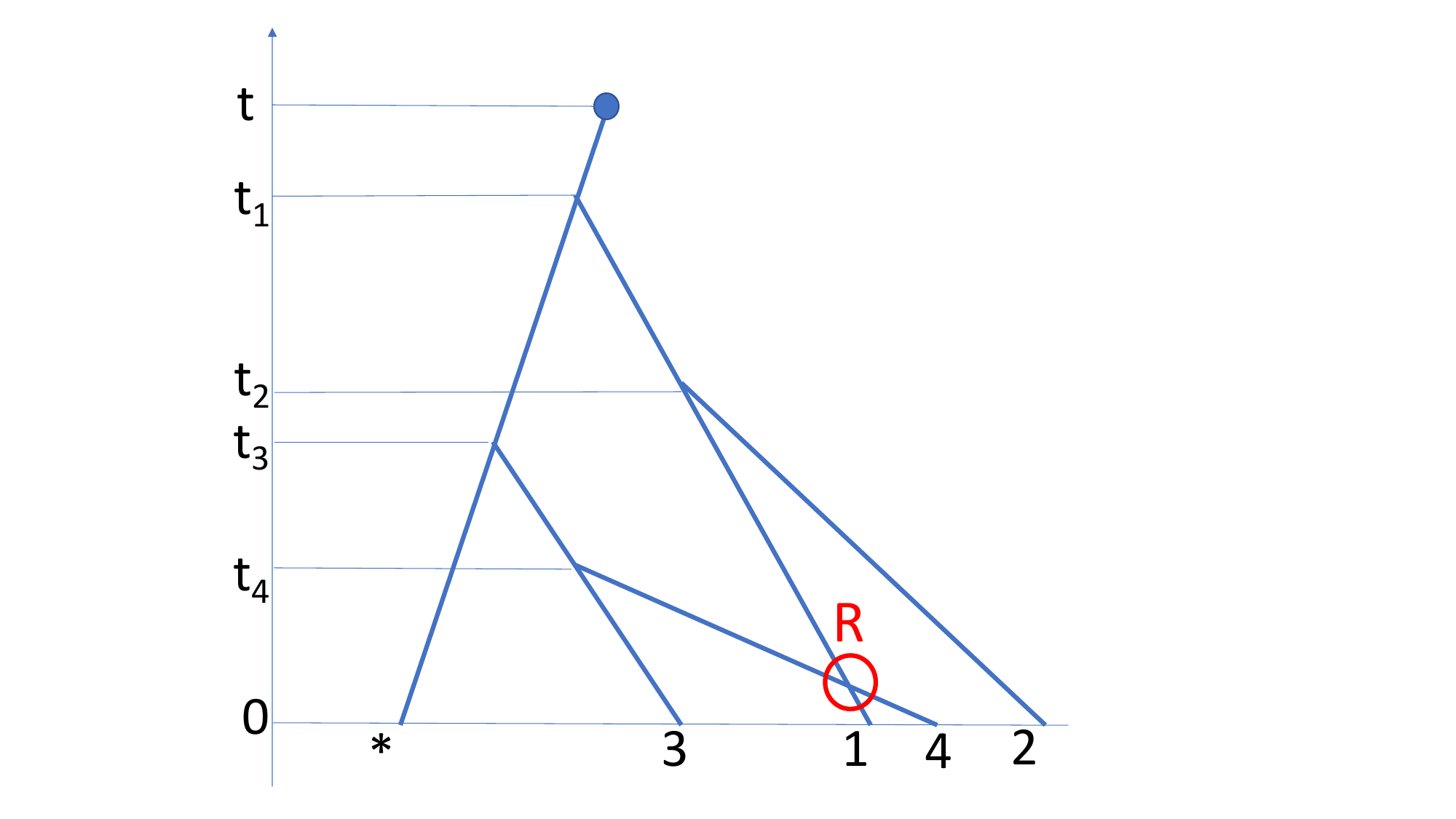} 
\caption{When two particles recollide in the backward flow for fixed $\eps$, their velocities are scattered and the resulting pseudo-dynamics is quite different from the Boltzmann pseudo-dynamics. The sets 
 $\widehat B^\eps_n$   are the sets of integration parameters  leading to at least a recollision within a pseudo-trajectory (as on the picture with $n = 4$).
}
\label{recollision-fig}
\end{figure}

The convergence of $F^\eps_1$ as $\mu_\eps \to \infty$  will then follow termwise. In this third step of the proof, we therefore fix the number $n$ of branchings, as well as the collision tree $a \in \cA_{1,n}$. 
One goal is to understand the asymptotic behavior of the pseudo-trajectories $\Psi_{1,n} ^\eps$. Going back to their definition, we see that it is natural to define limit pseudo-trajectories $\Psi_{1,n}$ 
  (when $\mu_\eps$ tends to $\infty$)  as follows
\begin{itemize}
\item on $]t_i, t_{i-1}[$ the group of $i$ particles is transported by the backward free  flow (since the particles become pointwise in the limit, they cannot see each other);
\item at time $t_i$, particle $i$ is added at position $x_{a_i} (t_i^+) $, with velocity $v_i$ (the spatial shift at the creation time disappears);
\item if the velocities $(v_i, v_{a_i} (t_i^+))$ are post-collisional, then they are scattered (with deflection angle $\omega_i$).
\end{itemize}
Note that in the limit, all collision parameters are admissible (since the non overlap condition disappears).
With this definition of $\Psi_{1,n}$, we see that there is a very natural coupling between $\Psi^\eps_{1,n}$ and $\Psi_{1,n}$~: in most cases, the velocities are exactly equal and the positions differ at most by $n \eps $.
The only problem is when two particles of size $\eps$ recollide (see Figure \ref{recollision-fig}) in the backward flow on some interval $]t_i, t_{i-1}[$~: in this case they are deflected, and the pseudo-trajectory $\Psi^\eps_{1,n}$  is no longer close to  
$\Psi_{1,n}$  on $[0, t_{i-1}]$.
  We therefore split the set of collision parameters $(T_n, V_n, \Omega_n)$ into two parts (and correspondingly split each term in~\eqref{F1-series} into two integrals) : the first subset corresponds to admissible integration parameters such that there is no recollision in $\Psi_{1,n}^\eps$, and the second subset, denoted by $\widehat B^\eps_n$ corresponds either to non admissible integration parameters (leading to some overlap) or to integration parameters for which $\Psi_{1,n} ^\eps$ has at least one recollision.
Using the coupling between $\Psi^\eps_{1,n}$ and $\Psi_{1,n}$ and the regularity of the initial limiting correlation functions (which are nothing else than~$(f^0) ^{\otimes (1+n)}$), we obtain easily the convergence of the first integral.
It remains then to prove that the set $\widehat B^\eps_n$ has vanishing measure so that the corresponding integral has a negligible contribution. The recollision (or overlap) condition implies that the relative velocity between the two recolliding particles $j_1$ and $j_2$  has to be in a small cone, which imposes strong constraints on the last creation involving either $j_1$ or $j_2$. 
We do not detail  these geometric estimates here, but they are quite explicit and provide the following rate of convergence for $t$ sufficiently small (independently of $\eps$)
$$ 
\| F^\eps_1 (t) -  F_1 (t) \|_\infty \leq C \eps ^\alpha \hbox{ for any } \alpha <1
\,,
$$
provided that $f^0$ is Lipschitz.
This concludes the proof, as the series expansion defining $F_1$ turns out to be the (unique) solution of the Boltzmann equation with initial data $f^0$. 
Note that the convergence still holds if $f^0$ is only continuous, but, in that case, we lose the explicit  rate of convergence.  
 \begin{Rmk}
 \label{Propag-chaos-quantitative}
           Actually one can prove (see \cite{BGSS})  the following quantitative propagation of chaos, where the sets~$ \cB_k^\eps$ have vanishing measure
\begin{equation}
\label{eq: factorisation Fk}
\sup_{ t \leq T_L} 
\sup_{Z_k \not \in \cB_k^\eps} \left|
F_k^\eps (t,Z_k)  - \prod_{i=1}^k  f (t,z_i) \right| \leq C^k \eps^\alpha\,,
\end{equation}
for some $\alpha >0$ and a constant $C$ depending on the initial measure $f^0$.   
This is a much stronger notion of convergence than the one stated in Theorem~\ref{lanford-thm}.      
\end{Rmk}

\subsection{On the irreversibility}
\label{irreversibility}

In this paragraph, we are going to argue that the answer to the irreversibility paradox is hidden in the 
chaos assumption \eqref{eq: molecular chaos} which holds only for specific configurations.
Understanding the range of validity of the chaos assumption will be the key to derive not only the Boltzmann equation, but also the stochastic corrections.

Actually  the notion of convergence which appears in the statement of Theorem~\ref{lanford-thm} differs slightly from the one   used in the proof (see Section~\ref{lanford-proof})~: Theorem~\ref{lanford-thm}  states the convergence of observables $ \la \pi^\eps_t, h\ra$, that is a convergence in the sense of measures since the test function $h$ has to be continuous. 
This convergence is rather weak and is actually not enough to   ensure the stability of the  collision term in the Boltzmann equation since this term involves traces.
In the proof of Lanford's Theorem, one actually considers all the correlation functions $F_k^\eps$ introduced in \eqref{eq: Fk}, and  one shows that each one of these correlation functions converges uniformly outside a set $\cB_k^\eps$ of vanishing measure when~$\mu_\eps$ tends to infinity  (see Remark~\ref{Propag-chaos-quantitative}).
\begin{figure}[h] 
\centering
\includegraphics[width=4in]{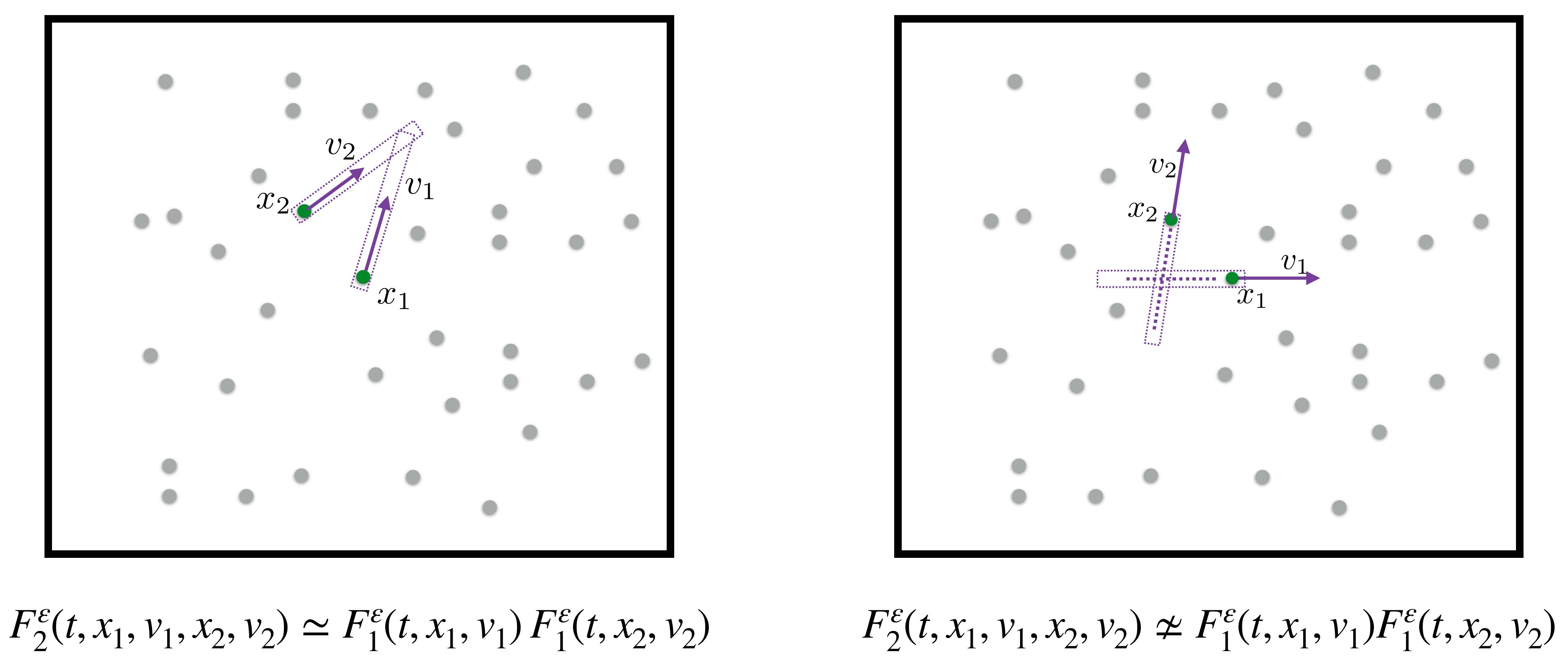} 
\caption{On the left figure, particles $1$ and $2$ will encounter in the future so   they are likely not to have  collided in the past and we expect that the correlation function $F^\eps_2$ factorises in the limit~$\mu_\eps \to \infty$. On the right figure, the particle coordinates belong to the bad  set $\cB_2^\eps$, meaning that they have met in the past. 
In this case,  microscopic  correlations have been built dynamically and the factorisation \eqref{eq: molecular chaos} should not  be valid.
The sets leading to a forward or a backward collision have a similar geometric structure and a similar size which vanishes with respect to the Lebesgue measure when $\eps$ tends to 0. However they play  different roles : the memory of the system is encoded in the sets $\cB_2^\eps$, instead the forward sets are the only ones relevant for the chaos assumption.
The sets $\cB_k^\eps$ are built similarly in terms of the backward flow of $k$ particles (see \cite{BGSS}).
}
\label{fig: bad sets}
\end{figure} 
Moreover the set $\cB_k^\eps$ of bad microscopic configurations  $(t,Z_k)$ (on which $F^\eps_k$ is not converging) is somehow transverse to the set of precollisional configurations 
(as can be seen in Figure \ref{fig: bad sets}, two particles in~$\cB_2^\eps$ tend to move far apart so that they  are  unlikely to collide).
The convergence defect is therefore not an obstacle to taking limits in the collision term, however, these singular sets $\cB_k^\eps$ carry  important information on the time correlations~: 
in particular they encode the memory of the evolution and by neglecting them it is no longer possible to reverse time and to retrace the dynamics backwards.
Thus by discarding  the microscopic information encoded in $\cB_k^\eps$, one can only recover an irreversible kinetic description which is far from describing the complete microscopic dynamics.
The singular sets $\cB_k^\eps$ have been described in \cite{vanbeijeren_lanford_lebowitz_spohn,  denlinger, BGSS} and their complex structure has been made more precise 
in \cite{BGSS2} by means of the cumulants which will be introduced in Section \ref{cumulant-proof}.


\section{CORRELATIONS AND FLUCTUATIONS}

\subsection{From instability to stochasticity}
\label{sec: From instability to stochasticity}

In order to understand the specific features of the hard sphere dynamics {\color{black}  in the low density regime (dilute Boltzmann-Grad limit), it is worthwhile to compare its behaviour to the mean field dynamics.
For this let us consider more general  microscopic dynamics interpolating between the short range and the mean field regimes.
For a given number~$N$ of particles, we set
\begin{equation*}
\forall i \leq N, \qquad
\frac{d}{dt} x_i = v_i, \qquad \frac{d}{dt} v_i = 
-  \frac{1}{N \, \lambda^d} \, \sum_j  \;  \nabla \Phi \left( \frac{x_i - x_j}{\lambda} \right)  ,
\end{equation*}
for some {\color{black}  smooth repulsive (radial decreasing)  potential} $\Phi : [0,1]^d \to \bbR^+$ and a fixed parameter $\lambda \in (0,1]$. 
 This  dynamics is Hamiltonian and 
 by choosing $\lambda = \eps$ (with $N \eps^{d-1} =1$), one recovers  dynamics with a short range potential
which behaves qualitatively as the hard sphere gas and which follows a Boltzmann equation in the   limit \cite{GSRT,Pulvirenti_Saffirio_Simonella}.
For fixed $\lambda$ however, say $\lambda =1$, the limiting behaviour is mean field like and the typical density follows the  Vlasov equation \cite{Braun_Hepp}}
\begin{equation*}
\partial_t f (t,x,v) + v \cdot \nabla_x f (t,x,v) 
=  \left( \int dy dw f(t,y,w) \nabla \Phi ( x - y) \right)  \cdot \nabla_v f(t,x,v).
\end{equation*}
The Vlasov equation has very different properties from the Boltzmann equation, in particular it is reversible, as the microscopic dynamics. 
Furthermore, contrary to the hard sphere dynamics, the precise structure of the initial data plays no role in the limiting behaviour and it has even been shown in \cite{Braun_Hepp} that the fluctuations of the initial data are simply transported by the linearised Vlasov equation. 
Finally we stress the fact that the chaos assumption \eqref{eq: molecular chaos} is known to be propagated in a very strong sense 
for the mean field dynamics \cite{Golse2016, Jabin_Wang2016}.


A drastic difference between the two regimes comes from the fact that the mean field dynamics is not sensitive to a small shift of the coordinates, as the function $\Phi$ is smooth for fixed $\lambda$. This is not the case for the choice $\lambda = \e$ in the Boltzmann-Grad limit. Indeed in the latter situation the scattering behaves qualitatively as in Figure 2, where asymptotically for $\e$ small the deflection parameter decouples completely from the positions and becomes random (cf. Section \ref{lanford-proof-heuristics}). This gives a probabilistic flavour to the surface integral in Boltzmann's collision operator.
%
As we shall see in Theorem~\ref{BGSS-thm1}, the corrections to the limiting Boltzmann equation 
are driven by a stochastic noise which is also generated by the dynamical instabilities.
Thus the limiting structure of the hard sphere dynamics behaves qualitatively as a stochastic process, combining free transport and a random jump process in the velocity space.
  Notice that in the mean field regime, some instability remains for large times  $O(\mu_\eps)$ and this is expected to lead to the Lenard-Balescu stochastic correction \cite{Nota_Velazquez_Winter, duerinckx_saint-raymond}.

The crucial role of randomness in the low density limit was understood by Mark Kac. He devised a purely stochastic process  \cite{Kac_kinetic} whose limiting distribution is  a solution to the homogeneous Boltzmann equation. 
Mathematically, at the microscopic level, this model has a very different structure from the Hamiltonian dynamics previously mentioned.
Indeed, it is a Markov chain restricted only to  particle velocities  and the collisions are modelled by a jump process with a random deflection parameter. 
For Kac's model, the chaos assumption has been derived in a very strong sense \cite{Mischler_Mouhot}.

In the following sections, we are going to argue that the hard sphere dynamics shares, however, many similarities with Kac's model, not only at the typical level, but also at the level of the fluctuations and of the large deviations. 
In this respect, random modelling is an excellent approximation of the hard sphere dynamics.
The key step to accessing  these refined statistical informations will be to understand more precisely the chaos assumption~\eqref{eq: molecular chaos}.

\subsection{Defects in the chaos assumption}
\label{sec: Defects in the chaos assumption}

\begin{figure}[h] 
\centering
\includegraphics[width=3in]{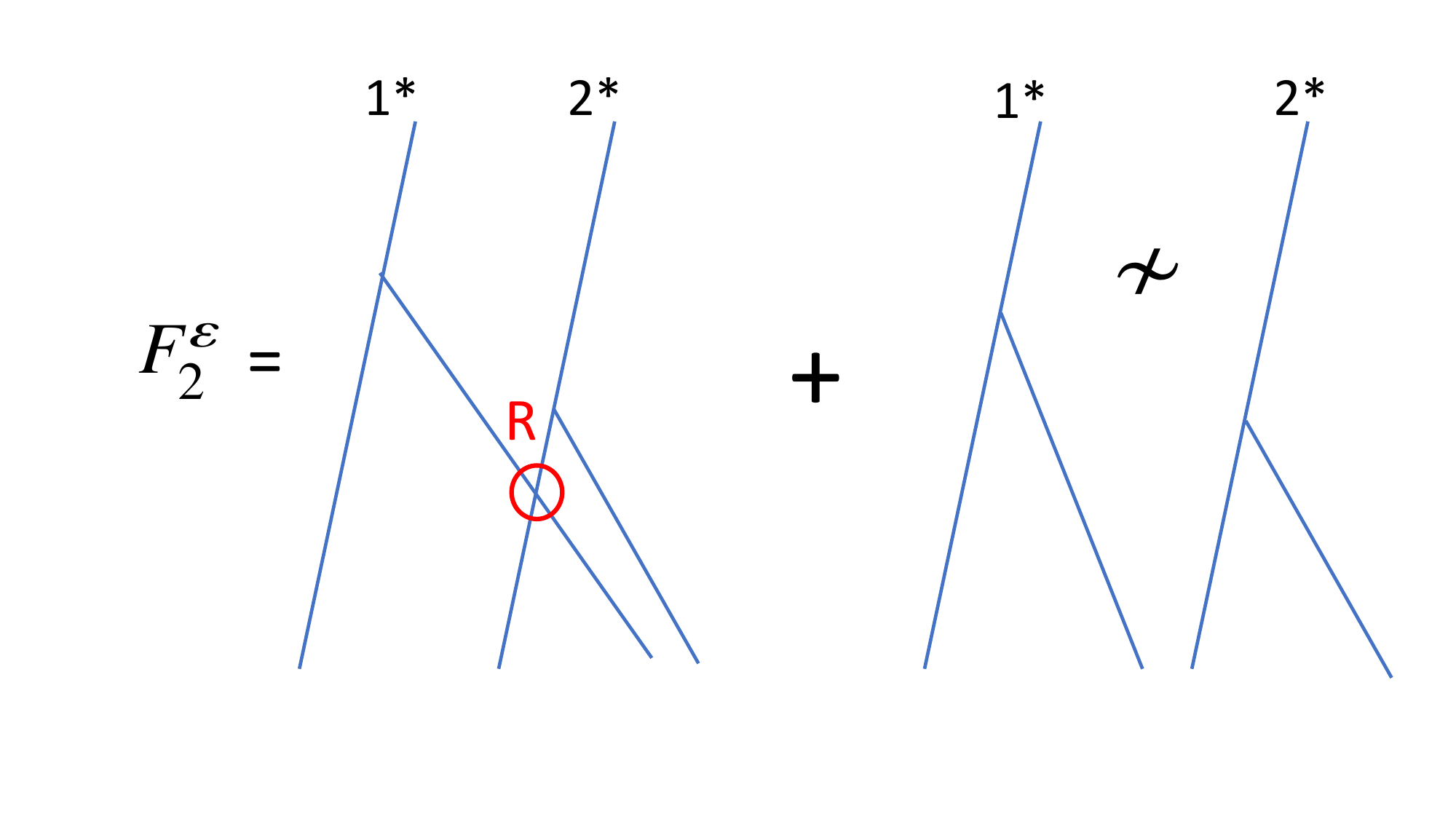} 
\caption{
Among the pseudo-dynamics describing $F_2^\eps$, 
we separate those having a recollision between trees $1^*$ and $2^*$, 
and those where particles from tree $1^*$ and particles from tree $2^*$ 
remain at a distance greater than $\eps$, which will be denoted by $\nsim$. 
In this picture, $n_1= n_2 = 1$.
}
\label{classificationF2-fig}
\end{figure}

Going back to the equation \eqref{F1-eq} on $F_1^\eps$, one can see that up to the small spatial shifts in the collision term (known as Enskog corrections to the Boltzmann equation), deviations from the Boltzmann dynamics are due to the defect of factorization $F_2^\eps- F_1^\eps \otimes  F_1^\eps$, the so-called second order cumulant. In terms of our geometric interpretation, this corresponds to pseudo-trajectories which are correlated. Recall that  $F_2^\eps$ can be described by interacting collision trees with two roots, say labeled by $1^*$ and $2^*$, and $n_1+n_2$ branchings (see Remark~\ref{Fk-rmk}), while the tensor product is described by two independent collision trees each with one root, and $n_1$, resp. $n_2$ branchings. 
The main difference when building the  pseudo-dynamics  corresponding to $F_2^\eps$  is that particles from tree $1^*$ and $2^*$ may (or may not) interact. 
We start by extracting the  pseudo-trajectories of $F_2^\eps$  having at least one interaction between the two trees, which will be called an  \textit{external recollision} (see Figure \ref{classificationF2-fig}) 
in contrast with a recollision inside a collision tree which will be called 
\textit{internal}.

\begin{figure}[h] 
\centering
\includegraphics[width=3in]{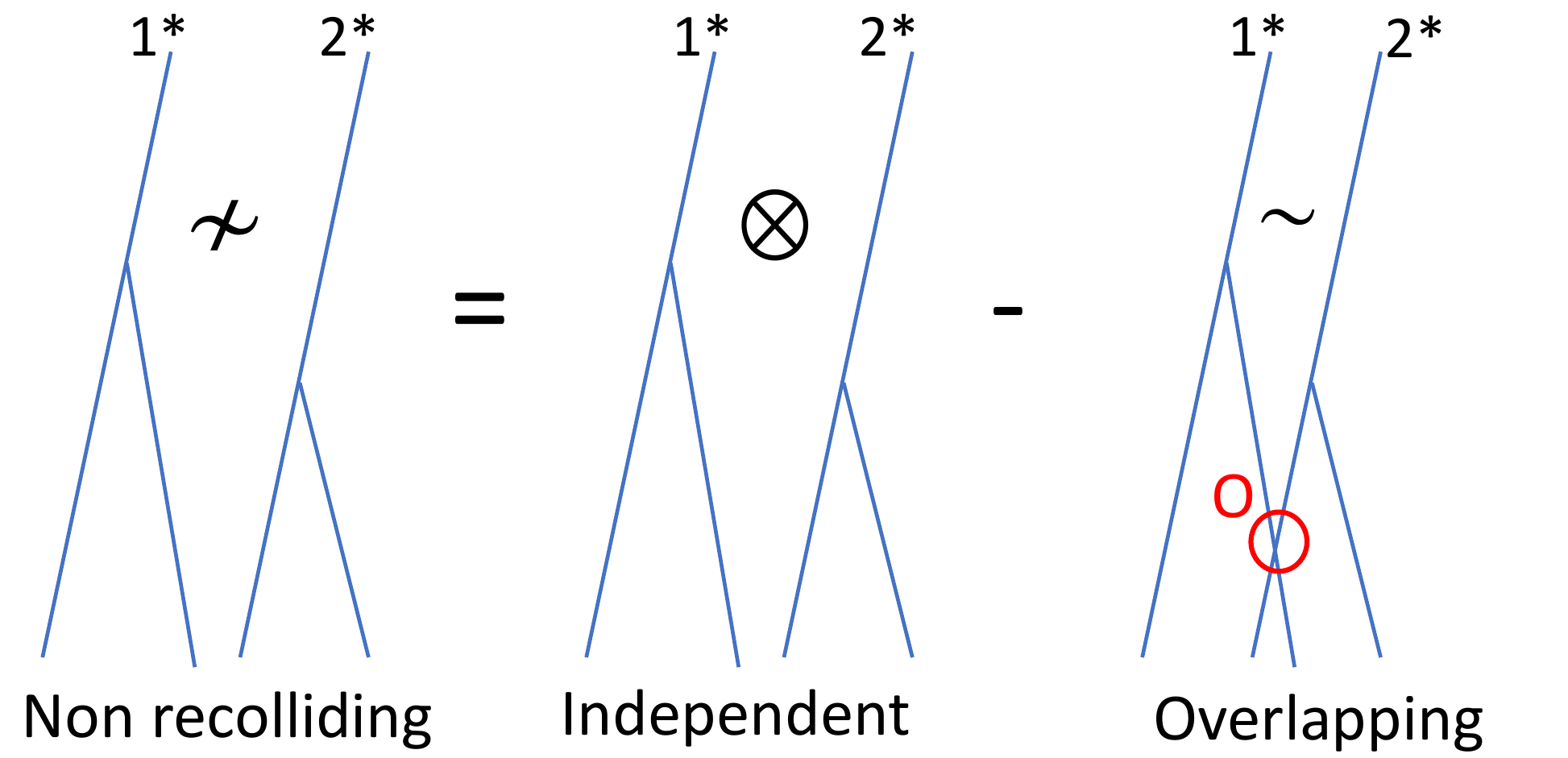} 
\caption{Expanding the dynamical exclusion condition leads to the definition of   overlaps.}
  \label{overlap-fig}\end{figure}

We stress that pseudo-dynamics without external recollision are not  independent since they  satisfy a dynamical exclusion condition. We  therefore decompose the exclusion condition $ \indc_{ 1^* \not\sim 2^*} = 1 -  \indc_{ 1^* \sim 2^*}$ (see Figure \ref{overlap-fig}). 

\begin{figure}[h] 
\centering
\includegraphics[width=3in]{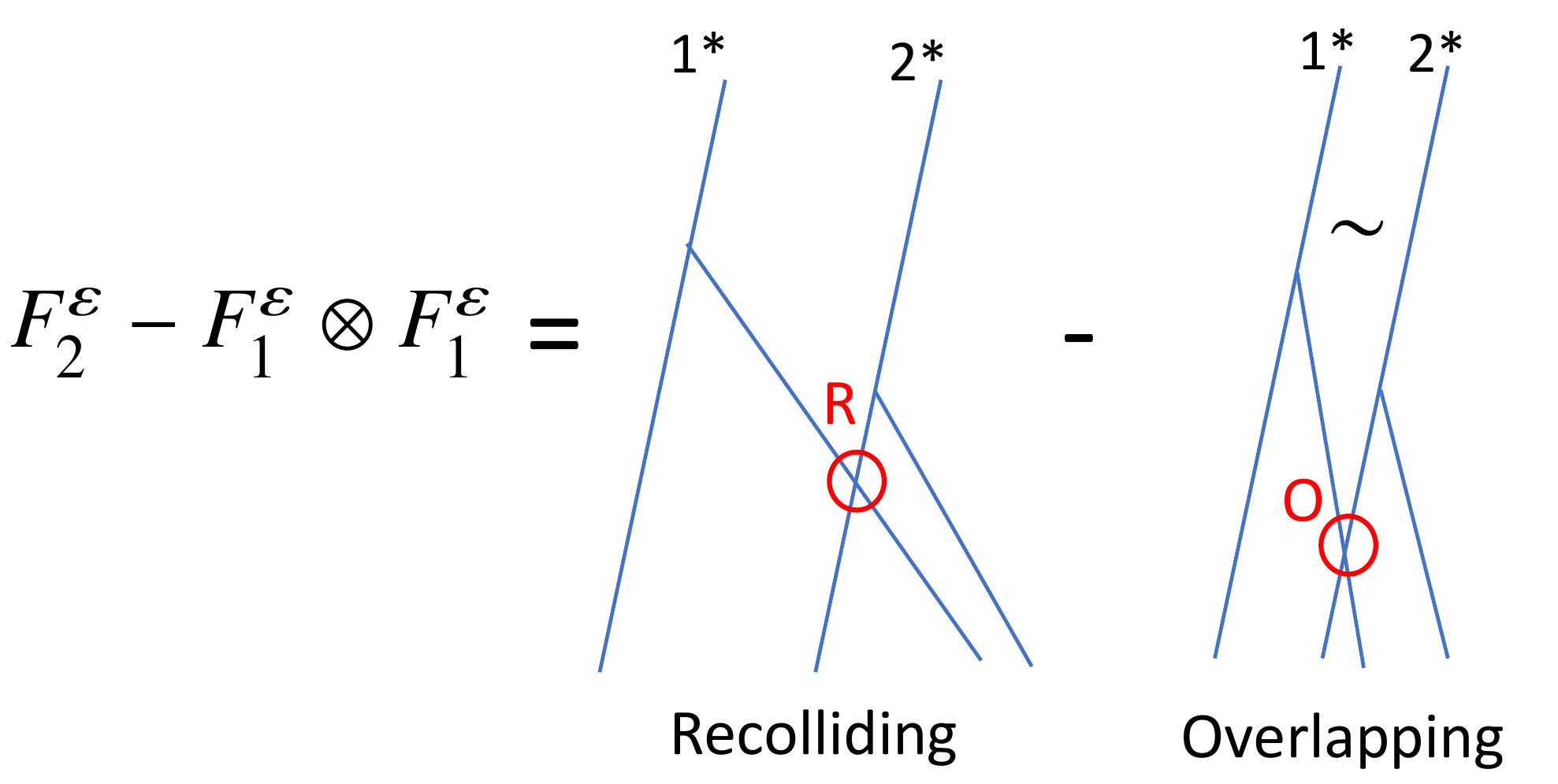} 
\caption{The second order cumulant corresponds to pseudo-trajectories with at least one external recollision or overlap.}
\label{fig: f2}
\end{figure}

Note that this decomposition  is a pure mathematical artefact to compare pseudo-dynamics without external recollision with independent pseudo-dynamics. In particular, the overlapping condition $1^* \sim 2^*$  does not affect the dynamics itself (overlapping particles are not scattered!). If we ignore the correlation encoded in the initial data, we then end up with a representation of the second order cumulant  by trees which are coupled  by  external recollisions or overlaps (see Figure \ref{fig: f2}).

\begin{Rmk} 
Recall that the initial measure does not factorise exactly $F_2^{\eps, 0} \neq F_1^{\eps, 0} \otimes F_1^{\eps, 0}$ due to the exclusion condition. Thus the initial data induces also a small correlation which is actually much smaller than the dynamical correlations (by a factor $\eps$), so we will neglect it in the following. 
\end{Rmk}

\medskip

Recolliding and overlapping pseudo-trajectories should provide a contribution of order 1 in $L^\infty$  to $F_2^\eps- F_1^\eps \otimes  F_1^\eps$. For $n_1 = n_2 = 0$, i.e. for collision trees without branchings, this defines the bad set of configurations $\cB^\eps_2$ (mentioned in Sections \ref{lanford-proof}-\ref{irreversibility})
encoding the collisions between  two particles in the backward flow (see Figure \ref{fig: bad sets}).
In particular, by choosing $z_{1^*}$ and $z_{2^*}$ at time $t$ such that 
$\big| x_{1^*} - x_{2^*} - (v_{1^*} - v_{2^*}) (t-s) \big| \leq \eps$ for some~$s \leq t$, the 
contribution to the cumulant of the pseudo-dynamics with $n_1 = n_2 = 0$ is expected to be non zero (except at equilibrium when recollisions and overlaps almost compensate).
Smallness of the second cumulant $F_2^\eps- F_1^\eps \otimes  F_1^\eps$ actually comes from the size of its support. The right norm to measure the smallness of correlations is thus the $L^1$ norm
and the quantity to be studied asymptotically is the rescaled second order cumulant
\begin{equation}
\label{eq: 2nd cumulant}
f_2^\eps = \mu_\eps (F_2^\eps- F_1^\eps \otimes  F_1^\eps)\,.
\end{equation}
With this scaling, we expect that $f_2^\eps$ has a limit $f_2$ in the sense of measures. 
The set supporting the function $f_2^\eps$ records the correlation between two pseudo-trajectories 
(rooted in $1^*$ and $2^*$) via a recollision or an overlap. 
On the other hand, once the two pseudo-trajectories are correlated by a recollision or an overlap then any additional recollision, overlap or internal recollision will impose stronger geometric constraints and they can be discarded in the limit as in Lanford's proof (see Figure \ref{recollision-fig}).
Therefore the limit $f_2$ corresponds to pseudo-trajectories with exactly one (external) recollision or 
overlap on $[0,t]$.

\medskip

\begin{figure}[h] 
\centering
\includegraphics[width=3.5in]{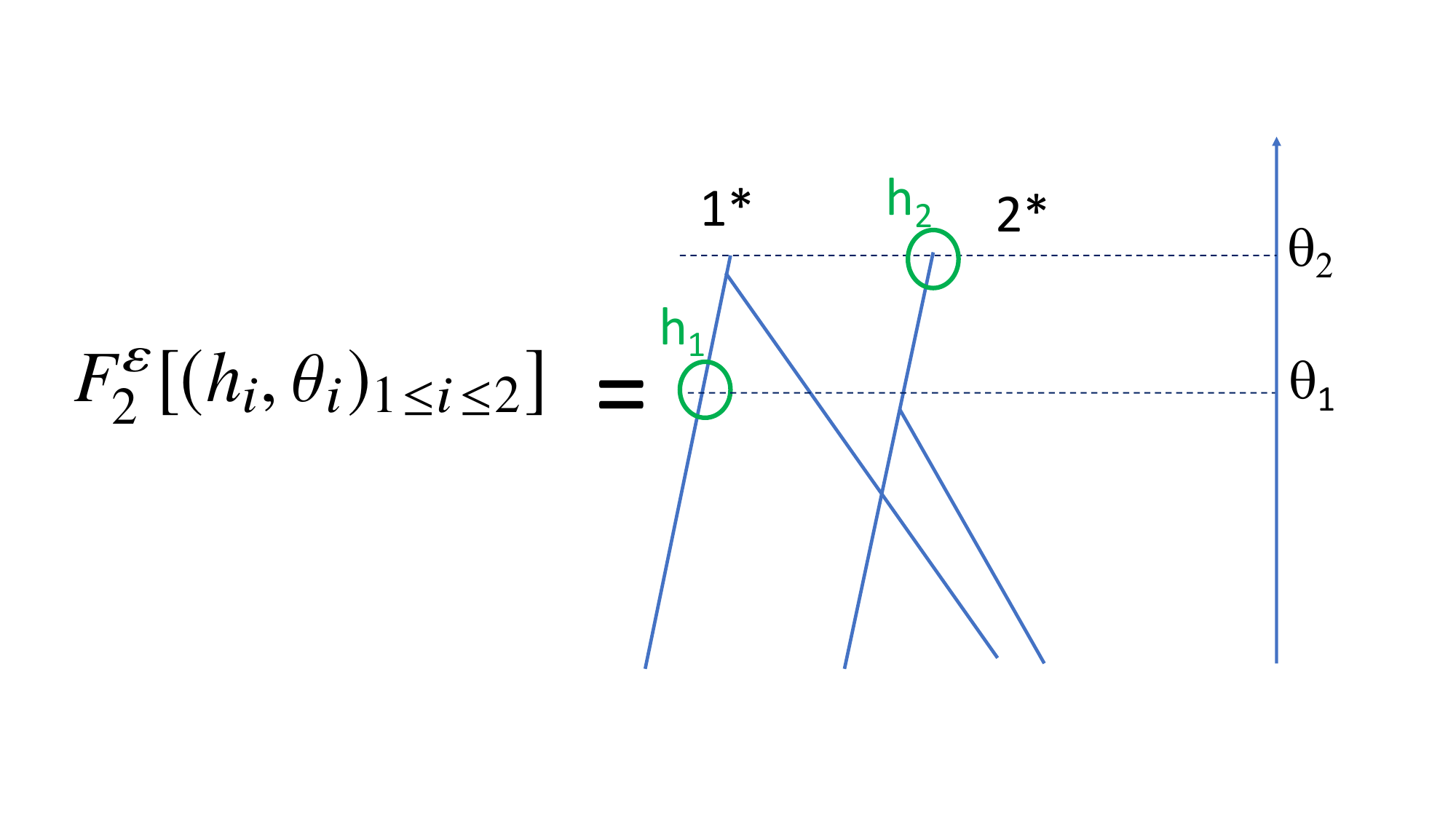} 
\caption{Time correlations \eqref{eq: time correlations} can be computed  by introducing weights along the pseudo-trajectories.}
  \label{weight-fig}\end{figure}

In order to understand fluctuations with respect to the Boltzmann dynamics, we also need to understand time correlations. To characterize these time correlations, one can proceed exactly in the same way, using a kind of duality method with weighted pseudo-trajectories. Recall that $F_2^\eps$ is by definition
$$ \int F_2^\eps (t, z_{1^*}, z_{2^*}) h_1(z_{1^*}) h_2(z_{2^*}) dz_{1^*} dz_{2^*} 
= \bbE_\eps \left( \frac1{\mu_\eps^2} \sum_{(i_1,i_2)} h_1\big( z^\e_{i_1}(t)\big) h_2\big( z^\e_{i_2}(t)\big) \right),
$$
meaning that there is a weight $h_1(z_{1^*}) h_2(z_{2^*})$ at time $t$ in the geometric representation. 
The counterpart for the time correlations 
\begin{equation}
\label{eq: time correlations}
 F^\eps_2[ (h_i, \theta_i)_{ i \leq 2}]    = 
 \bbE_\eps \left( \frac1{\mu_\eps^2} \sum_{(i_1,i_2)} h_1\big( z^\e_{i_1}(\theta_1)\big) h_2\big( z^\e_{i_2}(\theta_2)\big) \right)
\end{equation}
is to construct the same pseudo-trajectories $\Psi^\eps_{2,n}$ starting from some $\theta_2 > \theta_1$,  and to evaluate the weight $h_1$ on the resulting configuration of particle $1^*$ at time $\theta_1$ and the weight $h_2$ on the resulting  configuration of particle $2^*$ at time $\theta_2$ (see Figure \ref{weight-fig}).  
  
We then define the rescaled weighted second order cumulant
 \begin{equation}
 \label{weighted-cum}
 f^\eps_2[ (h_i, \theta_i)_{ i \leq 2}] = \mu_\eps \Big( F^\eps_2[ (h_i, \theta_i)_{ i \leq 2}]  - F^\eps_1[ h_1, \theta_1] F^\eps_1[ h_2, \theta_2] \Big) \,,
 \end{equation}
 and performing the same geometric analysis as before, the cumulant $f^\eps_2[ (h_i, \theta_i)_{ i \leq 2}]$ at different times converges also to a limit $f_2[ (h_i, \theta_i)_{ i \leq 2}]$ as $\mu_\eps$ diverges.

\subsection{Higher order correlations and exponential moments}
\label{cumulant-proof}

For a Gaussian process, the first two correlation functions $F_1^\eps$, $F_2^\eps$ determine completely all other correlation functions $F_k^\eps$, but in general part of the information is encoded in the (scaled) cumulants of higher order defined by (restricting here for simplicity to only one time)
\begin{equation*}
 f_k^\eps (t, Z_k)  = \mu_\eps^{k-1} \sum_{\ell = 1}^k 
\sum _{\sigma \in \cP^\ell_k} (-1)^{\ell-1} (\ell  -1)! \prod_{i = 1 } ^\ell   F^\eps_{|\sigma_i|} ( t, Z_{\sigma_i}) \,,
\end{equation*}
where $\cP^\ell_k$ is the set of partitions of $\{1, \dots, k\}$ in $\ell$ parts with $\sigma = \{ \sigma_1, \dots, \sigma_\ell \}$, $|\sigma_i|$ stands for the cardinality of the set $\sigma_i$ and $Z_{\sigma_i} = \left( z_j\right)_{j \in \sigma_i}$.
Each cumulant encodes finer and finer correlations. Contrary to correlation functions $F_k^\eps$, they do not duplicate the information which is already encoded at lower orders.

\begin{figure}[h] 
\centering
\includegraphics[width=3.5in]{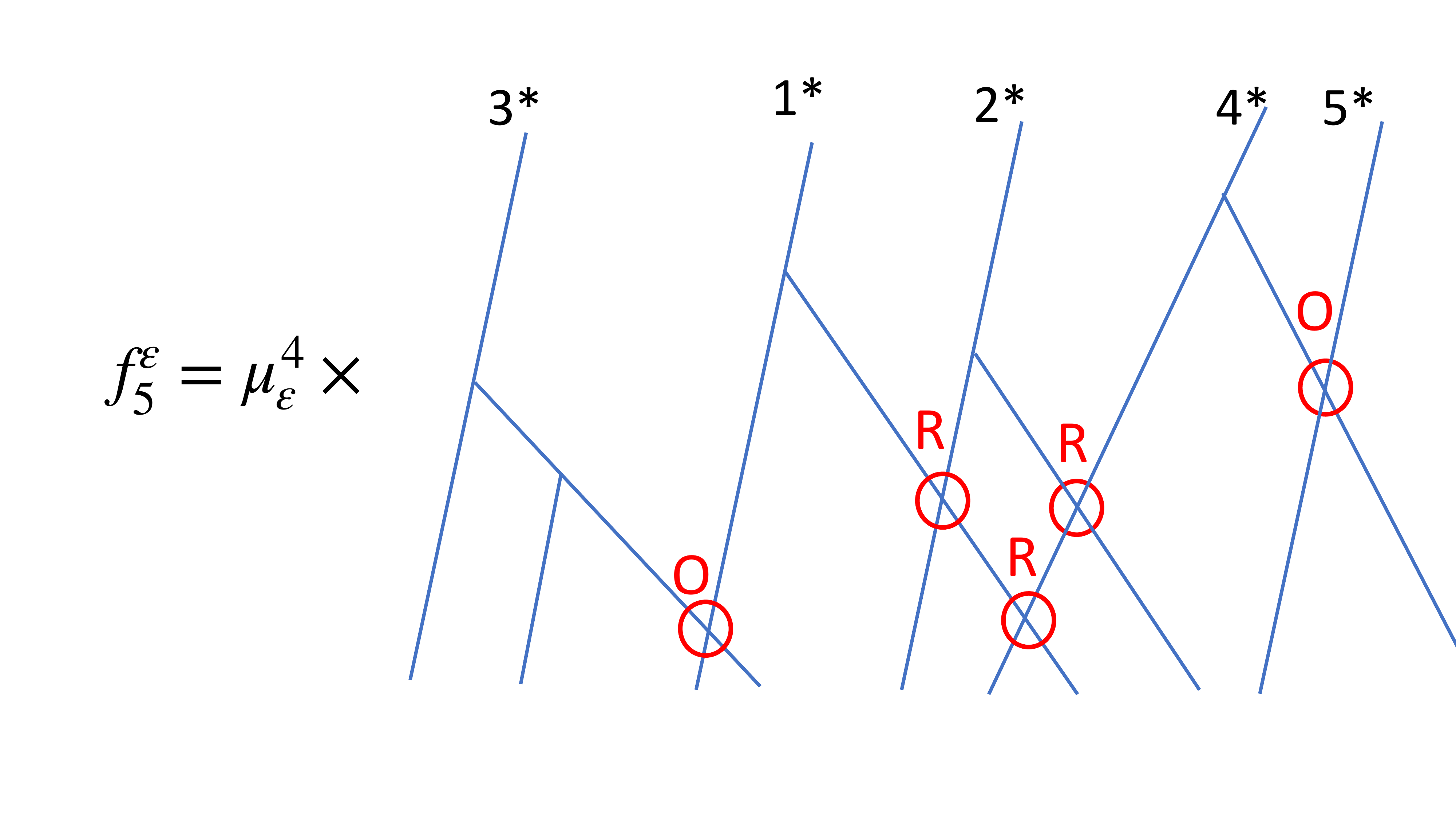} 
\caption{The cumulant of order $k$ corresponds to pseudo-trajectories issued from $z_{1^*}, \dots, z_{k^*}$ completely connected by external recollisions or overlaps.}
\label{cumulantk-fig}
\end{figure}

From the geometric point of view, one can extend the analysis of the previous paragraph and show that the cumulant of order $k$ can be represented by $k$ pseudo-trajectories which are completely connected either by external recollisions or by overlaps (see Figure~\ref{cumulantk-fig}).

\begin{figure}[h] 
\centering
\includegraphics[width=3.5in]{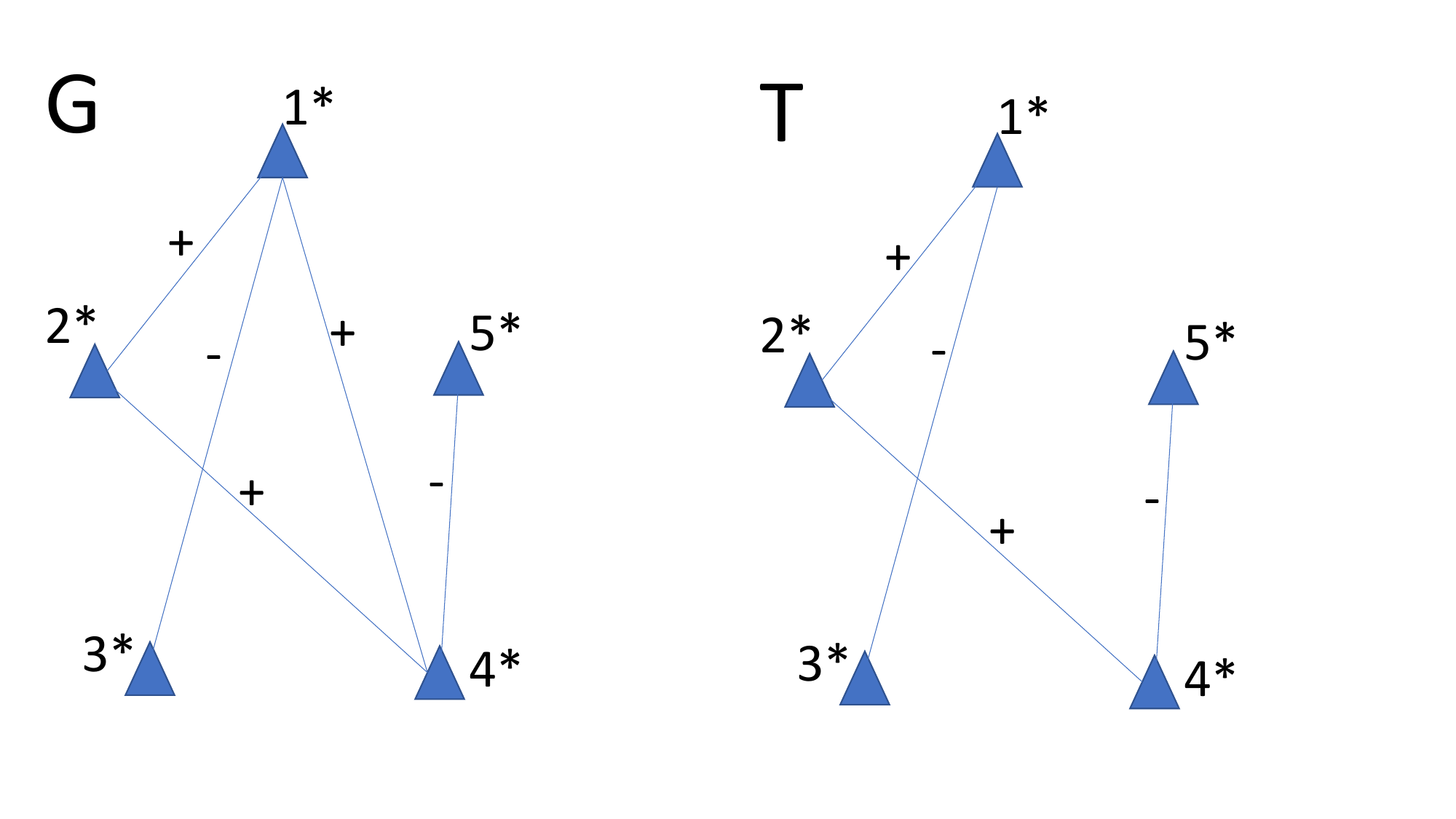} 
\caption{All recollisions and overlaps from the pseudo-trajectories depicted in Figure \ref{cumulantk-fig} are encoded  in the graph $G$. Only recollisions/overlaps which do not create a cycle (going backward in time) are kept in the tree $T$.}
\label{fig : graph}
\end{figure}

One can classify these completely connected pseudo-trajectories by associating them with a dynamical graph $G$ with $k$ vertices representing the different trees encoding the external recollisions (edge with a + sign) and the overlaps (edges with a - sign). 
Furthermore, one can define a systematic procedure to extract from this connected graph $G$ a minimally connected graph $T$ by identifying $k-1$ ``clustering recollisions'' or ``clustering overlaps'' (see Figure \ref{fig : graph}).
Here we use a {\color{black} cluster expansion reminiscent} of the method originally developed by Penrose to deal with correlations in the grand canonical Gibbs measure~\cite{penrose1967convergence, poghosyan2009abstract}.

We then expect the scaled cumulant $f_k^\eps$ to decompose in a sum of $2^{k-1}  k^{k-2}$ terms obtained by grouping all pseudo-trajectories compatible with each one of the signed minimally connected graphs $T$ (recall that $k^{k-2}$ is the number of trees on $k$ labelled vertices, known as Cayley's formula).
For each given signed minimally connected graph, the recollision/overlap conditions can be written as $k-1$ ``independent'' constraints on the configuration $z_{1^*}, \dots, z_{k^*}$ at time $t$. 
Therefore neglecting the velocity dependence as in \eqref{eq: size k-1 collisions}, 
this contribution to the cumulant $f_k^\eps$ has a support of size $O \left( ( t / \mu_\eps) ^{k-1} \right)$
with respect to Lebesgue measure and from this we deduce the expected $L^1$ estimate
\begin{equation}
\label{cum-est}
 \| f_k^\eps\|_{L^1} \leq 
 \underbrace{\mu_\eps ^{k-1}}_{\text{scaling}}  
\;  \times \;
 \underbrace{2^{k-1} k^{k-2}}_{\text{number of signed trees}}  
 \;  \times \;  
 \underbrace{\Big({C  t \over \mu_\eps}\Big) ^{k-1}}_{\text{support size}} 
 \leq k! \; (C  t)^{k-1}\,.
 \end{equation}
Furthermore, a  geometric argument similar to the one developed in Lanford's proof (see Section \ref{lanford-proof}) and already  used in the study of the second order cumulant allows to show that $f_k^\eps$ converges to some limiting cumulant $f_k$ and that only the pseudo-trajectories having exactly $k-1$ recollisions or overlaps (and no cycle) contribute in the limit. 

\medskip

This geometric approach allows to characterise \textit{all corrections to the chaos assumption, up to exponential order}, at least for times of the same order as $T_L$ \cite{BGSS1,BGSS2}.
Actually a classical and rather straightforward computation (based on the series expansions of the exponential and logarithm) shows that cumulants are nothing else than the coefficients of the series expansion of the exponential moment
\begin{align}
\label{exp-moment}
\cI^\eps_t (h) = &
\frac1{\mu_\eps} \log \bbE_\eps 
\Big[ \exp \big( \mu_\eps \la \pi^\eps_t, h\ra \big) \Big]  
=  \frac1{\mu_\eps} \log 
\bbE_\eps \left[ \exp \left(  \sum_i h \big( z^\eps_i(t) \big) \right) \right] \\
= &
\sum_{k = 1} ^\infty {1\over k!} \int f_k^\eps (t, Z_k)  \prod_{i = 1} ^k (e^{h(z_i) } - 1) dZ_k \,.
\nonumber
 \end{align}
The quantity ${\mathcal I}^\eps_t (h)$ is referred to as the \emph{cumulant generating function}.  
Estimate (\ref{cum-est}) provides the analyticity of $\cI^\eps_t (h)$ as a functional of  $e^h$, and this uniformly with respect to  $\eps$ (small enough). 
The limit $\cI_t$  of $\cI^\eps_t$ can then be determined as a series in terms of the limiting cumulants  $f_k$.

\medskip

Instead of using the cumulant expansion, we present a heuristic approach to characterise the limit $\cI_t$ as the solution of
{\color{black}
 the \emph{Hamilton-Jacobi equation} \eqref{HJ-eq}.
 At first reading, this formal derivation can be skipped and the reading can be resumed at Equation \eqref{HJ-eq}.}
We proceed as in Section \ref{lanford-proof-heuristics} for the Boltzmann equation \eqref{boltz}
and write the formal equation satisfied by $\cI^\eps_t (h)$ for fixed $\eps$. 
Considering an  evolution for a short time $\delta$ as in \eqref{eq: short time delta} and then taking a formal limit $\delta \to 0$, 
we get (after deriving the logarithm)
\begin{align*}
&
\bbE_\eps \left[ \exp \left(  \sum_i h \big( z^\eps_i(t) \big) \right) \right]  \partial_t \cI^\eps_t (h)
=  
\bbE_\eps \left[ \left( \frac{1}{\mu_\eps} \sum_j {d x_j^\eps \over dt}  \cdot \nabla_x h \big( z^\eps_j(t) \big)\right) \exp \left( \sum_i h \big( z^\eps_i(t) \big) \right) \right] \\
&  \qquad \qquad 
+ 
\; \int d \omega \; \bbE_\eps \left[   \frac{1}{\mu_\eps^2} \sum_{j_1 \neq j_2}  
 \delta_{ x^\eps_{j_2} (t)- x^\eps_{j_1}(t) - \eps\omega } 
\big( (v^\eps_{j_2} (t) - v^\eps_{j_1} (t) ) \cdot \omega \big)_+
\right.\\
& \qquad  \qquad \qquad \qquad \left.  \times \left( e^{h ( z^\eps_{j_1} (t^+)) +h ( z_{j_2}^\eps (t^+) )} -e^{  h ( z_{j_1}^\eps (t^-)) +h ( z_{j_2}^\eps (t^-) ) } \right)
 \exp \left( \sum_{i \not = j_1,j_2} h \big( z^\eps_i(t) \big) \right) \right] ,
\end{align*}
{\color{black}
where $\omega$ becomes a random parameter after changing variables at the collision time  as in~\eqref{eq: jacobien}. We used the Dirac notation as in \eqref{eq: weak boltzmann} to stress that $x^\eps_{j_2} (t) = x^\eps_{j_1}(t) + \eps \omega$ at the collision.}
Denoting by $\pi_{2,t}^\eps$ the generalised empirical measure depending on $2$ arguments (see \eqref{eq: Fk}), we get
\begin{equation}
\label{time-derivative}
\begin{aligned}
 \bbE_\eps & \left[ \exp \left(  \sum_i h \big( z^\eps_i(t) \big) \right) \right]    \partial_t \cI^\eps_t (h)
=\bbE_\eps \left[ \pi_t^\eps \{v\cdot \nabla_x h\} \exp \left( \mu_\eps \la \pi^\eps_t, h\ra   \right) \right] \\
& \qquad \qquad \qquad 
+ \frac{1}2 
{\color{black} \int d \omega} \bbE_\eps \left[ \pi_{2,t}^\eps \left\{ \delta_{x_2-x_1 - \eps\omega  }
\left( e^{ \Delta h (z_1,z_2, \omega) } - 1 \right) \right\}  
\exp \left( \mu_\eps \la \pi^\eps_t, h\ra    \right) \right]  ,
\end{aligned}
\end{equation}
where $\Delta h (z_1,z_2, \omega) = h \big( x_1, v_1'  \big) + h \big( x_2, v_2' \big) 
- h \big( z_1 \big) - h \big( z_2 \big)$ was already introduced in \eqref{eq: omega collision}.
To obtain a closed equation, it remains to find the counterparts of the correlation functions
$F^\eps_1$ and $F^\eps_2$ which describe the distribution under the measure tilted  by the exponential weight $\la \pi^\eps_t, h\ra$.

Differentiating the exponential moment \eqref{exp-moment} at $h$ in the direction $\gp$, we 
recover the quantity $\la \pi_t^\eps , \gp \ra$
\begin{align*}
\la {\d \mathcal{I}^\eps_t \over \d h}  (h) , \gp \ra 
= &
\lim_{\delta \to 0} \frac{1}{\delta} (  \mathcal{I}^\eps_t (h + \delta \gp) -  \mathcal{I}^\eps_t (h))\\
= & \frac{1}{ \bbE_\eps \left[ \exp \left( \mu_\eps \la \pi^\eps_t, h \ra  \right) \right]}
 \bbE_\eps 
\left[  \la \pi_t^\eps , \gp \ra \; \exp \left( \mu_\eps \la \pi^\eps_t, h\ra  \right) \right]  .
\end{align*}
Thus the transport term has the form $\big\la {\d  \mathcal{I}^\eps_t \over \d h} (h) , v\cdot \nabla_x h\big\ra$.  By taking a second derivative, the 
tilted distribution of the two-point correlations can be identified in terms of  
$$
\frac{1}{\mu_\eps}   {\d^2 \mathcal{I}^\eps_t \over \d h^2}  (h) 
+ {\d \mathcal{I}^\eps_t \over \d h}   (h)  \otimes  {\d \mathcal{I}^\eps_t \over \d h} (h).
$$
The collision term is singular, but formally   the right-hand side of (\ref{time-derivative}) 
can be rewritten as
\begin{equation}
\label{HJeps-eq}
\begin{aligned}
\d_t \mathcal{I}^\eps_t (h) 
=&  \frac{1}2 
 \Big\la {\d  \mathcal{I}^\eps_t \over \d h}   (h)  \otimes  {\d  \mathcal{I}^\eps_t \over \d h} (h) , 
 {\color{black} \int d \omega }
 ((v_2-v_1)\cdot \omega)_+ \delta_{x_2-x_1 - \eps \omega}
\Big( e^{ \Delta h (z_1,z_2, \omega)}-1 \Big) \Big \ra\\
&+
 \frac{1}{2 \mu_\eps}  \Big\la  {\d^2  \mathcal{I}^\eps_t \over \d h^2}  (h) , 
  {\color{black} \int d \omega} \, ((v_2-v_1)\cdot \omega)_+ \delta_{x_2-x_1 - \eps \omega}
\Big( e^{ \Delta h (z_1,z_2, \omega) }-1 \Big) \Big\ra\\
&+ \Big\la {\d  \mathcal{I}^\eps_t \over \d h}  (h) , v\cdot \nabla_x h\Big\ra .
\end{aligned}
\end{equation}
 We recognize here a kind of Hamilton-Jacobi equation, with a small  ``viscous" term (involving derivatives of order 2 with respect to $h$, but without a definite sign). 
Thus the limiting functional $\mathcal{I}_t$ has  to satisfy the following Hamilton-Jacobi equation obtained
by formally taking the limit  $\mu_\eps \to \infty$
 \begin{equation}
 \label{HJ-eq}
 \begin{aligned}
\d_t \mathcal{I}_t (h) = & \frac{1}2  \Big\la {\d \mathcal{I}_t \over \d h} (h)  \otimes
  {\d \mathcal{I}_t \over \d h} (h), \, \int d \omega \, ((v_2-v_1)\cdot \omega)_+ \delta_{x_2-x_1}
\left( e^{ \Delta h (z_1,z_2, \omega) }-1 \right) \Big \ra\\
&+ \Big\la {\d\over \d h} \mathcal{I}_t (h), v\cdot \nabla_x h\Big\ra .
\end{aligned}
\end{equation}
The structure of this Hamilton-Jacobi equation is reminiscent of the Boltzmann equation~\eqref{HJ-eq}, with a collision term and a transport term. However it encodes a much more complete description of the hard sphere dynamics, including in particular  the structure of the exponentially small correlations and of the large deviations (see Theorem \ref{thm2}).

As in \eqref{eq: time correlations},  further information on the correlations in a time interval $[0,t]$ can be obtained by generalising \eqref{exp-moment}
\begin{align}
\label{eq: exponential moment interval}
\cI^\eps_{[0,t]} (H) = \frac1{\mu_\eps} \log  \bbE_\eps \left[ \exp \left(  \sum_i H \Big( z^\e_i ([0,t]) \Big) \right) \right] ,
\end{align}
for functions $H$ depending on the trajectory of a particle in $[0,t]$.
For example, a sampling at different times $\theta_1 < \theta_2 < \dots < \theta_k \leq t$ by test functions 
$( h_\ell)_{\ell \leq k}$ is obtained by considering 
\begin{equation}
\label{eq: time sampling}
H \Big( z([0,t]) \Big) = \sum_{\ell =1}^k h_\ell \big( z(\theta_\ell) \big).
\end{equation}

\begin{Rmk} 
The  procedure described here allows to obtain easily the limiting equation~(\ref{HJ-eq}) without having to guess how to combine the different cumulant terms (which happens to be quite technical). However the weak understanding we have on this equation does not allow to use it to justify the limit as $\mu_\eps \to \infty$ (without going through the cumulant analysis
of \cite{BGSS2}).
\end{Rmk}

\begin{Rmk} 
In the  absence of spatial inhomogeneities, one can discard the transport term and retrieve asymptotically the same cumulant generating function as for the Kac model, i.e. the dynamics in which collisions are given by a random jump process  \cite{Leonard, heydecker2021large, Rezakhanlou_GD, BBBO}. This indicates that in the limit $\mu_\eps \to \infty$, both models are indistinguishable (up to exponentially small corrections). In other words, the Hamilton Jacobi equation \eqref{HJ-eq} conserves the stochastic reversibility, but not the deterministic reversibility~: one cannot hope for any strong convergence result.
\end{Rmk}

\subsection{A complete statistical picture for short times}
\label{BGSS-sec}

As mentioned in the previous paragraph, the cumulant generating function provides a complete statistical picture of the hard sphere dynamics. We  now explain how it can be used to answer the main questions raised in Section \ref{questions-sec} (on a short time $T^\star$, of the same order as  Lanford's time $T_L$  in Theorem \ref{lanford-thm}).

\medskip

As a first consequence of  the uniform estimates on the cumulant generating function~${\mathcal I}^\eps_{[0,t]}$, the convergence of the fluctuation field, defined by (\ref{zeta-def}) and recalled below 
$$ 
\la \zeta^\eps_t ,h \ra = \sqrt{\mu_\eps}  \Big(  \la \pi^\eps_t, h\ra   - \bbE_\eps ( \la \pi^\eps_t, h\ra) \Big),
$$
can be obtained.

At time $0$, it is known 
that, under the grand-canonical measure introduced Page \pageref{Def: initial data}, the fluctuation field $\zeta^\eps_0$ converges in the Boltzmann-Grad limit to a Gaussian 
field $\zeta_0$ with covariance
\begin{equation}
\label{eq: Gaussian field time 0}
\bbE \Big ( \zeta_0 (h) \; \zeta_0(g) \Big) = 
\int dz \; f^0 (z) h(z) g(z).
\end{equation}
The following theorem controls the dynamical fluctuations.

\begin{Thm}[Bodineau, Gallagher, Saint-Raymond, Simonella \cite{BGSS3}]
\label{BGSS-thm1}
Under the assumptions on the initial data stated Page \pageref{Def: initial data},
 the fluctuation field $\zeta^\eps_t $ of the hard sphere system converges,  in the Boltzmann-Grad limit ($\mu_\eps \to \infty$ with $\mu_\eps \eps^{d-1} = 1$), on a time interval $[0,T^\star]$ towards a process $\zeta_t$, solution to the fluctuating Boltzmann equation~:
\begin{equation}
\label{boltz-fluct Heq}
\left\{ 
\begin{aligned}
&d\zeta_t =  \underbrace{\cL_t \zeta_t  dt}_{\mbox{\footnotesize{linearized Boltzmann operator} }} + \underbrace{d\eta_t }_{\mbox{\footnotesize{Gaussian noise} }} \\
&\cL_t h  =  \underbrace{- v\cdot \nabla_x h  }_{\mbox{\footnotesize{transport} }} + \underbrace{C (f_t, h) + C(h,f_t)    }_{\mbox{\footnotesize{ linearized collision operator} }}
\end{aligned}\right.
\end{equation}
where $f_t$ denotes the solution at time $t$ to the Boltzmann equation (\ref{boltz}) with initial data $f^0$, and $d \eta_t$ is  a centered Gaussian noise delta-correlated in $t,x$  with covariance
$$
\Cov_t (h_1, h_2) =  \frac12 \int  dz_1 dz_2 d\omega 
 \; ((v_2-v_1)\cdot \omega)_+ \delta_{x_2-x_1}
 f(t,z_1) f(t,z_2) \Delta h_1\Delta h_2 (z_1, z_2, \omega)
$$
with $\Delta h (z_1, z_2, \omega ) =   h (z'_1) + h (z'_2) - h(z_1) - h( z_2)$
as in \eqref{eq: omega collision}.
\end{Thm}

As hinted in Section \ref{sec: Defects in the chaos assumption}, the limiting noise is a consequence of the asymptotically unstable structure of the microscopic dynamics (see Figure \ref{fig: instability}).
The randomness of the initial configuration is transported deterministically by the dynamics and generates a white noise in space and time through a particular class of collisions. The velocity scattering mechanism is coded in the covariance of the noise.

If the system starts initially from an equilibrium measure, i.e. with particle positions spatially  independent (up to the exclusion) and velocities  identically distributed according to the Maxwell-Boltzmann equilibrium distribution
\begin{equation}
\label{eq: Maxwell}
f^0 (x,v) = M (v) = \frac{1}{(2 \pi)^{d/2}} \exp \left( -  \frac{|v|^2}{2} \right),
\end{equation} 
  then $f_t = f^0$ so that the linearised operator is time independent and it will be denoted by $\cL_{eq}$. 
The limiting stochastic partial differential equation $d\zeta_t =  \cL_{eq} \, \zeta_t + d\eta_t$ satisfies the fluctuation/dissipation relation : the dissipation from the linearised operator $\cL_{eq}$ is exactly compensated by the noise 
$\eta_t$. As the equilibrium measure is time invariant, it was expected on physical grounds 
 that a stochastic correction should emerge in order to keep this invariance in time. 
In fact, the equation governing the covariance of the limiting process $\Cov( \zeta_t)$ away from equilibrium was obtained, and the full fluctuating equation  for $(\zeta_t)_{t \geq 0}$  conjectured, in the pioneering works by Spohn \cite{Spohn_fluctuations, Spohn_review_fluctuations, spohn2012large}. 
In particular, it was already understood in  \cite{Spohn_fluctuations} that out of equilibrium, a non-trivial contribution to $\Cov( \zeta_t)$ is provided   by the second order cumulant \eqref{eq: 2nd cumulant}.
Note that the predictions on the stochastic corrections from the Kac model 
\cite{logan1976fluctuations, Meleard_fluctuations, Rezakhanlou_fluctuations} fully agree with the stochastic equation emerging from the deterministic hard sphere dynamics. 
{\color{black}
Thus from a phenomenological point of view, it is equivalent to consider a stochastic model (including as well the positions as in \cite{Rezakhanlou_fluctuations}) or a deterministic evolution.}
We refer also to the work by Ernst and Cohen \cite{ernst1981nonequilibrium} for further discussion on the time correlations and the fluctuations.

Note that equilibrium fluctuations for a microscopic evolution with spatial coordinates and  stochastic collisions have been derived in \cite{Rezakhanlou_fluctuations} for arbitrary long times.
We will see in Theorem \ref{thm3} that the convergence time of the previous theorem 
can be greatly improved at equilibrium.

\medskip

Out of equilibrium, although the solution $f$ to the Boltzmann equation (describing the averaged dynamics) is very smooth on~$[0,T^\star]$,  the fluctuating Boltzmann equation is quite singular~: the linearized operator $\cL_t$  is non autonomous, non self-adjoint, and the corresponding semigroup is not a contraction. 
Thus we consider a very weak notion of solution of \eqref{boltz-fluct Heq}, requiring only that
\begin{itemize}
\item the process $\zeta_t$ is Gaussian;
\item its covariance defined, for  test functions $h_1,h_2$ and times $\theta_1, \theta_2$,  as
{\color{black}  
\begin{equation}
\label{eq: covariance}
\cC ( \theta_1, h_1, \theta_2, h_2)
= \lim_{\eps \to 0} \bbE_\eps \left( \la \zeta^\eps_{\theta_1} , h_1 \ra \, \la \zeta^\eps_{\theta_2} , h_2 \ra \right)
\end{equation}
satisfies a set of equations governed by the linearised Boltzmann equation.}
\end{itemize}
The convergence of the process $( \zeta^\eps_t )_{t \leq T^\star}$ 
can be derived in 3 steps :
\begin{itemize}
\item \textit{The convergence of the time marginals to a Gaussian process.}

The characteristic function of the process tested at times 
$\theta_1 <  \dots < \theta_k \leq T^\star$ by  functions $(h_\ell)_{\ell \leq k}$ is encoded by the exponential moment \eqref{eq: exponential moment interval} by choosing~$H \Big( z([0, T^\star ]) \Big) = \frac{i}{\sqrt{\mu_\eps}} \sum_{\ell =1}^k h_\ell \big( z(\theta_\ell) \big)$ as in 
\eqref{eq: time sampling}
\begin{align}
\label{eq: fonction caracteristique}
\log \bbE_\eps \left[ \exp \left( i \sum_{\ell =1}^k \left( \la \zeta^\eps_{\theta_\ell}, h_\ell \ra
+\sqrt{\mu_\eps} \bbE_\eps \left(  \la \pi^\eps_{\theta_\ell}, h_\ell \ra \right)\right)\right)
 \right] 
= \mu_\eps \; \cI^\eps_{[0, T^\star ]} (H) \,.
\end{align}
The cumulant expansion \eqref{exp-moment} combined with sharp controls on the cumulants ensure that $\cI^\eps_{[0,T^\star]} (H)$ is an analytic function of $H$ in a neighbourhood of 0 so that complex values can also be handled. Furthermore in the scaling considered for the fluctuations, $H$ is of order $\frac{1}{\sqrt{\mu_\eps}}$. Thus in the cumulant expansion \eqref{exp-moment},
the term of order $n$ scales as 
$$
f^\eps_n \Big( (e^{H} - 1)^{\otimes n}  \Big)   \simeq \frac{1}{\mu_\eps^{n/2} },
$$
so that the asymptotics of the characteristic function \eqref{eq: fonction caracteristique}
is only determined by the cumulants of order less than $2$.
This implies that the Wick rule holds and therefore the limiting variables are Gaussian.
\item \textit{The characterisation of the limit covariance}.

The evolution equation of the covariance $\cC ( \theta_1, h_1, \theta_2, h_2)$ can be recovered from the equations satisfied by the first two cumulants.
As already pointed out in \cite{Spohn_fluctuations}, we stress that the behaviour of the covariance $\cC ( \theta_1, h_1, \theta_2, h_2)$ is determined by means of a careful analysis of the second cumulant
$f^\eps_2[ (h_\ell, \theta_\ell)_{ \ell \leq 2}]$ introduced in \eqref{weighted-cum}. 
Out of equilibrium, the cumulant of order 2 takes into account the contribution of one external recollision 
or of one overlap (as explained in Section \ref{sec: Defects in the chaos assumption}).
Even though the contribution of the recollisions vanishes when deriving the Boltzmann equation (recall the chaos assumption \eqref{eq: molecular chaos}), it plays an important role in the stochastic corrections. 

\item \textit{The tightness of the sequence $(\zeta^\eps_t)_{\eps >0}$}. 

This is the most technical part of the proof as it requires to control uniform estimates in time for a wide class of test functions $h$
$$
\bbE_\eps  \big[  \sup_{|s-s'|\leq \delta } | \la \zeta^\eps_s , h \ra  - \la \zeta^\eps_{s'} , h \ra | \big].
$$
We will not discuss further this point and refer to \cite{BGSS2} for details.
\end{itemize}
Note that Theorem \ref{BGSS-thm1}, which is a kind of central limit theorem, does not use the fine structure of cumulants~:  a sufficient decay of the correlations is enough to control the typical fluctuations (which are of size $O(1/\sqrt{\mu_\eps})$).

\medskip

The strength of the cumulant generating function appears at the level of large deviations, 
i.e. for very unlikely trajectories which are at a ``distance'' $O(1)$ from the averaged dynamics.
The counterpart of the large deviation statement \eqref{eq: LD iid} for independent variables can be rephrased, in a loose way, as follows  : observing an empirical particle distribution  close to 
the density $\gp (t,x,v)$ during the time interval $[0,T^\star]$ decays exponentially fast with a rate quantified by the large deviation functional $\cF$
\begin{equation*}
\bbP_\eps \left( \pi^\eps_t \simeq \gp_t, \quad \forall t \leq T^\star  \right)
\sim \exp \big( - \mu_\eps \cF(\varphi)  \big) \,.
\end{equation*}
Notice that at time $0$, under the grand-canonical measure introduced Page \pageref{Def: initial data}, it is known that the large deviations around a density $\varphi^0$
can be informally stated as follows
\begin{equation*}
\bbP_\eps \left( \pi^\eps_0 \simeq \gp^0   \right)
\sim \exp \big( - \mu_\eps H(\varphi^0| f^0)  \big) \,,
\end{equation*}
with a static large deviation functional given by the relative entropy 
$$
H(\varphi^0| f^0) = \int \Big( \varphi^0 \log \frac{\varphi^0}{f^0} - (\varphi^0- f^0) \Big) dz. 
$$ 

More precisely, the distance between $\pi^\eps$ and $\gp$ is measured with respect to 
a weak topology on the Skorokhod space of measure valued functions. 
This topology is used in the theorem below.

\begin{Thm}[Bodineau, Gallagher, Saint-Raymond, Simonella \cite{BGSS2}]
\label{thm2}  
Under the assumptions on the initial data stated Page \pageref{Def: initial data},
there is a time $T^\star>0$ such that 
the empirical measure~$(\pi^\eps_t)_{t \leq T^\star}$ satisfies, in the Boltzmann-Grad limit $\mu_\eps \to \infty$ ($\mu_\eps \eps^{d-1} = 1$), 
 the following large deviation estimates
$$
\left\{ 
\begin{aligned}
&\limsup_{\mu_\eps \to \infty}  \; \frac{1}{\mu_\eps} \log \bbP_\eps \Big[ \pi^\eps \in K \hbox{ compact } \Big]  
\leq - \inf _{\varphi \in K}  \cF (\varphi) \\
&\liminf_{\mu_\eps \to \infty} \; \frac{1}{\mu_\eps} \log \bbP_\eps \Big[ \pi^\eps \in O \hbox{ open } \Big]  
\geq  - \inf _{\varphi \in O \cap \cR } \cF(\varphi) 
\end{aligned} \right.
$$
for some (nontrivial) restricted set $\cR$.  

The large deviation functional $\cF$ is defined by convex duality  from the cumulant generating function $\cI_{[0,T^\star]}$ (obtained as the limit of \eqref{eq: exponential moment interval}).
It  coincides on the restricted set~$\cR$ with 
\begin{equation}
\label{large-dev}
\left\{ \begin{aligned}
 \widetilde \cF(\varphi)  &=  \underbrace{H(\varphi^0| f^0)  }_{\mbox{\footnotesize{relative entropy of the initial data} }} 
 + \underbrace{ \sup_p \int _0^{T^\star} \left(\la p,( \d_t + v\cdot \nabla_x) \varphi \ra -  \cH (\varphi, p )\right)  }_{\mbox{\footnotesize{Legendre transform of the Hamiltonian} }} \\
 \cH (\varphi, p )
 & = \frac12 \int  dz_1 dz_2 \, d\omega \, 
 ((v_2-v_1)\cdot \omega)_+ \delta_{x_2-x_1}\varphi(z_1) \varphi(z_2) ( e^{\Delta p(z_1, z_2, \omega) } - 1) 
 \end{aligned}\right.
 \end{equation}
 with $\Delta p (z_1, z_2, \omega ) =   p (z'_1) + p (z'_2) - p(z_1) - p( z_2)$
as in \eqref{eq: omega collision}.
\end{Thm}

All the functionals appearing in the above statement are quite singular (notice that the Hamiltonian is defined by an integral over a manifold of codimension $d$ with a weight growing  for large velocities)
and our method is restricted to considering  very smooth and sufficiently decaying test functions. These restrictions on the functional spaces are the reason why we are not able to obtain a more precise large deviation principle, nor to identify clearly the large deviation functional.
We refer to \cite{BGSS2} for the proof which follows a quite standard path, once the limiting cumulant generating function $\cI_{[0,T^\star]}$ has been constructed.
The identification between $\cF$ and~$\widetilde \cF$ relies on the limiting Hamilton-Jacobi equation \eqref{HJ-eq}.

\begin{Rmk}
Note that the large deviation functional $\tilde \cF$ defined by (\ref{large-dev}) was conjectured in \cite{RezakhanlouLNM} and \cite{bouchet}. As already mentioned, it actually corresponds to the large deviation functional for stochastic microscopic processes, such as the Kac model (in the absence of transport) \cite{Leonard, heydecker2021large}, or intermediate models  (with transport and stochastic collisions) introduced by Rezakhanlou \cite{Rezakhanlou_GD}.
\end{Rmk}



\section{BEYOND LANFORD'S TIME}
\label{time-sec}

Up to a short time, Theorems \ref{BGSS-thm1} and \ref{thm2} provide a good statistical description of the hard sphere dynamics  in the Boltzmann-Grad limit 
($\mu_\eps \to \infty$ with $\mu_\eps \eps^{d-1} = 1$). 
The stochastic corrections to the Boltzmann equation emerge from the complex interplay between the random initial data and the asymptotic instability of the dynamics.

However, these results are still far from being satisfactory as the time restriction is not expected from  physics~: it does not allow to understand the relaxation toward equilibrium (and the corresponding entropy cascades between cumulants), nor to derive fluid limits. This question remains quite open, and the goal of this last section is to discuss theoretical obstructions and methodological difficulties, as well as some recent progress close to equilibrium.

\subsection{Main difficulties} \label{sec:maindifficulties}

A natural way to address this problem is trying to understand what kind of convergence one can hope for beyond Lanford's time $T_L$. Recall that Lanford's theorem describes the approximation of a reversible system by an irreversible system, where a macroscopic part of the information is missing.
This excludes any kind of ``strong'' convergence in terms of relative entropy.
This implies in particular that one will hardly  use the fine knowledge one might have on the solution to the Boltzmann equation to obtain a robust notion of stability 
which would be as well compatible with the microscopic system.

\begin{Rmk}
In the framework of fluid limits, these types of methods, referred to as modulated energy or modulated entropy methods, are among the most powerful  to prove convergence theorems \cite{Yau, Golse_Levermore_Saint-Raymond, SR}
since they require very few properties on the original system, typically
\begin{itemize}
\item an energy/entropy inequality satisfied by weak solutions;
\item the consistency of the approximation (meaning that the limiting equations are the ones inferred from the  formal asymptotics);
\item some bootstrap estimates controlling (nonlinear) fluxes in terms of the modulated energy/entropy.
\end{itemize}
\end{Rmk}

\medskip
An alternative would be to establish some weak convergence $F^\eps_1 \rightharpoonup f$, which paradoxically requires better compactness estimates on the sequence $(F^\eps_1)_\eps$. In this framework, the best one can do in general is to retrieve the structure of the limiting equation and its good (weak) stability properties from the solutions $F^\eps_1$ for fixed $\eps$, and this uniformly in $\eps$. The problem here, as mentioned in Section \ref{lanford-sec}, is that the Boltzmann equation does not have such a weak stability. Two ingredients are necessary to construct solutions satisfying only physical bounds (mass, energy and entropy estimates)~:
\begin{itemize}
\item a renormalization procedure to tame the possible singularity {\color{black} (concentration in $x$)}
in the loss collision term $   f(t,x,v)  \times \int f(t,x,v_1) b(v-v_1, \omega) d\omega dv_1$;
\item a bound on the entropy dissipation to control the gain term by the loss term.
\end{itemize}
These ingredients have been used in \cite{Rezakhanlou_BG} to recover the Boltzmann equation from a microscopic dynamics with stochastic collisions, but they do not seem to have a clear counterpart for a deterministic microscopic evolution.

\medskip

{\color{black}

The Hamilton Jacobi equation (\ref{HJeps-eq}) retains much more information on the system, thus 
the convergence of $\cI^\eps_t$ to $\cI_t$, in a sense to be understood, could provide a more stable
framework to study the kinetic limit for large times.
This would then imply the convergence to the Boltzmann equation.
}

%

\subsection{Close to equilibrium}
\label{sec: Close to equilibrium}

{\color{black} An easier setting to control the long time evolution is to consider a perturbation of an equilibirum measure. }Here the stationarity of the equilibrium becomes a key tool in order to provide uniform estimates in time and to control the pathological behaviours previously mentioned. In a series of recent works \cite{BGSS3, BGSS4}, we took advantage of the equilibrium structure to extend Theorem \ref{BGSS-thm1} to arbitrarily long kinetic times, and even slowly diffusive times.

\begin{Thm}[Bodineau, Gallagher, Saint-Raymond, Simonella \cite{BGSS3,BGSS4}] 
\label{thm3}
Consider a system of hard spheres initially at equilibrium, i.e. with a spatially uniform distribution and with a Maxwell-Boltzmann distribution $M$ in velocities as in  \eqref{eq: Maxwell} (Gibbs grand-canonical ensemble, $f^0 = M$ in \eqref{eq: initial measure}).

Then, in the Boltzmann-Grad limit $\mu_\eps \to \infty$ ($\mu_\eps \eps^{d-1} = 1$), the fluctuation field~$(\zeta^\eps_t)_{t \geq 0}$ of the hard sphere system converges on any time interval $[0,T_\eps]$, with $T_\eps = O(\log\log \log \mu_\eps)$, towards the process $(\zeta_t)_{t \geq 0}$, solution to the fluctuating Boltzmann equation~:
\begin{equation}
\label{boltz-fluct}
\left\{ 
\begin{aligned}
&d\zeta_t =  \underbrace{\cL_{eq} \,  \zeta_t dt  }_{\mbox{\footnotesize{linearized Boltzmann operator} }} + \underbrace{d\eta_t }_{\mbox{\footnotesize{Gaussian noise} }} \\
&\cL_{eq} h  =  \underbrace{- v\cdot \nabla_x h  }_{\mbox{\footnotesize{transport} }} + \underbrace{C( h , M) +C(M, h )    }_{\mbox{\footnotesize{ linearized collision operator} }}
\end{aligned}\right.
\end{equation}
where the linearised  operator $\cL_{eq}$ is   time independent and $\eta$ is  a Gaussian noise delta-correlated in $t,x$  with a time independent covariance
$$
\Cov (h_1, h_2) =  \frac12 \int  dz_1 dz_2 d\omega ((v_2-v_1)\cdot \omega)_+ \delta_{x_2-x_1}M(v_1)M(v_2) \Delta h_1\Delta h_2 (z_1, z_2, \omega)\,,
$$
with $\Delta h (z_1, z_2, \omega ) =   h (z'_1) + h (z'_2) - h(z_1) - h( z_2)$
as in \eqref{eq: omega collision}.
\end{Thm}

%
%

\medskip

Since the approximation holds true for very long times compared to the mean free time (diverging to infinity as $\log\log \log \mu_\eps$), it makes sense to look at fluid limits, i.e. at regimes when the collision process is much faster than the transport (density is still low but makes the collisions a bit more likely) $\mu_\eps \eps^{d-1} = \alpha^{-1}$ with $\alpha \gg \eps$, $\alpha \to 0$. Starting from the scaled linearised Boltzmann equation  
$$ 
\d_t h +v\cdot \nabla_x h = \frac1\alpha \Big( C(h, M) +C(M,h) \Big) \,,
$$
it is well known \cite{Bardos_golse_levermore}  that in the limit  $\alpha \to 0$, the gas will be close to a local thermodynamic equilibrium, with density, bulk velocity and  temperature satisfying the acoustic equations. Zooming out on longer times $O(1/\alpha)$, these acoustic waves become fast oscillating and thus converge weakly to 0, but the incompressible component has a diffusive behavior, satisfying the Stokes-Fourier equations.
This by now classical asymptotic analysis can be actually combined with Theorem \ref{thm3} to derive directly the Stokes-Fourier equations from the dynamics of hard spheres as in \cite{BGSR2}. 
In a work in progress, we also take into account the noise, and get the corresponding fluctuating hydrodynamics (satisfying the fluctuation-dissipation principle).

\subsection{Some elements of the proof of Theorem \ref{thm3}}

As in the previous sections, we will not enter into the technicalities of the proof, which is actually quite involved. We will just focus here on some key arguments, providing a better understanding of  large time asymptotics.
We  work directly on moments of the fluctuation field, defined for any collection of times $\theta_1 < \dots \ < \theta_p$ by
\begin{equation}
\label{eq: moment P}
\bbE_\eps \Big[ \la \zeta^\eps_{\theta_1}, h_1\ra \;  \dots 
\;  \la \zeta^\eps_{\theta_p} , h_p \ra \Big] \,,
\end{equation}
 and we are going to prove their convergence to the moments of the field in the stochastic equation
 $d\zeta_t =   \cL_{eq} \zeta_t dt   +  d\eta_t$.  Combined with the tightness results from 
 \cite{BGSS2}, this fully characterises the convergence of the microscopic fluctuation field.

\medskip

\begin{figure}[h] 
\centering
\includegraphics[width=2.5in]{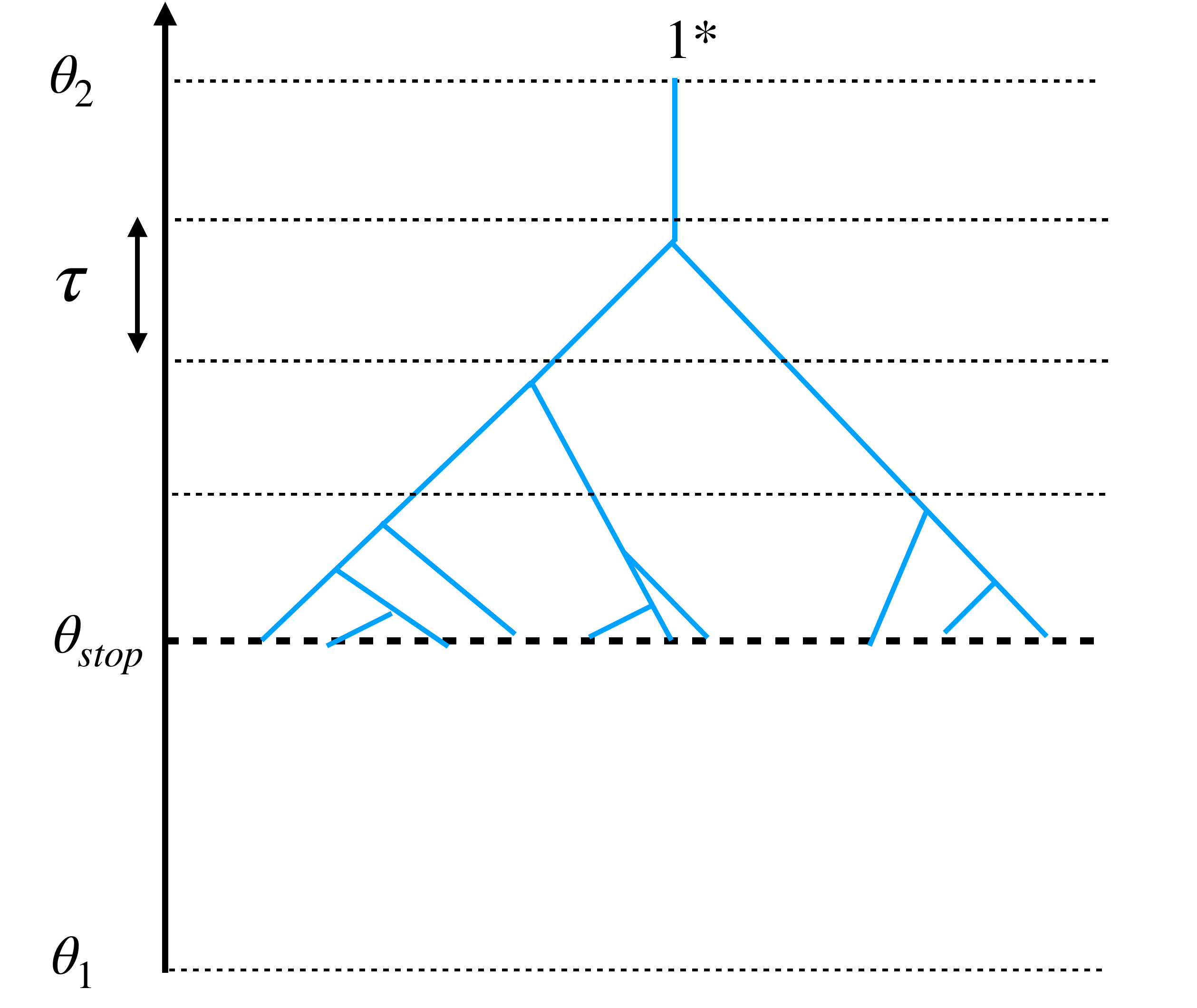} 
\hskip1cm
\includegraphics[width=.62in]{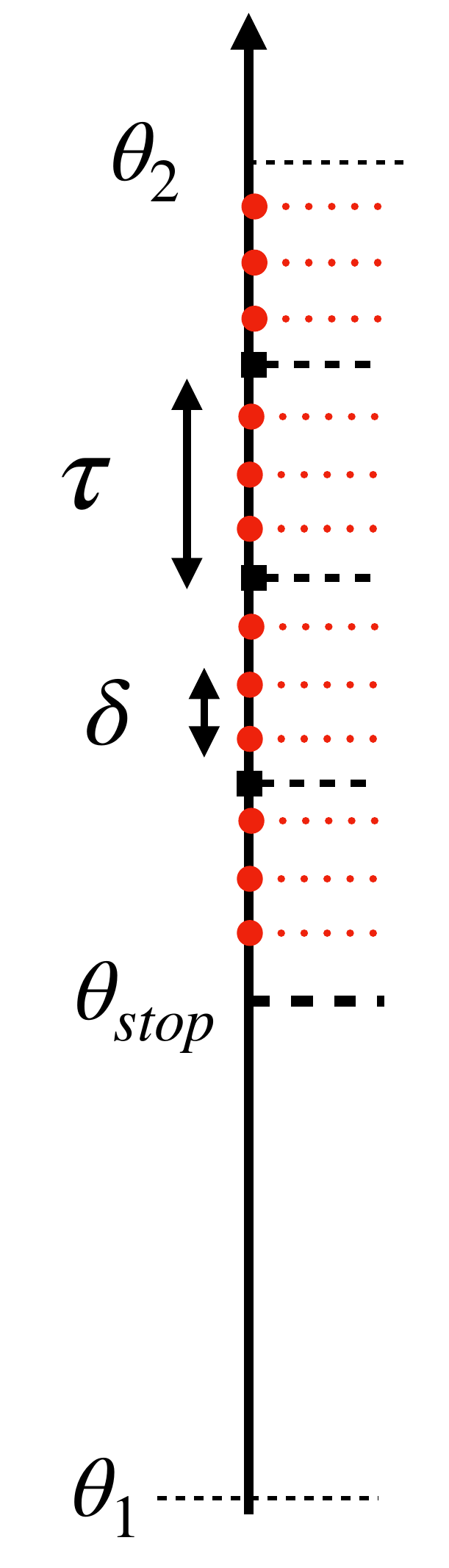} 
\caption{\color{black} Pseudo-trajectories are build iteratively on short time intervals of length $\tau$ starting from $\theta_2$. The procedure stops before reaching time $\theta_1$ if superexponential branchings occur in a time interval of length $\tau$. 
The corresponding pseudo-trajectories stop at time $\theta_{\text{stop}}$ and are then discarded.
A double sampling at scales $\delta \ll \tau \ll 1$, depicted on the right figure,  is implemented to control the recollisions. 
}
\label{superexponential-fig}
\end{figure}

Let us start with $p=2$ and compute the covariance $\bbE_\eps \Big[ \zeta^\eps_{\theta_1}(h_1) \, \zeta^\eps_{\theta_2} (h_2) \Big]$.
The idea is to pull back the observable $h_2$  from time $\theta_2$ to $\theta_1$ in order to reduce the estimates  at  a single time $\theta_1$. 
A similar strategy was presented in 
Sections \ref{lanford-proof} and \ref{sec: Defects in the chaos assumption} to transport  the correlation up to time 0 for which the distribution was known.
In particular, we have seen that the correlation functions at a time $\theta_2$
can be represented by backward pseudo-trajectories involving collision trees with a number $m$ of additional particles encoding the dynamical history during the time interval $[\theta_1,\theta_2]$.
The time restriction $T_L$ for the convergence to Boltzmann equation in Theorem \ref{lanford-thm} was due to the lack of control on the growth of the tree sizes $m$  at large times. 
Indeed dynamical correlations may develop and form giant components of correlated particles for very pathological trajectories. In order to reach larger time scales, one has to show that the contribution of these bad  trajectories with large $m$ remains negligible. 
For this we perform a \emph{time sampling}.
The idea is to build the pseudo-trajectories iteratively from $\theta_2$ to $\theta_1$  on  time steps of length $\tau \ll 1$ and to neglect the collision trees with a fast (superexponential) growth  during a time  $\tau$ (see Figure \ref{superexponential-fig}).
The large collision trees are therefore discarded before they reach the time $\theta_1$, i.e. before their sizes become uncontrollable. This can be achieved by using the time invariance property of the equilibrium measure which provides a priori controls on the statistics.
This kind of sampling was introduced for the first time in the context of the Boltzmann-Grad limit  in \cite{BGSR1, BGSR2}, but it is also an important ingredient in the weak coupling limit for quantum systems leading to quantum diffusion \cite{erdos-yau, erdos}.

Another key ingredient, to derive the convergence to the Boltzmann equation, is the
 procedure to neglect the  ``bad'' trajectories  involving recollisions  (see Section \ref{lanford-proof}). 
Controlling the growth of the collision trees is also essential to discard recollisions.
The idea is  to introduce a double sampling in time (with time scales 
$\delta \ll \tau \ll 1$,  see Figure \ref{superexponential-fig})
which takes care simultaneously of the recollisions and of the collision tree growth.
The backward  iteration is stopped  and the corresponding pseudo-trajectories are discarded
as soon as one of the following conditions is violated :
\begin{itemize}
\item there is at least one recollision on the last  very small interval of size $\delta = O(\eps^{1-\frac1{2d}})$;
\item on the last small interval of size $\tau = O( \log \, \log \mu_\eps)^{-1/2}$ the number of particles has been multiplied at least by 2.
\end{itemize}
Note that both conditions are entangled. 
On the one hand, the bigger the size of the system, the easier for recollisions to occur. On the other hand, it is rather  difficult to control the growth of the system if there are recollisions.

\begin{figure}[h] 
\centering
\includegraphics[width=4.5in]{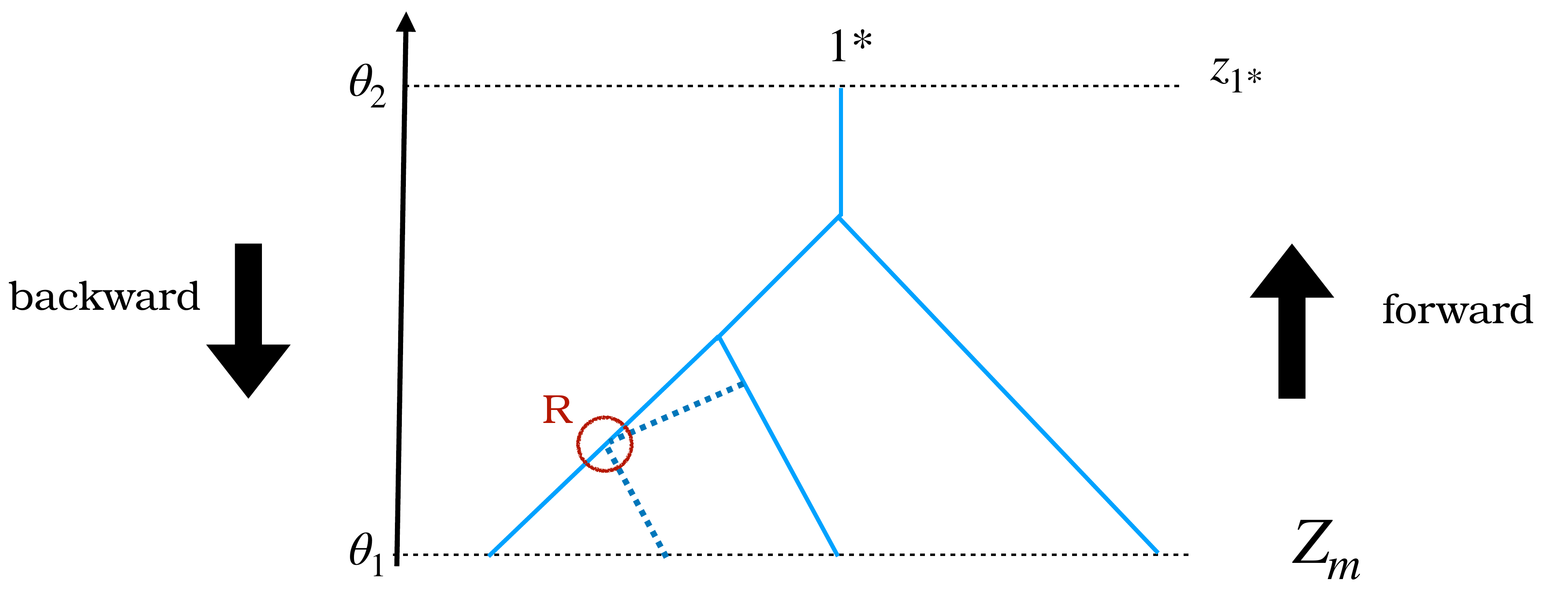} 
\caption{\color{black} Starting from $z_{1^*}$ at time $\theta_2$, the blue pseudo-trajectory is built backward and  leads to a configuration $Z_m$ at time $\theta_1$ (with $m = 3$ on the picture). 
The dual procedure goes forward, starting from $Z_m$ in order to reconstruct $z_{1^*}$ as a function of $Z_m$ at time $\theta_2$. Following the forward flow, a tree is built by removing one of the particles at each encounter between two particles. Notice that one has to choose which particle will be removed and if a scattering occurs. Thus there are potentially several ways to build forward trajectories, but their 
combinatorics is well under control.
This is no longer the case when recollisions can occur.
Indeed this adds the possibility that when two particles encounter in the forward flow, 
none of them disappears  so    when the number of recollisions 
is not bounded the combinatorics diverges.
 }
\label{forward-fig}
\end{figure}

Assuming  that the pseudo-trajectories can be controlled by the previous time sampling, 
let us now explain the \emph{weak convergence method} for  computing the covariance.
The two-time correlation  $\bbE_\eps \big[ \la \zeta^\eps_{\theta_1}, h_1 \ra \;  \la \zeta^\eps_{\theta_2} , h_2 \ra \big] $ can be rephrased as the expectation of two fluctuation fields at the same time $\theta_1$ 
\begin{equation}
\label{duality}
\bbE_\eps \Big[ \la \zeta^\eps_{\theta_1}, h_1 \ra \;  \la \zeta^\eps_{\theta_2} , h_2 \ra \Big] 
\; "=" \;
\sum_{m} \bbE_\eps \left[ \la \zeta^\eps_{\theta_1}, h_1 \ra \; 
\la   \zeta^\e_{m,\theta_1  } , \phi_{\theta_2-\theta_1}  \ra  \right]  \,,
\end{equation}
where the new test function $\phi_{\theta_2-\theta_1}(Z_m)$ is  obtained from $h_2$ by considering all possible forward flows starting from $Z_m$ at time $\theta_1$ and having only one particle left at time $\theta_2$
(see Figure \ref{forward-fig}). 
In this sense, \eqref{duality} is dual to the backward representation of the correlation functions
\eqref{F1-series}. The price to pay, to reduce the expectation at a single time,
is that the new test function $\phi_{\theta_2-\theta_1}$ depends on $m$ particles (a parameter related to the size of the collision trees  in the time interval $[\theta_1, \theta_2]$) so that the fluctuation field $\zeta^\eps_{m, \theta_1}$ has the form
$$
\la \zeta^\eps_{m, \theta_1} , \phi_{\theta_2-\theta_1} \ra 
= \sqrt{\mu_\eps}
\left( \frac1{\mu_\eps^m} \sum_{(i_1,\dots, i_m)} \phi_{\theta_2-\theta_1}  \big( z^\e_{i_1}(\theta_1), \dots, z^\e_{i_m}(\theta_1)\big)
- \bbE_\eps( \phi_{\theta_2-\theta_1} ) \right)\, ,
$$
which is related to  the generalized empirical measure defined in  \eqref{eq: Fk}, with the abbreviation 
$$\bbE_\eps( \phi_{\theta_2-\theta_1} ) = \bbE_\eps \left(  \la \pi_{m,\theta_1}^\eps , \phi_{\theta_2-\theta_1} \ra  \right)\;. $$
In the following, we will abusively forget the subscript $m$.


The difficulty to make sense of the pullback in (\ref{duality})  is that the forward flow is not a priori well-defined. Indeed different backward pseudo-trajectories may end up at time $\theta_1$ with the same particle configuration $Z_m$. Thus starting from $Z_m$, there are many possibilities to build the  forward flow from $\theta_1$ to $\theta_2$:
when two particles touch each other, we need to prescribe whether one of them will be deleted (corresponding to a creation in the backward flow) or not (corresponding to a recollision), and in the case of deletion whether there is scattering of the remaining particle
(see Figure \ref{forward-fig}). 
The combinatorics of these choices is diverging very fast if the number of recollisions is not under control.
The very short time sampling~$\delta$ is introduced so that the number of recollisions during a time $\delta$ is controlled with high probability under the equilibrium measure.

Then the pullback relation \eqref{duality} is obtained by successive iterations of the sampling time $\delta$. After the first elementary time step in the time interval $[\theta_2 - \delta, \theta_2]$, the pathological events are discarded and then the elementary pullback can be iterated. This means that, at each time~$\theta_2 - r\delta$, remainder terms due to recollisions are neglected, and that at each time~$\theta_2 - k\tau$, remainder terms due to superexponential growth can also be  discarded. 
Let~$\theta_{\text{stop}}\in [\theta_1,\theta_2]$, be the first time at which a pseudo-trajectory becomes pathological (see Figure \ref{superexponential-fig}).
The corresponding terms obtained by forward transport from the time $\theta_{\text{stop}}$ are generically denoted by $\phi_{\theta_{\text{stop}}}^{\text{bad}}$ and  are proved to be small by 
using  the time invariance of the equilibirum measure. Indeed  the \emph{time decoupling} follows from a Cauchy-Schwarz estimate
\begin{equation}
\label{eq: decoupling covariance}
\Big|  \bbE_\eps \Big[   \la \zeta^\eps_{\theta_1}, h_1 \ra \; 
\la \zeta^\eps_{\theta_{\text{stop}}} , \phi_{\theta_{\text{stop}}}^{\text{bad} } \ra \Big]  \Big| 
\leq  \bbE_\eps \Big[    \la \zeta^\eps_{\theta_1}, h_1 \ra^2 \Big] ^{1/2} 
\; \bbE_\eps \Big[    \la \zeta^\eps_{\theta_{\text{stop}}}, \phi_{\theta_{\text{stop}}}^{\text{bad}}  \ra^2 \Big] ^{1/2}\,,
\end{equation}
and from the strong geometric constraints on the corresponding pathological pseudo-trajectories which  
can be estimated under the equilibrium measure on can deduce that:
$$
\bbE_\eps \Big[   \la \zeta^\eps_{\theta_{\text{stop}}}, \phi_{\theta_{\text{stop}}}^{\text{bad}}  \ra^2  \Big]   \to 0 
\quad \text{as} \quad  
\mu_\eps \to \infty\,.
$$

%
%

\medskip

The last important step to prove that the   limiting process is Gaussian
boils down to showing that,  asymptotically when $\mu_\eps \to \infty$, the moments, defined in \eqref{eq: moment P}, are determined by the covariances according to \emph{Wick's rule}
\begin{equation} 
\label{eq: appariement} 
\lim_{\mu_\eps \to \infty} \left| \bbE_\eps \Big[ \la \zeta^\eps_{\theta_1}, h_1\ra \;  \dots 
\;  \la \zeta^\eps_{\theta_p} , h_p \ra \Big] 
-  \sum_{\eta \in {\mathfrak S}_p^{\text{pairs}}}
\ \prod_{\{i,j\}  \in \eta} \bbE_\eps \left[ \la \zeta^\eps_{\theta_i} ,  h_i \ra \; 
\la\zeta^\eps_{\theta_j} , h_j \ra  \right] \right| = 0 \, ,
\end{equation}
where ${\mathfrak S}_p^{\text{pairs}}$ is the set of partitions of $\{1, \dots , p \}$ made only of pairs. 
Notice that if $p$ is odd then ${\mathfrak S}_p^{\text{pairs}}$ is empty and the product of the moments is asymptotically 0.

To understand this pairing mechanism, let us start with a simpler example  for which explicit computations can be achieved.
Consider the moments of the fluctuation field at time~0, under the equilibrium measure with independently distributed particles.
This reduces to the case $\eps = 0$ and $\theta_1 = \dots = \theta_p =0$. Assuming furthermore that the test functions are of mean $\bbE_0 (h_i) = 0$ (we abusively write here $\bbE_0$ for this \textit{iid} case, not to be confused with $\bbE_\eps$ for $\eps = 0$), we get 
\begin{equation} 
\label{eq: appariement initial} 
 \bbE_0 \Big[ \prod_{\ell =1}^ p \la \zeta^\eps_0, h_\ell \ra \Big] =
\frac{1}{\mu_\eps^{p/2}}
 \bbE_0 \Big[ \prod_{\ell =1}^p \Big( \sum_{i_\ell} h_\ell (z_{i_\ell}) \Big)  \Big]
= \frac{1}{\mu_\eps^{p/2}}
 \bbE_0 \Big[ \sum_{i_1, \dots , i_p} \prod_{\ell =1}^p  h_\ell (z_{i_\ell})   \Big],
\end{equation}
where the sum is over all the possible choices (with repetition) among 
 $N$ particles (with~$N \simeq \mu_\eps$ under the grand canonical measure). 
As the mean of the test functions is assumed to be 0, each particle has to be chosen 
at least twice, otherwise by the independence of the variables the expectation is equal to 0.
Thus in the sum $i_1, \dots , i_p$ the number $k$ of different particles is such that 
$k \leq p/2$.
Choosing $k$ different particles gives a combinatorial factor $\mu_\eps^k$ so that 
only the pairings with $k = p/2$ and $p$ even contribute to the limiting moment.
In this way, one recovers the Wick decomposition 
\eqref{eq: appariement} in terms of pairings.
Note that for $\eps>0$, a similar result holds (at time zero) in the Boltzmann Grad limit, but a cluster expansion of the equilibrium measure is necessary to control the (weak) correlations of the Gibbs measure.

{\color{black} For time dependent fluctuation fields, the pairing cannot be achieved in one step as in the previous example. One has instead to proceed iteratively. Let us revisit the computation above  to explain the idea first in this simple setting.
We start by focusing on the product of two fields and decompose it as follows  
\begin{equation} 
\label{eq: appariement 2 fields} 
 \la \zeta^\eps_0, h_p \ra
 \la \zeta^\eps_0, h_{p-1} \ra
 = \underbrace{\frac{1}{\mu_\eps} \sum_i h_p (z_i) h_{p-1} (z_i)}_{ = \; \Psi } 
 + \underbrace{\frac{1}{\mu_\eps} \sum_{i \not = j} h_p (z_i) h_{p-1} (z_j)}_{
= \; \la \zeta^\eps_0, h_p \ra  \otimes
 \la \zeta^\eps_0, h_{p-1} \ra}.
\end{equation}
The pairing between $\la \zeta^\eps_0, h_p \ra$ and $\la \zeta^\eps_0, h_{p-1} \ra$ is coded by
the function $\Psi$  which is called a \emph{contracted product} as the variables are repeated. 
As the variables are independent, the covariance between $h_p$ and $h_{p-1}$ is given by 
\begin{equation} 
\label{eq: covariance 2 fields} 
\bbE_0 \Big[ \la \zeta^\eps_0, h_p \ra \, \la \zeta^\eps_0, h_{p-1} \ra \Big]  =  \bbE_0 [ \Psi ].
\end{equation}
From the central limit theorem, $\Psi$ can be interpreted as a small fluctuation around the covariance
\begin{equation} 
\label{eq: Psi 2 fields} 
\Psi =  \bbE_0 [ \Psi ] + \frac{1}{\sqrt{\mu_\eps}} \widehat \Psi
\quad \text{with} \quad 
\widehat \Psi = \frac{1}{\sqrt{\mu_\eps}}  \left( \sum_i  h_p (z_i) h_{p-1} (z_i) - \mu_\eps \bbE_0 [ h_p \, h_{p-1} ]  \right),
\end{equation}
where $\widehat \Psi$ behaves as a random variable with finite covariance (uniformly in $\eps$).
The second term in \eqref{eq: appariement 2 fields} will be called a $\otimes$-product and denoted by $ \la \zeta^\eps_0, h_p \ra  \otimes  \la \zeta^\eps_0, h_{p-1} \ra$. It behaves qualitatively as a fluctuation field as the variables are not repeated.

Returning to~(\ref{eq: appariement initial}), to extract the pairing between $\la \zeta^\eps_0, h_p \ra$ and $\la \zeta^\eps_0, h_{p-1} \ra$, we write 
\begin{equation} 
\label{eq: appariement 1} 
 \bbE_0 \Big[ \prod_{\ell =1}^ p \la \zeta^\eps_0, h_\ell \ra \Big] =
 \underbrace{\bbE_0 \Big[ \Big( \prod_{\ell =1}^{p-2} \la \zeta^\eps_0, h_\ell \ra \Big) \Psi  \Big]}_{\text{pairing of 
$h_p, h_{p-1}$}}
+ \underbrace{ \bbE_0 \Big[ \Big( \prod_{\ell =1}^{p-2} \la \zeta^\eps_0, h_\ell \ra \Big) \;
\Big(  \la \zeta^\eps_0, h_p \ra  \otimes
 \la \zeta^\eps_0, h_{p-1} \ra  \Big) \,  \Big]}_{\text{product of $p-1$ fields}}.
\end{equation}
The second term can be seen as a product of  $p-1$ fields which will be treated recursively at the next step.
The pairing between $\la \zeta^\eps_0, h_p \ra$ and $\la \zeta^\eps_0, h_{p-1} \ra$ can be extracted from the first term as follows.
Using the decomposition \eqref{eq: Psi 2 fields} , we get 
\begin{align*} 
\bbE_0 \Big[ \Big( \prod_{\ell =1}^{p-2} \la \zeta^\eps_0, h_\ell \ra \Big) \Psi  \Big]
& = \bbE_0 \Big[   \prod_{\ell =1}^{p-2} \la \zeta^\eps_0, h_\ell \ra \Big]
\bbE_0 [ \Psi ] + \frac{1}{\sqrt{\mu_\eps}} \bbE_0 \Big[ \Big( \prod_{\ell =1}^{p-2} \la \zeta^\eps_0, h_\ell \ra \Big) \; \widehat \Psi   \Big]\\
& =  \bbE_0 \Big[   \prod_{\ell =1}^{p-2} \la \zeta^\eps_0, h_\ell \ra    \Big]
\; \bbE_0 \Big[ \la \zeta^\eps_0, h_p \ra \, \la \zeta^\eps_0, h_{p-1} \ra \Big] 
+ O \Big(  \frac{1}{\sqrt{\mu_\eps}} \Big),
\end{align*}
where the smallness of the last term follows from H\"older's inequality
\begin{align} 
\label{eq: Holder} 
\left| \bbE_0 \Big[ \Big( \prod_{\ell =1}^{p-2} \la \zeta^\eps_0, h_\ell \ra \Big) \; \widehat \Psi   \Big] \right|
\leq 
\bbE_0 \Big[\widehat \Psi^2   \Big]^{\frac{1}{2}} \; 
\prod_{\ell =1}^{p-2} \bbE_0 \Big[ \la \zeta^\eps_0, h_\ell \ra^{2(p-2)}  \Big]^{\frac{1}{2(p-2)}} 
 ,
\end{align}
provided bounds on the moments of single fields can be obtained.
For independent variables, this procedure is far from optimal, however it will be extremely useful to decouple fields at different times. 
In this way,  the pairing between $\la \zeta^\eps_0, h_p \ra$ and $\la \zeta^\eps_0, h_{p-1} \ra$  can be extracted  without investigating the correlations between these two fields and the $p-2$ other fields. 
Note that a time decoupling inequality similar to \eqref{eq: Holder}  was used in the computation of the covariance 
\eqref{eq: decoupling covariance} to neglect bad pseudo-trajectories.
Finally, it  remains to iterate this procedure with $\bbE_0 \Big[   \prod_{\ell =1}^{p-2} \la \zeta^\eps_0, h_\ell \ra    \Big]$ and the second term in \eqref{eq: appariement 1}  which involves a product of at most $p-1$ fluctuation fields.
}

We turn now to the time dependent case \eqref{eq: appariement} and 
proceed backward in time to achieve the pairing step by step.
First the fluctuation at time $\theta_p$  is pulled back at time~$\theta_{p-1}$
as a sum of (more complicated) fluctuations  by the same duality method as for the covariance~\eqref{duality}.
Using analogous notation as in~\eqref{duality}, the test function $h_p$ is transformed into a function $\phi_{\theta_p - \theta_{p-1}} ^{(p)}$ with $m$ variables.
Forgetting for a moment the product $\prod_{\ell =1}^{p-2} \la \zeta^\eps_{\theta_\ell}, h_\ell \ra$,
we focus on the product of the fields at time $\theta_{p-1}$ 
\begin{equation}
\label{eq: produit en p-1}
\la  \zeta^\eps_{\theta_{p-1}} , h_{p-1} \ra \; 
\la  \zeta^\e_{m,\theta_{p-1} } , \phi_{\theta_p - \theta_{p-1}} ^{(p)} \ra
\end{equation}
and decompose it as in \eqref{eq: appariement 2 fields}  according to the repeated indices in the spirit of the example above.
This leads to two types of contributions~:
\begin{itemize}
\item 
{\color{black}
a ``contracted product'' (by analogy with the function $\Psi$) which records all the repeated indices in the product \eqref{eq: produit en p-1} at time $\theta_{p-1}$.
By H\"older's inequality as in \eqref{eq: Holder}, this term can be decoupled from the rest of the weight 
formed by the moments $\prod_{\ell =1}^{p-2} \la \zeta^\eps_{\theta_\ell} , h_\ell \ra$.
This strategy is particularly relevant for time dependent fields as 
it reduces the estimates to computing moments of fields at a single time.
In an equilibrium regime, the  moments of the field at a single time  can be easily analysed as the distribution is time  invariant. 
In this way the moments at $\theta_p$ and $\theta_{p-1}$ are paired and their covariance
$\bbE_\eps [ \la  \zeta^\eps_{\theta_{p-1}} , h_{p-1} \ra \; 
\la   \zeta^\eps_{\theta_p} , h_p \ra ]$  is recovered.
It remains then to study the remaining moments $\bbE_\eps \Big[ \prod_{\ell =1}^{p-2} \la \zeta^\eps_{\theta_\ell} , h_\ell \ra \Big]$.}
\item a ``$\otimes$-product'', which by definition takes into account the non-repeated indices, 
and which can be interpreted as a product of two independent fluctuations at time $\theta_{p-1}$.
In a very loose way, we have to evaluate now the following structure 
$$
\bbE_\eps \Big[  \Big( \prod_{\ell =1}^{p-2} \la \zeta^\eps_0, h_\ell \ra \Big)
\; \la  \zeta^\eps_{\theta_{p-1}} , h_{p-1} \ra \otimes 
\la  \zeta^\e_{m,\theta_{p-1} } , \phi_{\theta_p - \theta_{p-1}} ^{(p)} \ra
\Big],
$$
with a more complicated fluctuation field at time $\theta_{p-1}$.
\end{itemize}
 The key point here is that using the cumulant techniques introduced in Section \ref{cumulant-proof}, one can then prove that the tensorized structure $\otimes$ is essentially preserved by the pullback of test functions : 
the configurations for which the $\otimes$-product breaks can be neglected. 
Thus with high probability the fields $\la  \zeta^\eps_{\theta_{p-1}} , h_{p-1} \ra \otimes 
\la  \zeta^\e_{m,\theta_{p-1} } , \phi_{\theta_p - \theta_{p-1}} ^{(p)} \ra$ 
can be pulled back up to time~$\theta_{p-2}$  as if they were independent.
Then we apply the pairing procedure at time~$\theta_{p-2}$.
This leads to new pairings between $\la  \zeta^\eps_{\theta_{p-2}} , h_{p-2} \ra$ and the pulled-back fields.
In particular, the covariances $ \bbE_\eps [ \la  \zeta^\eps_{\theta_{p-2}} , h_{p-2} \ra \; 
\la   \zeta^\eps_{\theta_p} , h_p \ra ]$ and $ \bbE_\eps [ \la  \zeta^\eps_{\theta_{p-2}} , h_{p-2} \ra \; 
\la   \zeta^\eps_{\theta_{p-1}} , h_{p-1} \ra ]$ can be identified.
The non repeated variables at time $\theta_{p-2}$ build new $\otimes$-products involving the fluctuation fields (or their pullbacks) from times $\theta_{p-2}$, $\theta_{p-1}$ and $\theta_p$.

Iterating this procedure up to time $\theta_1$ all the pairings can be recovered and 
the Wick's decomposition \eqref{eq: appariement} is obtained in the limit $\mu_\eps \to \infty$.
This shows that the limiting process is Gaussian, thus achieving the proof of Theorem~\ref{thm3}.

\section{OPEN PROBLEMS AND PERSPECTIVES}
\label{sec: open problems}

The research program that we conducted during this last decade and which is presented in this survey has led to two important breakthroughs compared to the state of the art after Lanford's theorem:
\begin{itemize}
\item 
an extended statistical picture of the dynamics of hard-sphere gases for  short times, including fluctuations and large deviations;
\item a complete answer to Hilbert's sixth problem connecting the three levels of modeling (atomistic, kinetic and fluid) for  linear equations of dilute hard-sphere gases close to equilibrium.
\end{itemize}
Nevertheless, the problem of the axiomatization of gas dynamics remains largely open, even in dilute regimes. We propose in this final section to review some important directions to be explored in the future.   We choose to discuss here only kinetic limits, involving a separation of scales, for which an enterprise in the spirit of the one discussed above is conceivable (albeit possibly hard).

\subsection{Long time behavior for dilute gases}

The only case in which we have a complete picture of the transition from the atomistic description to fluid models is the \textbf{equilibrium} case. 
Nevertheless the diffusive scaling  considered in these linear regimes is 
sub-logarithmic (see e.g.\,\cite{BGSR1, BGSS3}). It would be interesting to reach more relevant physical scales, for which we expect the limiting picture to remain unchanged.

The law of large numbers in the equilibrium case is trivial, and the fluctuations are governed by linear models. In order to extend this analysis to gases which are initially \textbf{out of equilibrium}, a major obstruction is to define a good notion of stability for the nonlinear Boltzmann equation, which plays the role of pivot between the microscopic and macroscopic scales. In other words, this requires designing a good notion of convergence. The weak convergence method developed in the equilibrium case uses a topology which is a priori too weak to make sense of the nonlinear collision operator. 
 Based on our analysis, we believe that stronger convergence methods require a rather precise understanding of the mechanisms responsible for the entropy cascade through the cumulants,  retaining enough information in the limiting system. 
Note that this information is encoded in the supports of the cumulants, which have a finer and finer structure as the order of the cumulant increases.
This structure  might well be 
a key ingredient, as entropy and entropy dissipation play a crucial role in the stability of the Boltzmann equation.

Beyond the law of large numbers, it would be also natural to extend the analysis of fluctuations and large deviations for long kinetic times, and even diffusive times. This would allow to derive the fluctuating hydrodynamics (typically the fluctuating Navier-Stokes-Fourier equations). A fine understanding of the Hamilton-Jacobi equations and of the associated gradient structure would be certainly a major step in this direction.

\subsection{The role of microscopic interactions}

Our study is focused on the  case of hard sphere gases, for which the interaction is pointwise in time and the scattering law  is very simple. The papers \cite{King,GSRT,Pulvirenti_Saffirio_Simonella} have shown that, despite technical complications, the same average behaviour,  in the low density limit, is obtained for compactly supported potentials satisfying some suitable lower bound (thermodynamic stability). Only the collision cross section (i.e. the transition rate of the jump process in the velocity space) and consequently the hydrodynamic transport coefficients are modified. 
 One expects, and can prove for short times \cite{King}, that multiple collisions (three or more particles simultaneously interacting at a given time) are a correlation of higher order with respect to the dynamical correlations determining the fluctuation theory.  It is then very likely that the description of fluctuations and large deviations for short times can be also extended to this \textbf{short-range} case. 
Notice that the absence of monotonicity of the potential would require a more delicate treatment, as some trajectories can be trapped for a very long time \cite{Pulvirenti_Saffirio_Simonella}.

A problem of a much higher level of difficulty is to deal with  \textbf{long-range} interactions. We know that, as soon as the potential is not compactly supported, the collision cross-section (which can be computed by solving the two-body problem) has a non integrable divergence at grazing angles. It is therefore impossible to define solutions of the Boltzmann equation without taking into account the cancellations between the gain and loss terms in the collision operator, which would imply to find new ideas (in our methods dealing with microscopic systems, such cancellations are never used). Close to equilibrium, using a sampling to discard superexponential growth (as in Section \ref{time-sec} above), N.\,Ayi~\cite{Ayi} has proved a convergence result for very fast decaying potentials, but the method does not seem robust enough to deal with weaker decays nor systems out of equilibrium.

A natural idea, often used by physicists, would be to decouple the short range part (acting as ``collisions''), and the long range part of the interaction potential (to be dealt with by mean field methods). 
  However, from the mathematical point of view, this leads to a major issue: no analysis method is available so far, as the techniques used for the low density limit and for the mean field limit are completely different and apparently incompatible. This problem is investigated in \cite{DSS}, where a linear Boltzmann-Vlasov equation is derived rigorously for a simple (Lorentz gas) model system (see also \cite{Desvillettes_Pulvirenti}).
 
 A related issue is how to precisely identify and separate the long range and the collisional part for a given potential law, capturing the good scaling for both parts. There are some delicate aspects here involving the details of the potential and the dimension of the problem \cite{Nota_Velazquez_Winter,NVW2}. Formal considerations as in \cite{Bobylev_Dukes_Illner}  indicate that, in case of power law potentials $1/x^s$, the low density scaling should lead to a Boltmann equation for $s > d-1$, to a Boltzmann-Vlasov equation for $s=d-1$, and to a Vlasov equation (with Boltzmann's operator still describing the collisions as a long time correction) for~$s \in (d-2,d-1)$. For the Coulomb potential (and for smaller values of $s$), the Boltzmann operator has to be replaced by a diffusive variant of it (Landau, or Lenard-Balescu operator; see also Section \ref{sec:realmKL}). We refer to \cite{Nota_Velazquez_Winter} for details.

We remark that the combination of mean-field and collisions has an interest in connection with the problem of binary mixtures exhibiting phase segregation \cite{Bastea_Esposito_Lebowitz_Marra} (see also \cite{AMP21} on a derivation result  for mixtures).

%
%
%
%
%
%
%
%

We stress that Lanford's theorem can be seen as a   \textbf{propagation of chaos} result:  the strategy used to prove the kinetic limit boils down indeed to transferring the initial independence property to an independence at time $t$, thus recovering  the molecular chaos assumption~(\ref{eq: molecular chaos}). 
Another direction would be trying to study initial data with strong correlations, preventing the validity of our results at time zero or at short enough times. For instance, one could think of particles initially arranged on a lattice in position space; or construct pathological initial measures with a defect of convergence on precollisional configurations as done in \cite{BGSS}. In such cases, one could hope that, after some time, a different mixing mechanism takes place, producing a form of local equilibrium in a dynamical way.

\subsection{Non equilibrium   stationary states}

For short times, Lanford's theorem allows to consider particle systems which are initially put out of equilibrium, provided that their distribution is controlled in some sense by an equilibrium state. This assumption is a key argument to get uniform bounds (even for short times when the relaxation phenomenon cannot be observed). In this situation, one can use a comparison principle, because nothing forces the system to stay out of equilibrium, and the invariant measure is well known.

A natural extension is to deal with a gas evolving in a domain with \textbf{boundary conditions}, rather than the whole space or the periodic setting as considered previously.
In the case of boundary conditions ensuring conservation of energy, we still have a control by the invariant measure, and the main extra difficulty caused by the presence of boundaries lies in the  geometric analysis of recollisions. This has been discussed so far in the case of simple geometries \cite{Dolmaire, LeBihan} (see also \cite{EspositoGuoMarra} for the case of external forces).

A much more delicate situation is when the system of interacting particles is maintained out of equilibrium by a forcing or a boundary condition (reservoir, thermostat,\dots). One would like to derive, in this non-equilibrium framework, the Boltzmann equation and more generally the properties of the steady states. 
As exposed in \cite{Bonetto_Lebowitz_ReyBellet}, this question is a "challenge to theorists" and few quantitative results are known either for gas dynamics or for other mechanical systems such as chains of anharmonic oscillators. 
Even though, under reasonable assumptions on the non-equilibrium forces, 
the existence of a \emph{stationary measure} of the microscopic dynamics is expected, one does not know how to construct such a measure   or any exact solution which would  play the role of supersolution for the actual distribution of particles. In particular, a good starting point for the analysis of the low density limit seems to be missing at present. Finally, it is worth mentioning that the theory of stationary solutions for the Boltzmann equation with thermal reservoirs is still far from mature, see \cite{EspositoMarra_Stationary} for a recent review.

Beyond the derivation of the Boltzmann equation for boundary driven systems, it would be interesting to investigate the large deviations as they can provide some knowledge on the 
invariant measure \cite{Bertini_DeSole_Gabrielli, Derrida}. 
Also it is conjectured  \cite{BGSR_cras, Bonetto_Lebowitz_ReyBellet} that the Fourier 
law should be valid for  a dilute gas maintained out of equilibrium by reservoirs. To prove its validity  would require an analysis beyond the kinetic time scale in order to derive fluid equations out of equilibrium.


\subsection{A realm of kinetic limits} \label{sec:realmKL}

Besides the low density (Boltzmann-Grad) scaling discussed so far, there is a variety of interacting particle models admitting a kinetic limit and sharing many similarities with the classical Boltzmann gas \cite{spohn2012large}.
We shall only mention here the two main obvious modifications of our assumptions (which are reviewed in detail in \cite{Pu06}): 
(i) start from a microscopic description based on \textbf{quantum mechanics} instead of classical mechanics, namely replace the Newton equations by the $N$-body Schr\"odinger equation, including additional symmetry/antisymmetry constraints which take into account the specificity of bosons/fermions;
(ii) perform a high-density, \textbf{weak-coupling} scaling with potential $\e^{\alpha}\phi(x/\e)$, where $\alpha \in [0,1]$
and the particle density is correspondingly tuned as $-d+1-2\alpha$. For $\alpha \in (0,1)$, the latter scaling should lead to the diffusive Landau equation in the case of classical systems, and is suited to a description of collisions in plasmas. The diffusion emerges from a central limit type effect on an accumulation of many weak collisions. The limiting point $\alpha = 1$ is expected to capture the famous Lenard-Balescu correction. Conversely in the case of quantum systems, each value of $\alpha$ should lead to a quantum version of the Boltzmann equation. The amount and quality of quantum features surviving in the limit depends on the particular value of $\alpha$. For $\alpha = 0$, the collision operator contains the full quantum cross-section. On the other hand for $\alpha = 1/2$ (when only the first term of the Born series survives), one expects to get additional cubic terms in the collision operator, expressing the inclination of particles to aggregate (Bose-Einstein condensation) or to repel each other (Pauli's exclusion principle).

For such a variety of situations, no rigorous full derivation result is available at present, not even for short kinetic times; see however \cite{Pu06, Bobylev_Pulvirenti_Saffirio,Velazquez_Winter,Winter} for consistency results and attempts in this direction (full results are instead available for Lorentz type (linear) models, see \cite{KP,DGL} for the classical case and \cite{erdos} for a review in the quantum case).
When trying to reproduce Lanford's strategy, one stumbles indeed upon many difficulties. The construction of the equilibrium measure is delicate, and it is not completely clear how to identify the suitable functional spaces for the  study of the limit. The Wigner transform, which allows to compute observables, is non-positive and quadratic with respect to the wave function~: this implies that the combinatorics associated with the Duhamel series, which can be represented by Feynman diagrams is much worse than the combinatorics of collision trees. In general these formal series are never absolutely convergent.

All the open questions regarding the long time behavior, the structure of correlations and the deviations from the average dynamics, the role of microscopic interactions or the stationary non equilibrium case remain, also in these different settings, as challenges for the future.


\bigskip

\noindent
{\it Acknowledgments.}
We thank P. Dario, C. Garban, E. Ghys,  F. Golse and J. Marklof for their very useful comments on a preliminary version of this manuscript.
This work was partially supported by ANR-15-CE40-0020-01 grant LSD.


 \bibliographystyle{abbrv}
 \bibliography{biblio_Boltzmann}

\begin{thebibliography}{10}

\bibitem{AMP21}
I.~Ampatzoglou, J.~K. Miller, and N.~Pavlovi{\'c}.
\newblock A rigorous derivation of a {B}oltzmann system for a mixture of
  hard-sphere gases, 2021.

\bibitem{Ayi}
N.~{Ayi}.
\newblock {From Newton's law to the linear Boltzmann equation without cut-off}.
\newblock {\em {Commun. Math. Phys.}}, 350(3):1219--1274, 2017.

\bibitem{Bardos_golse_levermore}
C.~{Bardos}, F.~{Golse}, and C.~D. {Levermore}.
\newblock Fluid dynamic limits of kinetic equations. ii: {C}onvergence proofs
  for the {B}oltzmann equation.
\newblock {\em Commun. Pure Appl. Math.}, 46(5):667--753, 1993.

\bibitem{BBBO}
G.~{Basile}, D.~{Benedetto}, L.~{Bertini}, and C.~{Orrieri}.
\newblock {Large deviations for Kac-like walks}.
\newblock {\em {J. Stat. Phys.}}, 184(1):27, 2021.
\newblock Id/No 10.

\bibitem{Bastea_Esposito_Lebowitz_Marra}
S.~{Bastea}, R.~{Esposito}, J.~L. {Lebowitz}, and R.~{Marra}.
\newblock {Binary fluids with long range segregating interaction. I: Derivation
  of kinetic and hydrodynamic equations}.
\newblock {\em {J. Stat. Phys.}}, 101(5-6):1087--1136, 2000.

\bibitem{Bertini_DeSole_Gabrielli}
L.~Bertini, A.~De~Sole, D.~Gabrielli, G.~Jona-Lasinio, and C.~Landim.
\newblock Macroscopic fluctuation theory.
\newblock {\em Rev. Mod. Phys.}, 87:593--636, Jun 2015.

\bibitem{Bobylev_Dukes_Illner}
A.~V. {Bobylev}, P.~{Dukes}, R.~{Illner}, and H.~D. jun. {Victory}.
\newblock {On Vlasov-Manev equations. I: Foundations, properties, and nonglobal
  existence}.
\newblock {\em {J. Stat. Phys.}}, 88(3-4):885--911, 1997.

\bibitem{Bobylev_Pulvirenti_Saffirio}
A.~V. {Bobylev}, M.~{Pulvirenti}, and C.~{Saffirio}.
\newblock From particle systems to the {L}andau equation: a consistency result.
\newblock {\em {Commun. Math. Phys.}}, 319(3):683--702, 2013.

\bibitem{BGSR_cras}
T.~{Bodineau}, I.~{Gallagher}, and L.~{Saint-Raymond}.
\newblock {D}e la dynamique des sph\`eres dures aux \'equations de
  {S}tokes-{F}ourier: une analyse \({L}^2\) de la limite de {Boltzmann-Grad}.
\newblock {\em {C. R., Math., Acad. Sci. Paris}}, 353(7):623--627, 2015.

\bibitem{BGSR1}
T.~Bodineau, I.~Gallagher, and L.~Saint-Raymond.
\newblock The {B}rownian motion as the limit of a deterministic system of
  hard-spheres.
\newblock {\em Inventiones mathematicae}, 203(2):493--553, 2016.

\bibitem{BGSR2}
T.~Bodineau, I.~Gallagher, and L.~Saint-Raymond.
\newblock From hard sphere dynamics to the {S}tokes--{F}ourier equations: An
  analysis of the {B}oltzmann--{G}rad limit.
\newblock {\em Annals of PDE}, 3(1):2, 2017.

\bibitem{BGSS}
T.~Bodineau, I.~Gallagher, L.~Saint-Raymond, and S.~Simonella.
\newblock One-sided convergence in the {B}oltzmann--{G}rad limit.
\newblock {\em Annales de la Facult{\'e} des sciences de Toulouse:
  Math{\'e}matiques}, 27(5):985--1022, 2018.

\bibitem{BGSS1}
T.~Bodineau, I.~Gallagher, L.~Saint-Raymond, and S.~Simonella.
\newblock Fluctuation theory in the {B}oltzmann--{G}rad limit.
\newblock {\em Journal of Statistical Physics}, 180(1):873--895, 2020.

\bibitem{BGSS2}
T.~Bodineau, I.~Gallagher, L.~Saint-Raymond, and S.~Simonella.
\newblock Statistical dynamics of a hard sphere gas: fluctuating {B}oltzmann
  equation and large deviations.
\newblock {\em preprint arXiv:2008.10403}, 2020.

\bibitem{BGSS3}
T.~Bodineau, I.~Gallagher, L.~Saint-Raymond, and S.~Simonella.
\newblock Long-time correlations for a hard-sphere gas at equilibrium.
\newblock {\em preprint arXiv:2012.03813, to appear in {CPAM}}, 2021.

\bibitem{BGSS4}
T.~Bodineau, I.~Gallagher, L.~Saint-Raymond, and S.~Simonella.
\newblock Long-time derivation at equilibrium of the fluctuating boltzmann
  equation.
\newblock {\em arXiv:2201.04514}, 2022.

\bibitem{boltzmann}
L.~{Boltzmann}.
\newblock {Weitere Studien \"uber das W\"armegleichgewicht unter
  Gasmolec\"ulen.}
\newblock {\em {Wien. Ber.}}, 66:275--370, 1872.

\bibitem{Bonetto_Lebowitz_ReyBellet}
F.~{Bonetto}, J.~L. {Lebowitz}, and L.~{Rey-Bellet}.
\newblock {F}ourier's law: a challenge to theorists.
\newblock In {\em {M}athematical physics 2000.}, pages 128--150. London:
  Imperial College Press, 2000.

\bibitem{bouchet}
F.~Bouchet.
\newblock Is the {B}oltzmann equation reversible? {A} large deviation
  perspective on the irreversibility paradox.
\newblock {\em J. Stat. Phys.}, 181(2):515--550, 2020.

\bibitem{Braun_Hepp}
W.~{Braun} and K.~{Hepp}.
\newblock {The Vlasov dynamics and its fluctuations in the \(1/N\) limit of
  interacting classical particles}.
\newblock {\em {Commun. Math. Phys.}}, 56:101--113, 1977.

\bibitem{Cercignani_Gerasimenko_Petrina}
C.~Cercignani, V.~I. Gerasimenko, and D.~Y. Petrina.
\newblock {\em Many-particle dynamics and kinetic equations}, volume 420 of
  {\em Mathematics and its Applications}.
\newblock Kluwer Academic Publishers Group, Dordrecht, 1997.
\newblock Translated from the Russian manuscript by K. Petrina and V. Gredzhuk.

\bibitem{Cercignani_Illner_Pulvirenti}
C.~Cercignani, R.~Illner, and M.~Pulvirenti.
\newblock {\em The mathematical theory of dilute gases}, volume 106 of {\em
  Applied Mathematical Sciences}.
\newblock Springer-Verlag, New York, 1994.

\bibitem{Dembo_Zeitouni}
A.~Dembo and O.~Zeitouni.
\newblock {\em Large deviations techniques and applications}, volume~38 of {\em
  Stochastic Modelling and Applied Probability}.
\newblock Springer-Verlag, Berlin, 2010.
\newblock Corrected reprint of the second (1998) edition.

\bibitem{denlinger}
R.~Denlinger.
\newblock The propagation of chaos for a rarefied gas of hard spheres in the
  whole space.
\newblock {\em Arch. Ration. Mech. Anal.}, 229(2):885--952, 2018.

\bibitem{Derrida}
B.~Derrida.
\newblock Microscopic versus macroscopic approaches to non-equilibrium systems.
\newblock 2011(01):P01030, jan 2011.

\bibitem{Desvillettes_Pulvirenti}
L.~{Desvillettes} and M.~{Pulvirenti}.
\newblock The linear {B}oltzmann equation for long-range forces: A derivation
  from particle systems.
\newblock {\em {Math. Models Methods Appl. Sci.}}, 9(8):1123--1145, 1999.

\bibitem{DSS}
L.~Desvillettes, C.~Saffirio, and S.~Simonella.
\newblock Collisions in a mean-field: kinetic limit for the {L}orentz gas with
  long-range forces.
\newblock {\em In preparation}, 2022.

\bibitem{DiPerna_Lions}
R.~J. {DiPerna} and P.~L. {Lions}.
\newblock {On the Cauchy problem for Boltzmann equations: Global existence and
  weak stability}.
\newblock {\em {Ann. Math. (2)}}, 130(2):321--366, 1989.

\bibitem{Dolmaire}
T.~Dolmaire.
\newblock About {L}anford's theorem in the half-space with specular reflection.
\newblock {\em arXiv preprint arXiv:2102.05513}, 2021.

\bibitem{duerinckx_saint-raymond}
M.~Duerinckx and L.~Saint-Raymond.
\newblock Lenard--{B}alescu correction to mean-field theory.
\newblock {\em Probability and Mathematical Physics}, 2(1):27--69, 2021.

\bibitem{DGL}
D.~{D\"urr}, S.~{Goldstein}, and J.~L. {Lebowitz}.
\newblock {A}symptotic motion of a classical particle in a random potential in
  two dimensions: {L}andau model.
\newblock {\em {Commun. Math. Phys.}}, 113:209--230, 1987.

\bibitem{erdos}
L.~{Erd\H{o}s}.
\newblock {\em {Lecture notes on quantum Brownian motion}}, pages 3--98.
\newblock Oxford University Press, 2012.

\bibitem{erdos-yau}
L.~{Erd\H{o}s}, M.~{Salmhofer}, and H.-T. {Yau}.
\newblock {Quantum diffusion of the random Schr\"odinger evolution in the
  scaling limit}.
\newblock {\em {Acta Math.}}, 200(2):211--277, 2008.

\bibitem{ernst1981nonequilibrium}
M.~Ernst and E.~Cohen.
\newblock Nonequilibrium fluctuations in $\mu$ space.
\newblock {\em Journal of Statistical Physics}, 25(1):153--180, 1981.

\bibitem{EspositoGuoMarra}
R.~Esposito, Y.~Guo, and R.~Marra.
\newblock Validity of the {B}oltzmann equation with an external force.
\newblock {\em Kinetic \& Related Models}, 4(2):499--515, 2011.

\bibitem{EspositoMarra_Stationary}
R.~{Esposito} and R.~{Marra}.
\newblock {Stationary non equilibrium states in kinetic theory}.
\newblock {\em {J. Stat. Phys.}}, 180(1-6):773--809, 2020.

\bibitem{GSRT}
I.~{Gallagher}, L.~{Saint-Raymond}, and B.~{Texier}.
\newblock {\em From {N}ewton to {B}oltzmann: hard spheres and short-range
  potentials}.
\newblock Z\"urich: European Mathematical Society (EMS), 2013.

\bibitem{Golse2016}
F.~Golse.
\newblock {\em On the Dynamics of Large Particle Systems in the Mean Field
  Limit}, pages 1--144.
\newblock Springer International Publishing, Cham, 2016.

\bibitem{Golse_Levermore_Saint-Raymond}
F.~Golse, C.~D. Levermore, and L.~Saint-Raymond.
\newblock La m\'{e}thode de l'entropie relative pour les limites
  hydrodynamiques de mod\`eles cin\'{e}tiques.
\newblock In {\em S\'{e}minaire: \'{E}quations aux {D}\'{e}riv\'{e}es
  {P}artielles, 1999--2000}, S\'{e}min. \'{E}qu. D\'{e}riv. Partielles, pages
  Exp. No. XIX, 23. \'{E}cole Polytech., Palaiseau, 2000.

\bibitem{Grad}
H.~Grad.
\newblock {\em Principles of the kinetic theory of gases}, pages 205--294.
\newblock Springer-Verlag, Berlin-G\"{o}ttingen-Heidelberg, 1958.

\bibitem{heydecker2021large}
D.~Heydecker.
\newblock Large deviations of {K}ac's conservative particle system and energy
  non-conserving solutions to the {B}oltzmann equation: A counterexample to the
  predicted rate function.
\newblock {\em preprint arXiv:2103.14550}, 2021.

\bibitem{Jabin_Wang2016}
P.-E. Jabin and Z.~Wang.
\newblock Mean field limit and propagation of chaos for {V}lasov systems with
  bounded forces.
\newblock {\em Journal of Functional Analysis}, 271(12):3588--3627, 2016.

\bibitem{Kac_kinetic}
M.~Kac.
\newblock {\em Foundations of kinetic theory}, pages 171--197.
\newblock University of California Press, Berkeley and Los Angeles, 1956.

\bibitem{KP}
H.~{Kesten} and G.~C. {Papanicolaou}.
\newblock {A limit theorem for stochastic acceleration}.
\newblock {\em {Commun. Math. Phys.}}, 78:19--63, 1980.

\bibitem{King}
F.~G. King.
\newblock {\em {BBGKY} hierarchy for positive potentials}.
\newblock University of California, Berkeley, 1975.

\bibitem{Lanford}
O.~E. Lanford, III.
\newblock {\em Time evolution of large classical systems}, pages 1--111.
  Lecture Notes in Phys., Vol. 38.
\newblock 1975.

\bibitem{LeBihan}
C.~Le~Bihan.
\newblock Boltzmann-{G}rad limit of a hard sphere system in a box with
  diffusive boundary conditions.
\newblock {\em arXiv preprint arXiv:2104.04354, to appear in {Disc. Cont. Dyn
  Syst.}}, 2021.

\bibitem{Leonard}
C.~{L\'eonard}.
\newblock {On large deviations for particle systems associated with spatially
  homogeneous Boltzmann type equations}.
\newblock {\em {Probab. Theory Relat. Fields}}, 101(1):1--44, 1995.

\bibitem{logan1976fluctuations}
J.~Logan and M.~Kac.
\newblock Fluctuations and the {B}oltzmann equation.
\newblock {\em Physical Review A}, 13(1):458, 1976.

\bibitem{Meleard_fluctuations}
S.~M\'el\'eard.
\newblock Convergence of the fluctuations for interacting diffusions with jumps
  associated with {B}oltzmann equations.
\newblock {\em Stochastics Stochastics Rep.}, 63(3-4):195--225, 1998.

\bibitem{Mischler_Mouhot}
S.~{Mischler} and C.~{Mouhot}.
\newblock {K}ac's program in kinetic theory.
\newblock {\em {Invent. Math.}}, 193(1):1--147, 2013.

\bibitem{NVW2}
A.~Nota, J.~J.~L. Vel{\'a}zquez, and R.~Winter.
\newblock Interacting particle systems with long range interactions:
  approximation by tagged particles in random fields, 2021.

\bibitem{Nota_Velazquez_Winter}
A.~Nota, J.~J.~L. Vel\'{a}zquez, and R.~Winter.
\newblock Interacting particle systems with long-range interactions: scaling
  limits and kinetic equations.
\newblock {\em Atti Accad. Naz. Lincei Rend. Lincei Mat. Appl.},
  32(2):335--377, 2021.

\bibitem{penrose1967convergence}
O.~Penrose.
\newblock Convergence of fugacity expansions for classical systems.
\newblock page 101. Benjamin, New York, 1967.

\bibitem{poghosyan2009abstract}
S.~Poghosyan and D.~Ueltschi.
\newblock Abstract cluster expansion with applications to statistical
  mechanical systems.
\newblock {\em Journal of mathematical physics}, 50(5):053509, 2009.

\bibitem{Pu06}
M.~{Pulvirenti}.
\newblock {The weak-coupling limit of large classical and quantum systems}.
\newblock In {\em Proceedings of the international congress of mathematicians
  {(ICM)}, {Madrid, Spain}}, pages 229--256. Z\"urich: European Mathematical
  Society (EMS), 2006.

\bibitem{Pulvirenti_Saffirio_Simonella}
M.~Pulvirenti, C.~Saffirio, and S.~Simonella.
\newblock On the validity of the {B}oltzmann equation for short range
  potentials.
\newblock {\em Rev. Math. Phys.}, 26(2):1450001, 64, 2014.

\bibitem{Pulvirenti_Simonella}
M.~Pulvirenti and S.~Simonella.
\newblock The {B}oltzmann-{G}rad limit of a hard sphere system: analysis of the
  correlation error.
\newblock {\em Invent. Math.}, 207(3):1135--1237, 2017.

\bibitem{Rezakhanlou_fluctuations}
F.~Rezakhanlou.
\newblock Equilibrium fluctuations for the discrete {B}oltzmann equation.
\newblock {\em Duke Math. J.}, 93(2):257--288, 1998.

\bibitem{Rezakhanlou_GD}
F.~Rezakhanlou.
\newblock Large deviations from a kinetic limit.
\newblock {\em Ann. Probab.}, 26(3):1259--1340, 1998.

\bibitem{Rezakhanlou_BG}
F.~Rezakhanlou.
\newblock Boltzmann-{G}rad limits for stochastic hard sphere models.
\newblock {\em Comm. Math. Phys.}, 248(3):553--637, 2004.

\bibitem{RezakhanlouLNM}
F.~Rezakhanlou.
\newblock {Kinetic limits for interacting particle systems}.
\newblock In {\em Entropy methods for the Boltzmann equation. Lectures from a
  special semester at the Centre \'Emil Borel, Instititut H. Poincar\'e, Paris
  2001.}, pages 71--105. Berlin: Springer, 2008.

\bibitem{Ruelle_livre}
D.~Ruelle.
\newblock {\em Statistical mechanics}.
\newblock World Scientific Publishing Co., Inc., River Edge, NJ; Imperial
  College Press, London, 1999.
\newblock Rigorous results, Reprint of the 1989 edition.

\bibitem{SR}
L.~{Saint-Raymond}.
\newblock {\em {Hydrodynamic limits of the Boltzmann equation}}, volume 1971.
\newblock Berlin: Springer, 2009.

\bibitem{Spohn_fluctuations}
H.~Spohn.
\newblock Fluctuations around the {B}oltzmann equation.
\newblock {\em J. Statist. Phys.}, 26(2):285--305, 1981.

\bibitem{Spohn_review_fluctuations}
H.~Spohn.
\newblock Fluctuation theory for the {B}oltzmann equation.
\newblock In {\em Nonequilibrium phenomena, {I}}, volume~10 of {\em Stud.
  Statist. Mech.}, pages 225--251. North-Holland, Amsterdam, 1983.

\bibitem{spohn2012large}
H.~Spohn.
\newblock {\em Large scale dynamics of interacting particles}.
\newblock Springer Science \& Business Media, 2012.

\bibitem{vanbeijeren_lanford_lebowitz_spohn}
H.~van Beijeren, O.~E. Lanford, III, J.~L. Lebowitz, and H.~Spohn.
\newblock Equilibrium time correlation functions in the low-density limit.
\newblock {\em J. Statist. Phys.}, 22(2):237--257, 1980.

\bibitem{Velazquez_Winter}
J.~J.~L. {Vel\'azquez} and R.~{Winter}.
\newblock From a non-{M}arkovian system to the {L}andau equation.
\newblock {\em {Commun. Math. Phys.}}, 361(1):239--287, 2018.

\bibitem{Winter}
R.~{Winter}.
\newblock Convergence to the {L}andau equation from the truncated {BBGKY}
  hierarchy in the weak-coupling limit.
\newblock {\em {J. Differ. Equations}}, 283:1--36, 2021.

\bibitem{Yau}
H.-T. {Yau}.
\newblock {Relative entropy and hydrodynamics of Ginzburg-Landau models}.
\newblock {\em {Lett. Math. Phys.}}, 22(1):63--80, 1991.

\end{thebibliography}

%
%
%
%
%
%
%
%
%
%
%
%
%
%
%
%
%
%
%
%

\end{document}